\documentclass[12pt]{amsart}

\usepackage{amsmath, amssymb, amsthm}

\usepackage{graphicx}

\usepackage{hyperref}
\usepackage{geometry}
\usepackage{tikz}
\usetikzlibrary{shapes.misc}
\usepackage{caption}
\usepackage{booktabs}
\usepackage{enumitem}
\usepackage{float}
\usepackage{longtable}
\newtheorem{Def}{Definition}[section]

\newtheorem{Prop}[Def]{Proposition}
\newtheorem{Lem}[Def]{Lemma}

\newtheorem{Rem}[Def]{Remark}
\newtheorem{Pf}{Proof}
\newtheorem*{mainthm}{Main Theorem}

\newcommand{\R}{\mathbb{R}}

\newcommand{\Q}{\mathbb{Q}}
\newcommand{\C}{\mathbb{C}}
\newcommand{\Z}{\mathbb{Z}}

\newcommand{\rP}{\mathrm{P}}
\newcommand{\mr}[1]{\mathrm{#1}}
\newcommand{\wt}[1]{\widetilde{#1}}
\newcommand{\al}{\alpha}
\newcommand{\uQ}{\underline{\Q}}

\allowdisplaybreaks

\title{BOUNDARY COHOMOLOGY OF $\mr{Sp}_6(\Z)$: \\
TRIVIAL REPRESENTATION}
\author{RYUTO MITOMA}
\address{Joint Graduate School of Mathematics for Innovation, Kyushu University}
\email{mitoma.ryuto.491@s.kyushu-u.ac.jp}

\date{\today}

\begin{document}

\begin{abstract}
In this article, we compute the boundary cohomology of the arithmetic group $\mr{Sp}_6(\Z)$ with coefficients in the trivial representation. 
Our computation utilizes the Borel-Serre compactification and the associated spectral sequence.
\end{abstract}

\maketitle

\tableofcontents

 
 

\section{Introduction}

The cohomology of arithmetic groups plays a central role in modern number theory, particularly in the context of the Langlands program. 
In this article, we compute the boundary cohomology of the arithmetic group $\Gamma = \mr{Sp}_{6}(\Z)$ with coefficients in the trivial representation $\Q$. 
Our computation utilizes the Borel-Serre compactification $\overline{S_\Gamma}$ of the locally symmetric space $S_\Gamma$ associated with $\mr{Sp}_{6}$.

The boundary $\partial S_\Gamma$ of the Borel-Serre compactification admits a stratification by faces $\partial_{\rP, \Gamma}$ corresponding to the conjugacy classes of $\Q$-parabolic subgroups $\rP$. 
This geometric structure yields a spectral sequence converging to the cohomology of the boundary:
\[
  E_1^{p, q} = \bigoplus_{prk(P) = p+1} \mr{H}^q(\partial_{P, \Gamma}, \wt{\mathcal{M}_\lambda}) \Longrightarrow \mr{H}^{p+q}(\partial S_\Gamma, \wt{\mathcal{M}_\lambda}).
\]
In our case, since the $\Q$-split rank of $\mr{Sp}_{6}$ is $3$, the spectral sequence has non-trivial terms only for columns $p=0, 1, 2$, and it degenerates at the $E_{3}$-page.
We explicitly compute the $E_1$-terms using the cohomology of the faces, determine the differentials $d_1$ and $d_2$, and thereby obtain the full boundary cohomology.

The main result of this paper is summarized as follows:

\begin{mainthm}
The boundary cohomology groups $\mr{H}^{q}(\partial S_{\Gamma}, \Q)$ of $\mr{Sp}_{6}(\Z)$ with trivial coefficients are given by:
\begin{equation*}
  \mr{H}^{q}(\partial S_{\Gamma}, \uQ) = 
  \begin{cases}
    \Q & q = 0,2,4,7,9,11 \\
    \Q^2 & q = 5,6 \\
    0 & \text{otherwise}
  \end{cases}
\end{equation*}
\end{mainthm}

In Section 2, we review the construction of the Borel-Serre compactification and the associated spectral sequence. We also define the specific differentials $d_1$ (horizontal) and $d_v$ (vertical) arising from the double complex structure.
In Section 3, we analyze the cohomology of the faces associated with parabolic subgroups of various ranks, focusing on the parity conditions.
Finally, in Section 4, we carry out the explicit computation of the spectral sequence from the $E_1$-page to the $E_3$-page to prove the Main Theorem.

 
 
\section{Basic Notions}

\subsection{Structure theory}

\

In this subsection, we review the basic properties of $\mr{Sp}_6$ and fix the notation.
The symplectic group $\mr{Sp}_{6}(K)$ over a field $K$ is defined by
\[
  \mr{Sp}_{6}(K) = \{M \in \mr{GL}_{6}(K) \mid M^t J M = J \}
\]
where $M^t$ denotes the transpose of $M$, and
\[
  J =
  \begin{pmatrix}
    0 & I_3 \\
    -I_3 & 0
  \end{pmatrix} 
\]
where $I_3$ is the $3 \times 3$ identity matrix.
The unitary group $\mr{U}(3)$ is identified with the maximal compact subgroup of $\mr{Sp}_{6}(\R)$; 
we denote this maximal compact subgroup by $K_{\infty}$.
\[
  K_{\infty} = \left\{ 
    \begin{pmatrix}
      A & B \\
      -B & A
    \end{pmatrix}
    \mid A + iB \in \mr{U}_{3}
   \right\}
  \]
The quotient $\mr{Sp}_{6}(\R)/K_{\infty}$ is identified with the Siegel upper half-space $\mathcal{H}_{3}$
\[
  \mathcal{H}_{3} = \{Z \in \mr{M}_3(\C) \mid Z^t = Z, \mr{Im}(Z) > 0 \},
\]
where $\mr{Im}(Z) > 0$ means that the imaginary part of $Z$ is positive definite.
The group $\mr{Sp}_{6}(\R)$ acts on $\mathcal{H}_{3}$ by
\[
  g \cdot Z = (AZ+B)(CZ+D)^{-1}, \quad \text{for }
  g =
  \begin{pmatrix}
    A & B \\
    C & D
  \end{pmatrix}
  \in \mr{Sp}_{6}(\R) 
\]
with $A, B, C, D \in \text{M}_3(\R)$, and $Z \in \mathcal{H}_{3}$.
Let $\Gamma = \mr{Sp}_{6}(\Z)$ be the arithmetic subgroup. 
The quotient space $S_{\Gamma} = \Gamma \backslash \mathcal{H}_{3}$ is called the Siegel modular variety of degree $3$.

Let $\mr{T}$ be the maximal torus of $\mr{Sp}_6$ consisting of diagonal matrices:
\begin{align*}
  \mr{T} = \{\mr{diag}(t_1, t_2, t_3, t_1^{-1}, t_2^{-1}, t_{3}^{-1})
   \mid t_{i} \in \R^{\times} \}
\end{align*}
Let $\varepsilon_{i} \in \mr{X}^*(\mr{T})$ be the character such that $\varepsilon_i(\mr{diag}(t_1, t_2, t_3, t_1^{-1}, t_2^{-1}, t_{3}^{-1})) = t_{i}$. 
The root system $\Phi$ of type $\mr{C}_3$ is then described as:
\begin{align*}
  \Phi = \{ \pm \varepsilon_{i} \pm \varepsilon_{j} \mid 1 \leq i < j \leq 3 \} \cup \{\pm2 \varepsilon_{i} \mid 1 \leq i \leq 3 \}.
\end{align*}
We fix the set of positive roots $\Phi^{+}$ and simple roots $\pi$ as:
\begin{align*}
  \Phi^+ &= \{\varepsilon_1 - \varepsilon_2, \varepsilon_1 - \varepsilon_3, \varepsilon_2 - \varepsilon_3, \varepsilon_1 + \varepsilon_2, \varepsilon_1 + \varepsilon_3, \varepsilon_2 + \varepsilon_3, 2\varepsilon_1, 2\varepsilon_2, 2\varepsilon_3 \} \\
  \pi &= \{ \alpha_1 = \varepsilon_1-\varepsilon_2, \alpha_2 = \varepsilon_2-\varepsilon_3, \alpha_3 = 2\varepsilon_3 \}
\end{align*}
The fundamental dominant weights dual to these simple roots are
\[
  \{ \gamma_1 = \varepsilon_1, \quad
   \gamma_2 = \varepsilon_1 + \varepsilon_2, \quad 
   \gamma_3 = \varepsilon_1 + \varepsilon_2 + \varepsilon_3 \}
\]
The irreducible finite-dimensional representations of $\mr{Sp}_6$ are determined by their highest weights 
$\lambda = m_1\gamma_1 + m_2\gamma_2 + m_3\gamma_3$, where $m_i \in \Z_{\ge 0}$ are non-negative integers.
In this article, we focus exclusively on the trivial representation, $\mr{i.e.}$, the case where $\lambda = 0$. 
The Weyl group $\mathcal{W}$ is isomorphic to $(\Z/2\Z)^3 \rtimes \mathfrak{S}_3$.

There is a one-to-one correspondence between 
the set of proper standard $\Q$-parabolic subgroups and
the set of non-empty subsets of simple roots $\pi = \{\alpha_1, \alpha_2, \alpha_3\}$.
We denote the standard parabolic subgroup corresponding to a subset $I \subset \pi$ by $P_{I}$.

In this article, 
we adopt the convention that the cardinality $|I|$ equals the parabolic rank of $P_{I}$. 
Under this convention, the maximal parabolic subgroups correspond to the singleton sets $\{\alpha_1\}, \{\alpha_2\}$, and $\{\alpha_3\}$, 
while the Borel subgroup corresponds to the full set $\pi$.

We obtain the structure of the Levi quotient $\mr{M}_{P_I}$ from the Dynkin diagram of type $\mr{C}_3$. 
Specifically, $\mr{M}_{P_{I}}$ is, up to isogeny, the product of a semisimple group and a torus of dimension $|I|$. 
The semisimple part corresponds to the sub-diagram obtained by keeping the nodes in the complement $\pi \setminus I$. 
The correspondence between the subset $I$ and the Levi quotient $\mr{M}_{P_I}$ is summarized in Table \ref{table:levi_quotients}. 
For a detailed diagrammatic correspondence, see Appendix \ref{str of levi}.

\begin{table}[H]
  \centering
  \caption{Standard $\Q$-parabolic subgroups and Levi quotients}
  \label{table:levi_quotients}
  \begin{tabular}{lcl}
    \toprule
    Rank & Subset $I$ & Levi quotient $\mr{M}_{P_I}$ \\
    \midrule
    1 & $\{\alpha_1\}$ & $\mr{GL}_1 \times \mr{Sp}_4$ \\
           & $\{\alpha_2\}$ & $\mr{GL}_2 \times \mr{Sp}_2$ \\
           & $\{\alpha_3\}$ & $\mr{GL}_3$ \\
    \midrule
    2 & $\{\alpha_1, \alpha_2\}$ & $\mr{GL}_1 \times \mr{GL}_1 \times \mr{Sp}_2$ \\
           & $\{\alpha_1, \alpha_3\}$ & $\mr{GL}_1 \times \mr{GL}_2$ \\
           & $\{\alpha_2, \alpha_3\}$ & $\mr{GL}_2 \times \mr{GL}_1$ \\
    \midrule
    3 & $\pi$ & $\mr{GL}_1 \times \mr{GL}_1 \times \mr{GL}_1$ \\
    \bottomrule
  \end{tabular}
\end{table}

\subsection{Borel-Serre compactification and spectral sequence}
 
In this subsection, 
we describe the spectral sequence arising from the Borel-Serre compactification of the locally symmetric space. 
For a split semisimple group $\mr{G}$ over $\Q$, the maximal compact subgroup $\mr{K}_{\infty}$ of $\mr{G}(\R)$ and an arithmetic subgroup $\Gamma$, 
the corresponding locally symmetric space is 
\[
  S_{\Gamma} = \Gamma \backslash \mr{G}(\R) / \mr{K}_{\infty}.
\]
We consider the Borel-Serre compactification $\overline{S_\Gamma}$ of $S_\Gamma$ (\cite{BoSe}), 
whose boundary $\partial S_\Gamma = \overline{S_\Gamma} \backslash S_\Gamma$ is a union of spaces indexed by the $\Gamma$-conjugacy classes of $\Q$-parabolic subgroups of $\mr{G}$.
\[
  \partial S_\Gamma = \bigcup _{\rP \in \mathcal{P}_{\Q}(\mr{G})} \partial_{\rP, \Gamma}.
\]
Here $\mathcal{P}_{\Q}(\mr{G})$ denotes the set of the standard $\Q$-parabolic subgroups determined by the choice of a maximal split torus $\mr{T}$ of $\mr{G}$ and a system of positive roots $\Phi^+$.

Let $\mathcal{M}_\lambda$ be the irreducible representation of $\mr{G}$ with highest weight $\lambda$. 
This representation defines a sheaf $\wt{\mathcal{M}_{\lambda}}$ over $S_\Gamma$, 
which is defined over $\Q$ as follows:
\[
  \wt{\mathcal{M}_\lambda}(U) =
   \left\{
    f : \pi^{-1}(U) \to \mathcal{M}_\lambda \middle |
    \begin{array}{l}
     f \  \text{is locally constant,} \\
     f(\gamma u) = \gamma f(u)\  \text{for any}\  \gamma \in \Gamma, u \in \pi^{-1}(U)
    \end{array}  
   \right\}
\]
where $U$ is an open subset of $S_\Gamma$, 
and $\pi : \mr{G}(\R)/\mr{K}_\infty \to \Gamma \backslash \mr{G}(\R)/\mr{K}_\infty = S_\Gamma$ is the projection.
In the case of the trivial representation, 
$\mathcal{M}_0 = \Q$, 
the associated sheaf $\wt{{\mathcal{M}_0}}$ is canonically isomorphic to the constant sheaf $\uQ$.

By applying the direct image functor associated to the inclusion $i:S_\Gamma \hookrightarrow \overline{S}_\Gamma$, 
we obtain a sheaf on $\overline{S_{\Gamma}}$.
Since this inclusion is a homotopy equivalence (\cite{BoSe}), 
it induces an isomorphism
\[ 
\mr{H}^\bullet(S_\Gamma, \wt{\mathcal{M}_\lambda}) = \mr{H}^\bullet(\overline{S_\Gamma}, i_{\ast}(\wt{\mathcal{M}_\lambda})). 
\]
For simplicity, 
we denote the direct image sheaf $i_{\ast}(\wt{\mathcal{M}_\lambda})$ 
and its restriction to the boundary and its faces $\partial_{\rP, \Gamma}$
by the same symbol $\wt{\mathcal{M}_\lambda}$.

The stratification of the boundary yields a spectral sequence converging to the boundary cohomology:
\[
  E_1^{p, q} = \bigoplus_{prk(P) = p+1} \mr{H}^q(\partial_{P, \Gamma}, \wt{\mathcal{M}_\lambda}) \Longrightarrow \mr{H}^{p+q}(\partial S_\Gamma, \wt{\mathcal{M}_\lambda}).
\]
where $prk(P)$ denotes the parabolic rank of $P$.
In our convention, 
this rank coincides with the number of elements in the corresponding subset $I \subset \pi$,
that is, $prk(P_I) = |I|$.

This spectral sequence is derived from the double complex 
\[
  C^{p,q} = \bigoplus_{prk(P) = p+1} C^{q}(\partial_{P, \Gamma}, \wt{\mathcal{M}_\lambda}).
\]
We define two differentials for this double complex.
The vertical differential is the direct sum of the cochain differentials on each face:
\begin{align*}
  d_{v}^{p,q} = \bigoplus_{prk(P) = p+1} d_{P}^{q} : C^{p,q} \to C^{p,q+1},
\end{align*}
where $d_{P}^{q} : C^{q}(\partial_{P, \Gamma}, \wt{\mathcal{M}}_\lambda) \to C^{q+1}(\partial_{P, \Gamma}, \wt{\mathcal{M}}_\lambda)$ 
denotes the standard differential of the cochain complex for each face $\partial_{P, \Gamma}$.
The horizontal differential is defined by the sum of restriction maps induced by the inclusions of faces $\partial_{\mr{Q}, \Gamma} \hookrightarrow \partial_{P, \Gamma}$ 
for pairs of parabolic subgroups $\mr{Q} \subset P$ with $prk(\mr{Q}) = prk(P) + 1$:
\begin{align*}
  d_{h}^{p,q}:= \sum_{\substack{P:prk(P)=p+1, \\ \mr{Q}\subset P: prk(\mr{Q})=p+2}} \epsilon(P, \mr{Q})i_{\mr{Q}, P}^\bullet 
\end{align*}
where $i_{\mr{Q}, P}^\bullet : C^q(\partial_{P, \Gamma}) \to C^q(\partial_{\mr{Q}, \Gamma})$ is the restriction map and
 $\epsilon(P, \mr{Q})$ is a sign depending on the relative positions of $P$ and $\mr{Q}$.
These differentials satisfy $d_h^{p+1,q} \circ d_h^{p, q} = 0$ and $d_h^{p,q+1} \circ d_v^{p,q} = d_v^{p+1,q} \circ d_h^{p,q}$.
These differentials induce a spectral sequence converging to the boundary cohomology $\mr{H}^{\bullet}(\partial S_\Gamma, \wt{\mathcal{M}_\lambda})$.

The $E_1$-page is given by the cohomology of the vertical complexes:
\[
  E_{1}^{p,q} = \mr{H}^q(C^{p, \ast}, d_v).
\] 
The differential $d_1 : E_{1}^{p,q} \to E_{1}^{p+1,q}$ is the map induced on the cohomology by $d_h$.
The $E_{2}$-page is defined as the cohomology of the $E_1$-page with respect to $d_{1}$. 
The next differential $d_2$ is obtained as follows:
An element of $E_{2}^{p,q}$ can be identified with a pair $(a, b) \in C^{p,q} \times C^{p+1, q-1}$ 
satisfying $d_v(a) = 0$ and $d_h(a) + d_v(b) = 0$.
In this representation,
two pairs $(a,b)$ and $(a', b')$ are equivalent if their difference lies in 
the subgroup generated by the following types of pairs:
\begin{enumerate}
  \item $(d_v(x), d_h(x))$ for some $x \in C^{p, q-1}$
  \item $(0,d_v(y))$ for some $y \in C^{p+1,q-2}$
  \item $(d_h(c), 0)$ for some $c \in C^{p-1, q}$ such that $d_v(c) = 0$
\end{enumerate}
The first and second types represent the vanishing of elements at the $E_1$-level,
while the third type represents the boundaries of the $d_1$.

The mapping $(a, b) \mapsto (d_h(b), 0)$ determines a well-defined differential 
\[
  d_2:E_{2}^{p,q} \to E_{2}^{p+2, q-1}.
\]
The $E_3$-page is obtained as the cohomology of the $E_2$-page with respect to the $d_2$.
Since the $\Q$-split rank of $\mr{Sp}_6$ is $3$, 
the spectral sequence degenerates at this page,
and the boundary cohomology is determined as the direct sum of the $E_3$ terms.
\[
  \mr{H}^{n}(\partial \mr{S}_{\Gamma}, \wt{\mathcal{M}_{\lambda}}) = \bigoplus_{p+q = n} E_{3}^{p,q}.
\]

To compute the $E_1$-terms, 
we utilize the relationship with Lie algebra cohomology.
For a standard $\Q$-parabolic subgroup $P$,
let $\mr{U}_{P}$ be its unipotent radical and $\mr{M}_{P} = P/\mr{U}_{P}$ be its Levi quotient.
Let $\mathfrak{u}_{P} = \mr{Lie}(\mr{U}_{P})$ be the Lie algebra of $\mr{U}_{P}$.
We denote the images of $\Gamma \cap P(\Q)$ and $K_{\infty} \cap P(\R)$ under the canonical projection 
$P \to \mr{M}_{P}$ by $\Gamma_{\mr{M}_{P}}$ and $K_{\infty}^{\mr{M}_{P}}$, respectively.
To define the corresponding locally symmetric space, we consider the following subgroup:
\[
  ^\circ \mr{M} = \bigcap_{\chi \in X^{\ast}_{\Q}(\mr{M})}(\mr{Ker} \chi^2)
\]
where $X^{\ast}_{\Q}(\mr{M}_{P})$ denotes the set of $\Q$-characters of $\mr{M}_{P}$.
The locally symmetric space associated with the Levi quotient $\mr{M}_{P}$ is then defined by:
\[
  S_{\Gamma}^{\mr{M}_{P}} = \Gamma_{\mr{M}_{P}} \backslash ^\circ \mr{M}_{P}(\R) / K_{\infty}^{\mr{M}_{P}}.
\]
Let $\mathcal{W}^{P}$ be the set of Kostant representatives for the parabolic subgroup $P$, defined by
\[
  \mathcal{W}^{P} = \{ w \in \mathcal{W} \mid w(\Phi^-) \cap \Phi^{+} \subset \Phi^{+}(\mathfrak{u}_{P}) \}
\]
where $\Phi^{+}(\mathfrak{u}_{P})$ denotes the set of roots whose root spaces are contained in $\mathfrak{u}_{P}$.
We denote the half of the sum of the positive roots by $\rho$:
\[
  \rho = \frac{1}{2} \sum_{\alpha \in \Phi^{+}} \alpha.
\]
In the case of $\mr{Sp}_6$, we have $\rho = 3\varepsilon_1 + 2\varepsilon_2 + \varepsilon_3 = \gamma_1 + \gamma_2 + \gamma_3$.
For each $w \in \mathcal{W}^{P}$, the element
\[
  w \cdot \lambda = w(\lambda + \rho) - \rho
\]
defines a highest weight of an irreducible representation $\mathcal{M}_{w \cdot \lambda}$ of $^\circ \mr{M}_{P}$.
As in the case of $\mr{G}$, 
this representation induces a local system $\wt{\mathcal{M}_{w \cdot \lambda}}$ over the locally symmetric space $S^{\mr{M}_{P}}_{\Gamma}$.
These definitions allow us to relate the cohomology of each face to the cohomology of the locally symmetric spaces associated with the Levi quotients.

The cohomology of the face $\partial_{P, \Gamma}$ is computed by a spectral sequence whose $E_2$-page is given by:
\[
  E_2^{i,j} = \mr{H}^i(S^{\mr{M}_{P}}, \wt{\mr{H}^j(\mathfrak{u}_{P}, \mathcal{M}_{\lambda})}) \Longrightarrow \mr{H}^{i+j}(\partial_{P, \Gamma}, \wt{\mathcal{M}_{\lambda}}).
\]
In the case of $\mr{Sp}_6$, this spectral sequence degenerates at the $E_2$-page (cf. \cite{gsp2}). Thus, we obtain the following isomorphism of vector spaces:
\[
  \mr{H}^{q}(\partial_{P, \Gamma}, \wt{\mathcal{M}_{\lambda}}) \cong \bigoplus_{i+j=q} \mr{H}^i(S^{\mr{M}_{P}}, \wt{\mr{H}^j(\mathfrak{u}_{P}, \mathcal{M}_{\lambda})}).
\]
By applying the Kostant theorem 
\[
  \mr{H}^j(\mathfrak{u}_{P}, \mathcal{M}_{\lambda}) \cong \bigoplus_{w\in \mathcal{W}^{P}, l(w)=j} \mathcal{M}_{w \cdot \lambda},
\]
we obtain the explicit decomposition for each face:
\[
  \mr{H}^{q}(\partial_{P, \Gamma}, \wt{\mathcal{M}_{\lambda}}) \cong \bigoplus_{w\in \mathcal{W}^{P}}\mr{H}^{q-l(w)}(S^{\mr{M}_{P}}, \wt{\mathcal{M}}_{w \cdot \lambda}).
\]
By combining these results and substituting them into the definition of the spectral sequence, we obtain the following explicit formula for the $E_1$-terms:
\[
  E_1^{p,q} = \bigoplus_{prk(P)=p+1} \left( \bigoplus_{w\in \mathcal{W}^{P}}\mr{H}^{q-l(w)}(S^{\mr{M}_{P}}, \wt{\mathcal{M}_{w \cdot \lambda}}) \right).
\]

The differential $d_1^{p,q}: E_1^{p,q} \to E_1^{p+1, q}$ 
is induced by the horizontal differential $d_h$.
It is composed of the restriction maps between the faces as follows:
\[
  d_1^{p,q} = \bigoplus_{\substack{prk(\mr{Q})=p+1 \\ prk(P)=p+2 \\ \mr{Q} \subset P}} \epsilon(\mr{Q}, P) r_{\mr{Q}, P}^{p,q}, \quad 
  r_{\mr{Q}, P}^{p,q} = \bigoplus_{\substack{w \in \mathcal{W}^{\mr{Q}} \\ s \in \mathcal{W}_{\mr{M}_\mr{Q}} / \mathcal{W}_{\mr{M}_P} \\ sw \in \mathcal{W}^{P}}} r_{\mr{Q}, P}^{p,q}(w,s).
\]
where 
\[
r_{\mr{Q}, P}^{p,q}(w,s) : \mr{H}^{q-l(w)}(\mr{S}^{\mr{M}_{\mr{Q}}}, \wt{\mathcal{M}_{w \cdot \lambda}}) \to \mr{H}^{q-l(sw)}(\mr{S}^{\mr{M}_{P}}, \wt{\mathcal{M}_{sw \cdot \lambda}})
\]
and $\epsilon(\mr{Q}, P)$ is the sign defined as follows:
The sign $\epsilon(\mr{Q}, P)$ is determined by the relative position of the simple roots defining the parabolic subgroups. 
Let $I(\mr{Q}) = \{ \alpha_{i_1}, \dots, \alpha_{i_n} \}$ with $i_1 < \dots < i_n$. 
When $I(P)$ is obtained by adding a simple root $\alpha_k$ to $I(\mr{Q})$ such that 
\[
  I(P) = \{ \alpha_{i_1}, \dots, \alpha_{i_{j-1}}, \alpha_k, \alpha_{i_{j}}, \dots, \alpha_{i_n} \} \quad \text{with} \quad i_{j-1} < k < i_{j}
\]
we define the sign as $\epsilon(\mr{Q}, P) = (-1)^j$. 
This convention ensures the relation $d_1^2 = 0$ in the spectral sequence. 

For example, 
consider the case where $I(\mr{Q}) = \{ \alpha_2 \}$, 
which corresponds to a maximal parabolic subgroup. 
If we add a root to obtain a rank 2 parabolic subgroup, 
the signs are determined as:
\begin{itemize}
  \item For $I(P_1) = \{ \alpha_1, \alpha_2 \}$, 
  the added root $\alpha_1$ is at the first position, 
  thus $\epsilon(\mr{Q}, P_1) = (-1)^1 = -1$.
  \item For $I(P_2) = \{ \alpha_2, \alpha_3 \}$, 
  the added root $\alpha_3$ is at the second position, 
  thus $\epsilon(\mr{Q}, P_2) = (-1)^2 = 1$.
\end{itemize}

The second differential $d_{2}^{p,q} : E_{2}^{p,q} \to E_{2}^{p+2, q-1}$ is defined through a zig-zag process on the $E_1$-page. 
This map is composed of restriction maps and their lifts between the faces of three distinct parabolic ranks. 
Following the notation established for $d_1$, 
the differential $d_2$ is expressed as follows:
\[
  d_2^{p,q} = \bigoplus_{\substack{prk(\mr{Q})=p+1 \\ prk(P)=p+3 \\ \mr{Q} \subset P}} r_{\mr{Q}, P}^{p,q}, \quad 
  r_{\mr{Q}, P}^{p,q} = \bigoplus_{\substack{w \in \mathcal{W}^{\mr{Q}} \\ s \in \mathcal{W}_{\mr{M}_{\mr{Q}}} / \mathcal{W}_{\mr{M}_{P}} \\ sw \in \mathcal{W}^{P}}} r_{\mr{Q}, P}^{p,q}(w,s).
\]
The individual map $r_{\mr{Q}, P}^{p,q}(w,s)$ is
\[
r_{\mr{Q}, P}^{p,q}(w,s) \colon \mr{H}^{q-l(w)}(\mr{S}^{\mr{M}_{\mr{Q}}}, \wt{\mathcal{M}}_{w \cdot \lambda}) \to \mr{H}^{q-1-l(sw)}(\mr{S}^{\mr{M}_{P}}, \wt{\mathcal{M}}_{sw \cdot \lambda}).
\]

\subsection{Kostant representatives}

In this subsection, 
we identify the set of Kostant representatives $\mathcal{W}^{P}$ for each standard parabolic subgroup $P$. 
Throughout this paper, 
we denote the product of simple reflections $s_i s_j \cdots s_k$ by the abbreviated form $s_{ij\cdots k}$. 
For example, 
$s_{12}$ represents the element $s_1 s_2 \in \mathcal{W}$. 
To characterize these representatives,
 we use the following criterion:

\begin{Prop}\label{equi condi}
  Let $\Delta_M = \pi \setminus I$ be the set of simple roots for $M$. 
  Let $\Phi_M^+ = \Phi^+ \cap \mathrm{span}(\Delta_M)$ be the positive roots of $M$. 
  We can write $\Phi^{+}(\mathfrak{u}_{P}) = \Phi^+ \setminus \Phi_M^+$ as the set of positive roots corresponding to the unipotent radical of $P$.
  For an element $w \in \mathcal{W}$, the following two conditions are equivalent:
  \begin{enumerate}
      \item[(A)] $w^{-1}(\alpha) \in \Phi^{+}$ for all $\alpha \in \Delta_M$.
      \item[(B)] $w \in \mathcal{W}^{P}, \mr{i.e.}, \  w(\Phi^-) \cap \Phi^{+} \subseteq \Phi^{+}(\mathfrak{u}_{P})$ (or equivalently, $w(\Phi^-) \cap \Phi^{+} \cap \Phi_M^+ = \emptyset$).
  \end{enumerate}
  \end{Prop}
  
  \begin{Pf}\rm
  We prove the equivalence in two directions.
  
  \textbf{(A) $\implies$ (B):}
  Assume condition (A) holds: $w^{-1}(\alpha) \in \Phi^{+}$ for all $\alpha \in \Delta_M$.
  Let $\gamma$ be an arbitrary element in $w(\Phi^-) \cap \Phi^{+}$. This means there exists some $\beta \in \Phi^-$ such that $\gamma = w(\beta)$ and $\gamma \in \Phi^+$. We want to show that $\gamma \notin \Phi_M^+$.
  
  Assume that $\gamma \in \Phi_M^+$. 
  Since elements of $\Phi_M^+$ are non-negative integer linear combinations of simple roots in $\Delta_M$, we can write
  \[
    \gamma = \sum_{\alpha \in \Delta_M} c_\alpha \alpha
  \]
  where $c_\alpha \in \Z_{\ge 0}$ and at least one $c_\alpha > 0$.
  Applying $w^{-1}$ to both sides, we get
  \[
    w^{-1}(\gamma) = \sum_{\alpha \in \Delta_M} c_\alpha w^{-1}(\alpha)
  \]
  By assumption (A), each $w^{-1}(\alpha)$ is a positive root. Since the coefficients $c_\alpha$ are non-negative integers and not all zero, the right-hand side is a non-zero, non-negative integer linear combination of positive roots, which must itself be a positive root. Thus, $w^{-1}(\gamma) \in \Phi^{+}$.
  
  However, we started with $\gamma = w(\beta)$ where $\beta \in \Phi^-$. Applying $w^{-1}$ gives $w^{-1}(\gamma) = \beta$, which is a negative root. This contradicts our finding that $w^{-1}(\gamma) \in \Phi^{+}$.
  Therefore $\gamma \notin \Phi_M^+$.
  Since $\gamma \in \Phi^+$, this implies $\gamma \in \Phi^+ \setminus \Phi_M^+ = \Phi^{+}(\mathfrak{u}_{P})$.
  Thus, $w(\Phi^-) \cap \Phi^{+} \subseteq \Phi^{+}(\mathfrak{u}_{P})$.
  
  \textbf{(B) $\implies$ (A):}
  Assume condition (B) holds: $w(\Phi^-) \cap \Phi^{+} \cap \Phi_M^+ = \emptyset$. We want to show that $w^{-1}(\alpha) \in \Phi^{+}$ for all $\alpha \in \Delta_M$.
  
  Assume that there exists some $\alpha_0 \in \Delta_M$ such that $w^{-1}(\alpha_0) \notin \Phi^{+}$. 
  This implies $w^{-1}(\alpha_0) \in \Phi^-$. 
  Let $\beta_0 = w^{-1}(\alpha_0) \in \Phi^-$.
  Applying $w$ to both sides gives $w(\beta_0) = \alpha_0$.
  Since $\alpha_0$ is a simple root in $\Delta_M$, we have $\alpha_0 \in \Delta_M \subset \Phi_M^+ \subset \Phi^+$. Thus $w(\beta_0) \in \Phi^+$.
  
  Now, we have $\beta_0 \in \Phi^-$ and $w(\beta_0) \in \Phi^+$. Therefore, the element $\alpha_0 = w(\beta_0)$ belongs to the set $w(\Phi^-) \cap \Phi^{+}$.
  Furthermore, we know $\alpha_0 \in \Phi_M^+$.
  Combining these, we find that $\alpha_0 \in (w(\Phi^-) \cap \Phi^{+}) \cap \Phi_M^+$.
  This contradicts assumption (B), which states that this intersection is empty.
  Therefore, our initial assumption that there exists an $\alpha_0 \in \Delta_M$ with $w^{-1}(\alpha_0) \in \Phi^-$ must be false. Consequently, $w^{-1}(\alpha) \in \Phi^{+}$ must hold for all $\alpha \in \Delta_M$.
  \end{Pf}

Using Proposition \ref{equi condi} and the exhaustive table of the Weyl group in Appendix \ref{app:weyl_table}, 
we determine the sets of Kostant representatives $\mathcal{W}^{P_I}$. 
Recall that $\mathcal{W}^{P_I}$ provides a unique set of representatives for the cosets $\mathcal{W} / \mathcal{W}_{\mr{M}_{P}}$.

    \vspace{2mm}
      \subsubsection*{Rank $1$(\ $|I|=1$\ )}
      \begin{flalign*}
      & \bullet I = \{ \al_1 \}: 
        \Delta_{\mr{M}} = \{ \al_2, \al_3 \}, 
        \mr{M} \cong \mr{Sp}_4, 
        \mathcal{W}_{\mr{M}} = \{e, s_2, s_3, s_{23}, s_{32}, s_{232}, s_{323}, s_{2323}\}\\ 
        & \hspace{10mm} \mathcal{W}^{P_{\{ \al_1 \}}} = \{ e, s_{1}, s_{12}, s_{123},s_{1232}, s_{12321} \}  \\
      & \bullet I = \{ \al_2 \}: 
        \Delta_{\mr{M}} = \{ \al_1, \al_3 \},
        \mr{M} \cong \mr{SL}_2 \times \mr{Sp}_2,
        \mathcal{W}_{\mr{M}} = \{e, s_1, s_3, s_{13}\}\\ 
        & \hspace{10mm} \mathcal{W}^{P_{\{ \al_2 \}}} = \{ e, s_{2}, s_{21}, s_{23}, s_{213}, s_{232}, s_{2132}, s_{2321}, s_{21321}, s_{21323}, s_{213213}, s_{2132132} \} \\
      & \bullet I = \{ \al_3 \}:
        \Delta_{\mr{M}} = \{ \al_1, \al_2 \},
        \mr{M} \cong \mr{SL}_3,
        \mathcal{W}_{\mr{M}} = \{e, s_1, s_2, s_{12}, s_{21}, s_{121}\}. \\   
        & \hspace{10mm} \mathcal{W}^{P_{\{ \al_3 \}}} = \{ e, s_{3}, s_{32}, s_{321}, s_{323}, s_{3213}, s_{32132}, s_{321323} \}
      \end{flalign*}

      \subsubsection*{Rank $2$(\ $|I|=2$\ )}
      \begin{flalign*}
      &\bullet I = \{ \al_1, \al_2 \}:
        \Delta_{\mr{M}} = \{ \al_3 \},
        \mr{M} \cong \mr{Sp}_2, 
        \mathcal{W}_{\mr{M}} = \{e, s_3\} \\
        & \hspace{10mm} \mathcal{W}^{P_{\{ \al_1, \al_2 \}}} = \{ e, s_{1}, s_{2}, s_{12}, s_{21}, s_{23}, s_{121}, s_{123}, s_{213}, s_{232}, s_{1213}, s_{1232}, s_{2132}, s_{2321}, && \\
        & \hspace{15mm} s_{12132}, s_{12321}, s_{21321}, s_{21323}, s_{121321}, s_{121323}, s_{213213}, s_{1213213}, s_{2132132}, s_{12132132} \} && \\
      &\bullet I = \{ \al_1, \al_3 \}: 
        \Delta_{\mr{M}} = \{ \al_2 \},
        \mr{M} \cong \mr{SL}_2, 
        \mathcal{W}_{\mr{M}} = \{e, s_2\}  \\
        & \hspace{10mm} \mathcal{W}^{P_{\{ \al_1, \al_3 \}}} = \{ e, s_{1}, s_{3}, s_{12}, s_{13}, s_{32}, s_{123}, s_{132}, s_{321}, s_{323}, s_{1232}, s_{1321}, s_{1323}, s_{3213}, && \\
        & \hspace{15mm} s_{12321}, s_{12323}, s_{13213}, s_{32132}, s_{123213}, s_{132132}, s_{321323}, s_{1232132}, s_{1321323}, s_{12321323} \} && \\
      &\bullet I = \{ \al_2, \al_3 \}: 
        \Delta_{\mr{M}} = \{ \al_1 \},
        \mr{M} \cong \mr{SL}_2, 
        \mathcal{W}_{\mr{M}} = \{e, s_1\}  \\
        & \hspace{10mm} \mathcal{W}^{P_{\{ \al_2, \al_3 \}}} = \{ e, s_{2}, s_{3}, s_{21}, s_{23}, s_{32}, s_{213}, s_{232}, s_{321}, s_{323}, s_{2132}, s_{2321}, s_{2323}, s_{3213}, && \\
        & \hspace{15mm} s_{21321}, s_{21323}, s_{23213}, s_{32132}, s_{213213}, s_{232132}, s_{321323}, s_{2132132}, s_{2321323}, s_{21321323} \} &&
      \end{flalign*}

      \subsubsection*{Rank $3$(\ $|I|=3$\ )}
      \begin{flalign*}
      \bullet I = \{ \al_1, \al_2, \al_3 \} = \pi: \Delta_{\mr{M}} = \emptyset. \\
        \mathcal{W}^{P_{\pi}} = \mathcal{W} &&
      \end{flalign*}
    \vspace{2mm}

    To apply the results from previous work, 
    we express the Kostant representatives $w \cdot \lambda$ in terms of the fundamental dominant weights of the corresponding Levi component.
    The definitions of these weights in the standard basis $\{\varepsilon_1, \varepsilon_2, \varepsilon_3\}$ are given in the following list. 
    The explicit coefficients of $w \cdot \lambda$ and $w \cdot 0$ are detailed in Appendix \ref{app:weight_tables}.

    \begin{table}[H]
      \begin{tabular}{ccc}
        \toprule
        \textbf{type} & \textbf{simple roots} & \textbf{fundamental dominant weight} \\
        \midrule
        \midrule
        $\mr{SL}_2$ & $\{\varepsilon_1 - \varepsilon_2$\} & \{$\frac{\varepsilon_1 - \varepsilon_2}{2}\}$ \\
        $\mr{SL}_3$ & $\{\varepsilon_1 - \varepsilon_2, \varepsilon_2 - \varepsilon_3\}$ & $\{\varepsilon_1, \varepsilon_1 + \varepsilon_2 \}$ \\
        $\mr{Sp}_2$ & $\{2\varepsilon_1\}$ & $\{\varepsilon_1 \}$ \\
        $\mr{Sp}_4$ & $\{\varepsilon_1 - \varepsilon_2, 2\varepsilon_2\}$ & $\{\varepsilon_1, \varepsilon_1 + \varepsilon_2\}$ \\
        \bottomrule 
      \end{tabular}
    \end{table}

      \subsubsection*{Rank $1$($\ |I|=1\ $)}
      \ 

      $\bullet P_{\{\alpha_1\}}$:
      $\mr{M}_{P_{\{\alpha_1\}}} = \mr{GL}_1 \times \mr{Sp}_4$.
      \begin{align*}
        \gamma_1^{\{\al_1\}} &= \varepsilon_1,\\
        \gamma_2^{\{\al_1\}} &= \varepsilon_2,\\
        \gamma_3^{\{\al_1\}} &= \varepsilon_2 + \varepsilon_3
      \end{align*}
      \vspace{2mm}

      $\bullet P_{\{\alpha_2\}}$:
      $\mr{M}_{P_{\{\alpha_2\}}} = \mr{SL}_2 \times \mr{GL}_1 \times \mr{Sp}_2 = \mr{GL}_2 \times \mr{Sp}_2$.
      \begin{align*}
        \gamma_1^{\{\al_2\}} &= \frac{\varepsilon_1 - \varepsilon_2}{2},\\
        \gamma_2^{\{\al_2\}} &= \varepsilon_1 + \varepsilon_2,\\
        \gamma_3^{\{\al_2\}} &= \varepsilon_3
      \end{align*}
      \vspace{2mm}

      $\bullet P_{\{\alpha_3\}}$:
      $\mr{M}_{P_{\{\alpha_3\}}} = \mr{SL}_3 \times \mr{GL}_1 = \mr{GL}_3$.
      \begin{align*}
        \gamma_1^{\{\al_3\}} &= \varepsilon_1\\
        \gamma_2^{\{\al_3\}} &= \varepsilon_1 + \varepsilon_2,\\
        \gamma_3^{\{\al_3\}} &= \varepsilon_1 + \varepsilon_2 + \varepsilon_3
      \end{align*}

      \subsubsection*{Rank $2$($\ |I|=2$\ )}
      \ 

      $\bullet P_{\{\alpha_1, \alpha_2\}}$:
      $\mr{M}_{P_{\{\alpha_1, \alpha_2\}}} = \mr{GL}_1 \times \mr{GL}_1 \times \mr{Sp}_2$.
      \begin{align*}
        \gamma_1^{\{\al_1, \al_2\}} &= \varepsilon_1 \\
        \gamma_2^{\{\al_1, \al_2\}} &= \varepsilon_2 \\
        \gamma_3^{\{\al_1, \al_2\}} &= \varepsilon_3
      \end{align*}
      \vspace{2mm}

      $\bullet P_{\{\alpha_1, \alpha_3\}}$:
      $\mr{M}_{P_{\{\alpha_1, \alpha_3\}}} = \mr{GL}_1 \times \mr{SL}_2 \times \mr{GL}_1 = \mr{GL}_{1} \times \mr{GL}_2$. 
      \begin{align*}
        \gamma_1^{\{\al_1, \al_3\}} &= \varepsilon_1 \\
        \gamma_2^{\{\al_1, \al_3\}} &= \frac{\varepsilon_2 - \varepsilon_3}{2} \\
        \gamma_3^{\{\al_1, \al_3\}} &= \varepsilon_2 + \varepsilon_3
      \end{align*}
      \vspace{2mm}

      $\bullet P_{\{\alpha_2, \alpha_3\}}$:
      $\mr{M}_{P_{\{\alpha_2, \alpha_3\}}} = \mr{SL}_2 \times \mr{GL}_1 \times \mr{GL}_1 = \mr{GL}_2 \times \mr{GL}_1$.
      \begin{align*}
        \gamma_1^{\{\al_2, \al_3\}} &= \frac{\varepsilon_1 - \varepsilon_2}{2} \\
        \gamma_2^{\{\al_2, \al_3\}} &= \varepsilon_1 + \varepsilon_2 \\
        \gamma_3^{\{\al_2, \al_3\}} &= \varepsilon_3
      \end{align*}

      \subsubsection*{Rank $3$($\ |I|=3\ )$}
      \ 
      
      $\bullet P_{\pi}$:
      $\mr{M}_{P_{\pi}} = \mr{GL}_1 \times \mr{GL}_1 \times \mr{GL}_1$.
      \begin{align*}
        \gamma_1^{\pi} &= \varepsilon_1 \\
        \gamma_2^{\pi} &= \varepsilon_2 \\
        \gamma_3^{\pi} &= \varepsilon_3
      \end{align*}

    \vspace{2mm}

 
 
\section{Parity Conditions in Cohomology}
In this section, we determine the parity conditions for the coefficients of each $\gamma_{i}^{I}$
required for the non-vanishing of the associated local systems.

By the definition of the sheaf $\wt{\mathcal{M}}$ associated with the irreducible representation $\mathcal{M}$, 
any element in the intersection $\Gamma_{\mr{M}_{P}} \cap K_{\infty}^{\mr{M}_{P}} \cap Z(M_P)$ must act trivially on the representation space $\mathcal{M}$,
where $Z(M_P)$ denotes the center of $M_P$.
Indeed, 
for any local section $f \in \wt{\mathcal{M}}(U)$ and $u \in \pi^{-1}(U)$ for an open subset $U$ of $\mr{S}_{\Gamma} = \Gamma \backslash G / K_\infty$,
an element $\gamma \in \Gamma \cap K_{\infty} \cap Z(G)$ satisfies
\[
  \gamma f(u) = f(\gamma u) = f(u)
\]
where the last equality holds because $\gamma \in \Z(G)$ acts trivially on the symmetric space $G/K_\infty$ ($\mr{i.e.}, \gamma u = u$).
Thus, if $\gamma$ does not act trivially on $\mathcal{M}$,
the only possible section is $f = 0$,
which implies that the sheaf $\wt{\mathcal{M}} = 0$.

\subsection{Parabolic of rank $3$ (Borel subgroup)}

The Levi subgroup of the minimal parabolic subgroup
 $\mr{B} = P_{\pi}$ is the maximal torus $\mr{M}_{\pi} = \mr{T}$.
 The following three diagonal matrices are contained in
 $\Gamma_{T} \cap K_{\infty}^{T} \cap Z(T)$:
 \begin{align*}
   T_1 &=
   \begin{pmatrix}
     -1  & 0   & 0   & 0        & 0        & 0 \\
     0   & 1   & 0   & 0        & 0        & 0 \\
     0   & 0   & 1   & 0        & 0        & 0 \\
     0   & 0   & 0   & -1       & 0        & 0 \\
     0   & 0   & 0   & 0        & 1        & 0 \\
     0   & 0   & 0   & 0        & 0        & 1
   \end{pmatrix} \\
   T_2 &=
   \begin{pmatrix}
     1   & 0   & 0   & 0        & 0        & 0 \\
     0   & -1  & 0   & 0        & 0        & 0 \\
     0   & 0   & 1   & 0        & 0        & 0 \\
     0   & 0   & 0   & 1        & 0        & 0 \\
     0   & 0   & 0   & 0        & -1       & 0 \\
     0   & 0   & 0   & 0        & 0        & 1
   \end{pmatrix} \\
   T_3 &=
   \begin{pmatrix}
     1   & 0   & 0   & 0        & 0        & 0 \\
     0   & 1   & 0   & 0        & 0        & 0 \\
     0   & 0   & -1  & 0        & 0        & 0 \\
     0   & 0   & 0   & 1        & 0        & 0 \\
     0   & 0   & 0   & 0        & 1        & 0 \\
     0   & 0   & 0   & 0        & 0        & -1
   \end{pmatrix}.
 \end{align*} 
If the highest weight of the representation $\mathcal{M}$ on $T$ is of the form
\[
  \lambda' = n_1 \varepsilon_1 + n_2 \varepsilon_2 + n_3 \varepsilon_3,
\]
then the above three diagonal matrices act on the corresponding highest weight vector $v$ as
\begin{align*}
  T_1 \cdot v &= (-1)^{n_1} v, \\
  T_2 \cdot v &= (-1)^{n_2} v, \\
  T_3 \cdot v &= (-1)^{n_3} v.
\end{align*}
Therefore, $n_1, n_2$ and $n_3$ must be even; otherwise $\mathcal{M}$ vanishes.
Since the basis is given by
$\gamma_1^{\pi} = \varepsilon_1, 
\gamma_2^{\pi} = \varepsilon_2, 
\gamma_3^{\pi} = \varepsilon_3$,
we can deduce the following lemma.
\begin{Lem}
  Let $w \cdot \lambda = m_1 \gamma_1^{\pi} + m_2 \gamma_2^{\pi} + m_3 \gamma_3^{\pi}$.
  The local system $\wt{\mathcal{M}}_{w \cdot \lambda}$ is non-zero 
  only if $m_1, m_2$, and $m_3$ are all even.
\end{Lem}

\subsection{Parabolics of rank 2}

\begin{itemize}

  \item \textbf{Case $I = \{\alpha_1, \alpha_2\}$}: 
    The Levi subgroup is
    $\mr{M}_{P_{\{\alpha_1, \alpha_2\}}} = \mr{GL}_1 \times \mr{GL}_1 \times \mr{Sp}_2$.
    The basis for $\mr{M}_{\{\al_1, \al_2\}}$ is given by
    \[
      \gamma_1^{\{\al_1, \al_2\}} = \varepsilon_1, \quad
    \gamma_2^{\{\al_1, \al_2\}} = \varepsilon_2, \quad
    \gamma_3^{\{\al_1, \al_2\}} = \varepsilon_3.
    \]
    The following three elements are contained in $\Gamma_{\mr{M}_{P}} \cap K_{\infty}^{M_{P}} \cap Z(\mr{M}_{P})$: 
    \begin{align*}
      c_1 &= (-1, 1, I_2), \\
      c_2 &= (1, -1, I_2), \\
      c_3 &= (1, 1, -I_2).
    \end{align*}
    Let $\lambda' = m_1 \gamma_1^{\{\al_1, \al_2\}} + m_2 \gamma_2^{\{\al_1, \al_2\}} + m_3 \gamma_3^{\{\al_1, \al_2\}}$.
    The above three elements act on the corresponding weight vector $v$ as
    \begin{align*}
      c_1 \cdot v &= (-1)^{m_1} v, \\
      c_2 \cdot v &= (-1)^{m_2} v, \\
      c_3 \cdot v &= (-1)^{m_3} v.
    \end{align*}
    Therefore $m_1, m_2$, and $m_3$ must be even; otherwise $\mathcal{M}$ vanishes.
    We summarize this in the following lemma:
    \begin{Lem}
      Let $w \cdot \lambda = m_1 \gamma_1^{\{\al_1, \al_2\}} + m_2 \gamma_2^{\{\al_1, \al_2\}} + m_3 \gamma_3^{\{\al_1, \al_2\}}$ be the highest weight of the representation of $M_P$.
      The associated local system $\wt{\mathcal{M}}_{w \cdot \lambda}$ vanishes 
      if any of $m_1, m_2$, or $m_3$ is odd.
    \end{Lem}

    \item \textbf{Case $I = \{\alpha_1, \alpha_3\}$}: 
    The Levi subgroup is
    $\mr{M}_{P_{\{\alpha_1, \alpha_3\}}} = \mr{GL}_1 \times \mr{SL}_2 \times \mr{GL}_1 = \mr{GL}_1 \times \mr{GL}_2$.
    The basis for $\mr{M}_{\{\al_1, \al_3\}}$ is given by
    \[
      \gamma_1^{\{\al_1, \al_3\}} = \varepsilon_1, \quad
    \gamma_2^{\{\al_1, \al_3\}} = \frac{1}{2}(\varepsilon_2 - \varepsilon_3), \quad
    \gamma_3^{\{\al_1, \al_3\}} = \varepsilon_2 + \varepsilon_3.
    \]
    If the highest weight of the $\mr{GL}_2$ block is written as 
    \[
      m_2 \gamma_{2}^{\{\al_1, \al_3\}} + m_3 \gamma_{3}^{\{\al_1, \al_3\}} = \frac{m_2}{2}(\varepsilon_2 - \varepsilon_3) + m_3 (\varepsilon_2 + \varepsilon_3),
    \]
    then the corresponding representation is $\mathcal{V}_{m_2,m_3} = \mr{Sym}^{m_2}(\Q^2) \otimes \mr{det}^{(m_3 - \frac{m_2}{2})}$, since 
    $m_2 \gamma_{2}^{\{\al_1, \al_3\}} + m_3 \gamma_{3}^{\{\al_1, \al_3\}} = m_2 \varepsilon_2 + (m_3 - \frac{m_2}{2}) (\varepsilon_2 + \varepsilon_3)$.
    Note that the exponent $m_3 - \frac{m_2}{2}$ must be an integer, which implies $m_2$ is even.
    The following two elements are contained in $\Gamma_{\mr{M}_{P}} \cap K_{\infty}^{M_{P}} \cap Z(\mr{M}_{P})$:
    \begin{align*}
      c_1 &= (1, -I_2), \\
      c_2 &= (-1, I_2).
    \end{align*}
    Let $\lambda' = m_1 \gamma_1^{\{\al_1, \al_3\}} + m_2 \gamma_2^{\{\al_1, \al_3\}} + m_3 \gamma_3^{\{\al_1, \al_3\}}$.
    The above two matrices act on the corresponding weight vector $v$ as
    \begin{align*}
      c_1 \cdot v = (-1)^{m_2} v, \\
      c_2 \cdot v = (-1)^{m_1} v.
    \end{align*}
    For the local system to be non-zero, $m_1$ and $m_2$ must be even.

    If $m_2 = 0$, then $\mathcal{V}_{m_2, m_3} = \mr{det}^{m_3}$, which is one-dimensional.  
    This means that $\mr{H}^q(S^{\mr{GL}_2(\Z)}, \mathcal{V}_{m_2, m_3}) = 0$ for $q \not = 0$.
    For $q = 0$, we have $\mr{H}^0(S^{\mr{GL}_2(\Z)}, \mathcal{V}_{0, m_3}) = \mathcal{V}_{0, m_3}^{\mr{GL}_2(\Z)}$.
    The matrix $\mr{diag}(1, -1) \in \mr{GL}_2(\Z)$ must act trivially on $v$ in this invariant space.
    Since its action on $\mr{det}^{m_3}$ is given by $(-1)^{m_3}$, we deduce that $m_3$ is even.
    
    We summarize this result in the following lemma.
    \begin{Lem}
      Let $w \cdot \lambda = m_1 \gamma_1^{\{\al_1, \al_3\}} + m_2 \gamma_2^{\{\al_1, \al_3\}} + m_3 \gamma_3^{\{\al_1, \al_3\}}$ be the highest weight of the representation of $M_P$.
      The local system $\wt{\mathcal{M}}_{w \cdot \lambda}$ vanishes 
      if either $m_1$ or $m_2$ is odd.
      Moreover, if $m_2 = 0$ and $m_3$ is odd, then $\mr{H}^q(S^{\mr{M}_{P}}, \mathcal{M}_{w \cdot \lambda}) = 0$ for all $q$.
    \end{Lem}

  \item \textbf{Case $I = \{\alpha_2, \alpha_3\}$}:
    The Levi subgroup is
    $\mr{M}_{P_{\{\alpha_2, \alpha_3\}}} = \mr{SL}_2 \times \mr{GL}_1 \times \mr{GL}_1 = \mr{GL}_2 \times \mr{GL}_1$.
    The basis for $\mr{M}_{\{\al_2, \al_3\}}$ is given by
    \[
      \gamma_1^{\{\al_2, \al_3\}} = \frac{1}{2}(\varepsilon_1 - \varepsilon_2), \quad 
      \gamma_2^{\{\al_2, \al_3\}} = \varepsilon_1 + \varepsilon_2, \quad 
      \gamma_3^{\{\al_2, \al_3\}} = \varepsilon_3.
    \]
    By a similar argument to the case $I = \{\al_1, \al_3\}$, we can get
      \begin{Lem}
        Let $w \cdot \lambda = m_1 \gamma_1^{\{\al_2, \al_3\}} + m_2 \gamma_2^{\{\al_2, \al_3\}} + m_3 \gamma_3^{\{\al_2, \al_3\}}$ be the highest weight of the representation of $M_P$.
        The local system $\wt{\mathcal{M}}_{w \cdot \lambda}$ vanishes 
        if either $m_1$ or $m_3$ is odd.
        Moreover, if $m_1 = 0$ and $m_2$ is odd, then $\mr{H}^q(S^{\mr{M}_{P}}, \mathcal{M}_{w \cdot \lambda}) = 0$ for all $q$.
      \end{Lem}
  \end{itemize}

\subsection{Parabolics of rank 1}
\begin{itemize}
  \item \textbf{Case $I = \{\alpha_1\}$}: 
The Levi subgroup is
$\mr{M}_{P_{\{\alpha_1\}}} = \mr{GL}_1 \times \mr{Sp}_4$.
The basis for $\mr{M}_{\{\al_1\}}$ is given by
\[
  \gamma_1^{\{\al_1\}} = \varepsilon_1, \quad
  \gamma_2^{\{\al_1\}} = \varepsilon_2, \quad
  \gamma_3^{\{\al_1\}} = \varepsilon_2 + \varepsilon_3.
\]
The following matrices are contained in $\Gamma_{\mr{M}_{P}} \cap K_{\infty}^{M_{P}} \cap Z(\mr{M}_{P})$:
\begin{align*}
  c_1 &= (1, \mr{diag}(-1,-1,-1,-1)) \\
  c_2 &= (-1, \mr{diag}(1,1,1,1))
\end{align*}
Let $\lambda' = m_1 \gamma_1^{\{\al_1\}} + m_2 \gamma_2^{\{\al_1\}} + m_3 \gamma_3^{\{\al_1\}} = m_1\varepsilon_1 + (m_2 + m_3)\varepsilon_2 + m_3\varepsilon_3$.
The above three elements act on the corresponding weight vector $v$ as
\begin{align*}
  c_1 \cdot v &= (-1)^{m_2 + 2m_3} \\
  c_2 \cdot v &= (-1)^{m_1}
\end{align*}
Therefore $m_1$ and $m_2$ must be even, and we get
\begin{Lem}
  Let $w \cdot \lambda = m_1 \gamma_1^{\{\al_1\}} + m_2 \gamma_2^{\{\al_1\}} + m_3 \gamma_3^{\{\al_1\}}$ be the highest weight of the representation of $M_P$.
  The local system $\wt{\mathcal{M}}_{w \cdot \lambda}$ vanishes 
  if either $m_1$ or $m_2$ is odd.
\end{Lem}

\item \textbf{Case $I = \{\alpha_2\}$}:
The Levi subgroup is
$\mr{M}_{P_{\{\al_2\}}} = \mr{SL}_2 \times \mr{GL}_1 \times \mr{Sp}_2 = \mr{GL}_2 \times \mr{Sp}_2$.
The basis for $\mr{M}_{\{\al_2\}}$ is given by
\[
  \gamma_1^{\{\al_2\}} = \frac{1}{2}(\varepsilon_1 -\varepsilon_2), \quad 
  \gamma_2^{\{\al_2\}} = \varepsilon_1 + \varepsilon_2, \quad 
  \gamma_3^{\{\al_2\}} = \varepsilon_3.
\]
By a similar argument to the case $I = \{\al_2, \al_3\}$, we can get
      \begin{Lem}
        Let $w \cdot \lambda = m_1 \gamma_1^{\{\al_2\}} + m_2 \gamma_2^{\{\al_2\}} + m_3 \gamma_3^{\{\al_2\}}$ be the highest weight of the representation of $M_P$.
        The local system $\wt{\mathcal{M}}_{w \cdot \lambda}$ vanishes 
        if either $m_1$ or $m_3$ is odd.
        Moreover, if $m_1 = 0$ and $m_2$ is odd, then $\mr{H}^q(S^{\mr{M}_{P}}, \mathcal{M}_{w \cdot \lambda}) = 0$ for all $q$.
      \end{Lem}

\item \textbf{Case $I = \{\alpha_3\}$}:
The Levi subgroup is
$\mr{M}_{P_{\{\al_3\}}} = \mr{SL}_3 \times \mr{GL}_1 = \mr{GL}_3$.
The basis for $\mr{M}_{\{\al_3\}}$ is given by
\[
  \gamma_1^{\{\al_3\}} = \varepsilon_1, \quad
  \gamma_2^{\{\al_3\}} = \varepsilon_1 + \varepsilon_2, \quad 
  \gamma_3^{\{\al_3\}} = \varepsilon_1 + \varepsilon_2 + \varepsilon_3.
\]
The following element is contained in $\Gamma_{\mr{M}_{P}} \cap K_{\infty}^{M_{P}} \cap Z(\mr{M}_{P})$:
\begin{equation*}
  c = \mr{diag}(-1,-1,-1)
\end{equation*}
Let $\lambda' = m_1 \gamma_1^{\{\al_3\}} + m_2 \gamma_2^{\{\al_3\}} + m_3 \gamma_3^{\{\al_3\}}$.
The above element acts on the corresponding weight vector $v$ as
\begin{equation*}
  c \cdot v = (-1)^{m_1 + 2m_2 + 3m_3}
\end{equation*}
Therefore $m_1 + m_3$ must be even, and we get
\begin{Lem}
  The local system $\wt{\mathcal{M}}_{w \cdot \lambda}$ vanishes
  if $m_1 + m_3$ is odd.
\end{Lem}
\end{itemize}

\subsection{Summary of non-vanishing representatives}

We denote by $\overline{\mathcal{W}^{P_{I}}}$ the subset of Kostant representatives 
for which the corresponding local system satisfies the parity conditions and cohomology does not vanish for all degrees.

Based on the coefficients calculated in Appendix \ref{app:weight_tables}, 
we identify these subsets as follows:

\begin{align*}
  \overline{\mathcal{W}^{P_{\{\al_1\}}}} &= \{e, *s_{12}, *s_{123}, s_{12321}\},\\
  \overline{\mathcal{W}^{P_{\{\al_2\}}}} &= \{e, s_{232}, *s_{213}, *s_{232}, s_{2132}, *s_{2321}, s_{21323}\},\\
  \overline{\mathcal{W}^{P_{\{\al_3\}}}} &= \{e, s_3, *s_{32}, *s_{321}, *s_{323}, *s_{3213}, s_{32132}, s_{321323}\},\\
  \overline{\mathcal{W}^{P_{\{\al_1, \al_2\}}}} &= \{e, s_{121}, s_{232}, s_{1213}, s_{2132}, s_{12321}, s_{21323}, s_{12132132}\},\\
  \overline{\mathcal{W}^{P_{\{\al_1, \al_3\}}}} &= \{e, s_3, *s_{32}, *s_{132}, s_{1321}, *s_{1323}, s_{12321}, s_{13213}, s_{32132}, s_{123213}, s_{321323}, *s_{1232132}\},\\
  \overline{\mathcal{W}^{P_{\{\al_2, \al_3\}}}} &= \{e, s_3, *s_{213}, s_{232}, s_{2132}, *s_{2321}, s_{2323}, *s_{3213}, s_{21323}, *s_{23213}, s_{32132}, s_{321323}\},\\
  \overline{\mathcal{W}^{P_{\pi}}} &= \{e, s_{3}, s_{121}, s_{232}, s_{1213}, s_{1321}, s_{2132}, s_{2323}, s_{12321}, s_{13213}, s_{21323}, s_{32132},\\
  & \quad s_{123213}, s_{321323}, s_{12132132}, s_{121321323}\}
\end{align*}

 
 
\section{Boundary cohomology}
In this section, we calculate the cohomology of the boundary by using the spectral sequence associated with the stratification of the Borel-Serre compactification.
The boundary $\partial \overline{S}$ defines a spectral sequence in cohomology:
\[
  E_{1}^{p,q} \Rightarrow \mr{H}^{p+q}(\partial S, \uQ).
\]

Since the $\Q$-split rank of $\mr{Sp}_6$ is three, the spectral sequence consists of exactly three columns: $E_{1}^{0,q}, E_{1}^{1,q}$, and $E_{1}^{2,q}$.
We first consider the following sequence of $d_1$-differentials
\[
  0 \to E_{1}^{0,q} \overset{d_{1}^{0,q}}{\longrightarrow} E_{1}^{1,q} \overset{d_{1}^{1,q}}{\longrightarrow} E_{1}^{2,q} \to 0
\]
where $d_{1}^{p,q}$ are the differentials of the $E_1$-page. 
The terms on the $E_2$-page are given by
\begin{align*}
  E_{2}^{0,q} &= \mr{Ker}(d_{1}^{0,q}) \\
  E_{2}^{1,q} &= \mr{Ker}(d_{1}^{1,q})/\mr{Im}(d_{1}^{0,q}) \\
  E_{2}^{2,q} &= \mr{Coker}(d_{1}^{1,q})
\end{align*}

Next, we analyze the $d_2$-differentials,
\[
  0 \to E_{2}^{0,q} \overset{d_{2}^{0,q}}{\longrightarrow} E_{2}^{2, q-1} \to 0.
\]
The resulting $E_3$-terms are
\begin{align*}
  E_3^{0,q} &= \mr{Ker}(d_{2}^{0,q}), \\
  E_3^{1,q} &= E_2^{1,q}, \\
  E_3^{2,q-1} &= \mr{Coker}(d_{2}^{0,q}).
\end{align*}

Finally, all higher differentials $d_r$ ($r \geq 3$) vanish identically.
Consequently, the spectral sequence degenerates at the $E_3$-page,
and the boundary cohomology is determined by the direct sum
\[
  \mr{H}^{k}(\partial S, \uQ) = \bigoplus_{p+q = k}E_{3}^{p,q}
\]

\subsection{$E_{ 1}$-page}

The following is the set of non-vanishing Kostant representatives $\overline{\mathcal{W}^{P_{I}}}$ for each standard parabolic subgroup,
determined by the parity conditions established in the previous section.
\begin{align*}
  \overline{\mathcal{W}^{P_{\{\al_1\}}}} &= \{e, *s_{12}, *s_{123}, s_{12321}\},\\
  \overline{\mathcal{W}^{P_{\{\al_2\}}}} &= \{e, *s_{213}, s_{232}, s_{2132}, *s_{2321}, s_{21323}\},\\
  \overline{\mathcal{W}^{P_{\{\al_3\}}}} &= \{e, s_3, *s_{32}, *s_{321}, *s_{323}, *s_{3213}, s_{32132}, s_{321323}\},\\
  \overline{\mathcal{W}^{P_{\{\al_1, \al_2\}}}} &= \{e, s_{121}, s_{232}, s_{1213}, s_{2132}, s_{12321}, s_{21323}, s_{12132132}\},\\
  \overline{\mathcal{W}^{P_{\{\al_1, \al_3\}}}} &= \{e, s_3, *s_{32}, *s_{132}, s_{1321}, *s_{1323}, s_{12321}, s_{13213}, s_{32132}, s_{123213}, s_{321323}, *s_{1232132}\},\\
  \overline{\mathcal{W}^{P_{\{\al_2, \al_3\}}}} &= \{e, s_3, *s_{213}, s_{232}, s_{2132}, *s_{2321}, s_{2323}, *s_{3213}, s_{21323}, *s_{23213}, s_{32132}, s_{321323}\},\\
  \overline{\mathcal{W}^{P_{\pi}}} &= \{e, s_{3}, s_{121}, s_{232}, s_{1213}, s_{1321}, s_{2132}, s_{2323}, s_{12321}, s_{13213}, s_{21323}, s_{32132},\\
  & \quad s_{123213}, s_{321323}, s_{12132132}, s_{121321323}\}
\end{align*}

In this section,
$\mathcal{S}_{k}$ denotes the space of cusp forms of $\mr{SL}_2(\Z)$ of weight $k$,
$\overline{\mathcal{S}_{k}}$ denotes the space of anti-holomorphic cusp forms (which is isomorphic to $\mathcal{S}_k$),
and $\mathcal{E}_{k}$ denotes the space of Eisenstein series of $\mr{SL}_2(\Z)$ of weight $k$.
It is known that
\begin{align*}
  \dim_\C \mathcal{S}_{12l + 2 + i} &=
  \begin{cases}
    l - 1 & i = 0 \\
    l     & i = 2,4,6,8 \\
    l + 1 & i = 10 \\
    0     & i \  \text{is odd}
  \end{cases} \\
\dim_\C \mathcal{E}_{k} &= 
  \begin{cases}
    1 & k \text{ is even} \\
    0 & k \text{ is odd}
  \end{cases}
\end{align*}

\subsubsection{$p = 0$}

For $p = 0$, the $E_1$-term can be written as
\[
  E_1^{0,q} = \bigoplus_{i = 1}^{3} \mr{H}^q(\partial_{P_{\{\al_i\}}}, \uQ).
\]

We compute each face $\mr{H}^q(\partial_{P_{\{\al_i\}}}, \uQ)$.

\paragraph{(1) Case $I = \{\alpha_1\}$:}
In this case, the Levi component is $\mr{M} \cong \mr{GL}_1 \times \mr{Sp}_4$.
The set of non-vanishing Kostant representatives is 
\[
  \overline{\mathcal{W}^{P_{\{\al_1\}}}} = \{e, *s_{12}, *s_{123}, s_{12321}\}.
\]
For $w \cdot \lambda = m_1 \gamma_1^{\{\al_1\}} + m_2 \gamma_2^{\{\al_1\}} + m_3 \gamma_3^{\{\al_1\}}$, 
the pair $(m_2, m_3)$ corresponds to the highest weight of the representation of $\mr{Sp}_4$.
We have $\mr{H}^q(S^{\mr{Sp}_4}, \wt{\mathcal{M}}) = 0$ for all $q > 4$.
The cohomology of the face is
\begin{align*}
  \mr{H}^q(\partial_{P_{\{\al_1\}}}, \uQ) =
    \begin{cases}
      \mr{H}^0(S^{\mr{Sp}_4}, \uQ)_e  & q = 0 \\
      \vspace{2mm}
      \mr{H}^1(S^{\mr{Sp}_4}, \uQ)_e  & q = 1 \\
      \vspace{2mm}
      \mr{H}^2(S^{\mr{Sp}_4}, \uQ)_e \oplus \mr{H}^0(S^{\mr{Sp}_4}, \mathcal{M}_{(0,1)})_{s_{12}} & q = 2 \\
      \vspace{2mm}
      \mr{H}^3(S^{\mr{Sp}_4}, \uQ)_e \oplus \mr{H}^1(S^{\mr{Sp}_4}, \mathcal{M}_{(0,1)})_{s_{12}} \oplus \mr{H}^0(S^{\mr{Sp}_4}, \mathcal{M}_{(0,1)})_{s_{123}} & q = 3 \\
      \vspace{2mm}
      \mr{H}^4(S^{\mr{Sp}_4}, \uQ)_e \oplus \mr{H}^2(S^{\mr{Sp}_4}, \mathcal{M}_{(0,1)})_{s_{12}} \oplus \mr{H}^1(S^{\mr{Sp}_4}, \mathcal{M}_{(0,1)})_{s_{123}} & q = 4 \\
      \vspace{2mm}
      \mr{H}^0(S^{\mr{Sp}_4}, \uQ)_{s_{12321}} \oplus \mr{H}^3(S^{\mr{Sp}_4}, \mathcal{M}_{(0,1)})_{s_{12}} \oplus \mr{H}^2(S^{\mr{Sp}_4}, \mathcal{M}_{(0,1)})_{s_{123}}   & q = 5 \\
      \vspace{2mm}
      \mr{H}^1(S^{\mr{Sp}_4}, \uQ)_{s_{12321}} \oplus \mr{H}^4(S^{\mr{Sp}_4}, \mathcal{M}_{(0,1)})_{s_{12}} \oplus \mr{H}^3(S^{\mr{Sp}_4}, \mathcal{M}_{(0,1)})_{s_{123}} & q = 6 \\
      \vspace{2mm}
      \mr{H}^2(S^{\mr{Sp}_4}, \uQ)_{s_{12321}} \oplus \mr{H}^4(S^{\mr{Sp}_4}, \mathcal{M}_{(0,1)})_{s_{123}}    & q = 7 \\
      \vspace{2mm}
      \mr{H}^3(S^{\mr{Sp}_4}, \uQ)_{s_{12321}}  & q = 8 \\
      \vspace{2mm}
      \mr{H}^4(S^{\mr{Sp}_4}, \uQ)_{s_{12321}}  & q = 9 \\
      0 & \text{otherwise}
    \end{cases}
\end{align*}

Based on known results for $\mr{Sp}_4(\Z)$ \cite{sp4}, the Eisenstein cohomology and the interior cohomology satisfy
\begin{align*}
  \mr{H}^q_{\mr{Eis}}(S^{\mr{Sp}_4}, \uQ) &=
    \begin{cases}
      \Q & q = 0, 2 \\
      0  & \text{otherwise}
    \end{cases}, \\
  \mr{H}^q_{!}(S^{\mr{Sp}_4}, \uQ) &= 0 \quad \text{for all } q.
\end{align*}
and we have
\begin{align*}
  \mr{H}^q_{\mr{Eis}}(S^{\mr{Sp}_4}, \mathcal{M}_{(0,1)} \otimes \C) &=
    \begin{cases}
      \Q \oplus \mathcal{S}_{4} \oplus \overline{\mathcal{S}_{4}} \oplus \mathcal{S}_{6} = \Q & q = 3 \\
      0  & \text{otherwise}
    \end{cases}, \\
  \mr{H}^q_{!}(S^{\mr{Sp}_4}, \mathcal{M}_{(0,1)}) &= 0 \quad \text{for all } q.
\end{align*}
Consequently, we obtain:
\begin{align*}
  \mr{H}^q(\partial_{P_{\{\al_1\}}}, \uQ) &\cong
   \begin{cases}
    \Q_e & q = 0, 2 \\
    \Q_{s_{12321}} \oplus \Q_{s_{12}} & q = 5 \\
    \Q_{s_{123}} & q = 6 \\
    \Q_{s_{12321}} & q = 7 \\
    0 & \text{otherwise}
   \end{cases}.
\end{align*}

\paragraph{(2) Case $I = \{\alpha_2\}$:}
In this case, the Levi component is $\mr{M} \cong \mr{GL}_2 \times \mr{Sp}_2$.
The set of non-vanishing Kostant representatives is 
\[
  \overline{\mathcal{W}^{P_{\{\al_2\}}}} = \{e, *s_{213}, s_{232}, s_{2132}, *s_{2321}, s_{21323}\}.
\]
For $w \cdot \lambda = m_1 \gamma_1^{\{\al_2\}} + m_2 \gamma_2^{\{\al_2\}} + m_3 \gamma_3^{\{\al_2\}}$, 
the pair $((m_1, m_2), m_3)$ corresponds to the highest weights of $\mr{GL}_2$ and $\mr{Sp}_2$ respectively.
We have $\mr{H}^q(S^{\mr{GL}_2 \times \mr{Sp}_2}, \wt{\mathcal{M}}) = 0$ for $q > 2$, 
as $\mr{H}^q(S^{\mr{GL}_2}, \uQ) = \mr{H}^q(S^{\mr{Sp}_2}, \uQ) = 0$ for all $q > 1$.

Then the cohomology of the face is
\begin{align*}
  \mr{H}^q(\partial_{P_{\{\al_2\}}}, \uQ)
   &=
   \begin{cases}
    \vspace{2mm}
    \mr{H}^0(S^{\mr{GL}_2 \times \mr{Sp}_2}, \uQ)_e                                                   & q = 0 \\
    \vspace{2mm}
    \mr{H}^1(S^{\mr{GL}_2 \times \mr{Sp}_2}, \uQ)_e                                                      & q = 1 \\
    \vspace{2mm}
    \mr{H}^2(S^{\mr{GL}_2 \times \mr{Sp}_2}, \uQ)_e                                                      & q = 2 \\
    \vspace{2mm}
    \mr{H}^0(S^{\mr{GL}_2 \times \mr{Sp}_2}, \wt{\mathcal{M}}_{((4,-2),0)})_{s_{232}} \oplus \mr{H}^0(S^{\mr{GL}_2 \times \mr{Sp}_2}, \wt{\mathcal{M}}_{((2,-2),2)})_{s_{213}}  & q = 3 \\
    \vspace{2mm}
    \mr{H}^1(S^{\mr{GL}_2 \times \mr{Sp}_2}, \wt{\mathcal{M}}_{((4,-2),0)})_{s_{232}} \oplus \mr{H}^1(S^{\mr{GL}_2 \times \mr{Sp}_2}, \wt{\mathcal{M}}_{((2,-2),2)})_{s_{213}} \\ 
      \vspace{2mm}\hspace{10mm} \oplus \mr{H}^0(S^{\mr{GL}_2 \times \mr{Sp}_2}, \wt{\mathcal{M}}_{((2,-3),2)})_{s_{2132}} \oplus \mr{H}^0(S^{\mr{GL}_2 \times \mr{Sp}_2}, \wt{\mathcal{M}}_{((2,-3),2)})_{s_{2321}} & q = 4 \\
      \vspace{2mm}
    \mr{H}^2(S^{\mr{GL}_2 \times \mr{Sp}_2}, \wt{\mathcal{M}}_{((4,-2),0)})_{s_{232}} \oplus \mr{H}^2(S^{\mr{GL}_2 \times \mr{Sp}_2}, \wt{\mathcal{M}}_{((2,-2),2)})_{s_{213}} \\ 
      \vspace{2mm}\hspace{10mm}\oplus \mr{H}^1(S^{\mr{GL}_2 \times \mr{Sp}_2}, \wt{\mathcal{M}}_{((2,-3),2)})_{s_{2132}} \oplus \mr{H}^0(S^{\mr{GL}_2 \times \mr{Sp}_2}, \wt{\mathcal{M}}_{((2,-3),2)})_{s_{2321}}  \\
      \vspace{2mm}\hspace{10mm}\oplus \mr{H}^0(S^{\mr{GL}_2 \times \mr{Sp}_2}, \wt{\mathcal{M}}_{((0,-4),2)})_{s_{21323}}  & q = 5 \\
      \vspace{2mm}
    \mr{H}^2(S^{\mr{GL}_2 \times \mr{Sp}_2}, \wt{\mathcal{M}}_{((2,-3),2)})_{s_{2132}} \oplus \mr{H}^0(S^{\mr{GL}_2 \times \mr{Sp}_2}, \wt{\mathcal{M}}_{((2,-3),2)})_{s_{2321}} \\ 
      \vspace{2mm}\hspace{10mm} \oplus \mr{H}^1(S^{\mr{GL}_2 \times \mr{Sp}_2}, \wt{\mathcal{M}}_{((0,-4),2)})_{s_{21323}}                   & q = 6 \\
      \vspace{2mm}
    \mr{H}^2(S^{\mr{GL}_2 \times \mr{Sp}_2}, \wt{\mathcal{M}}_{((0,-4),2)})_{s_{21323}}                     & q = 7 \\
     \vspace{2mm}
    0 & \text{otherwise}
   \end{cases} \\
  &=
   \begin{cases}
    \vspace{2mm}
    \mr{H}^0(S^{\mr{GL}_2}, \uQ)_e 
        \otimes \mr{H}^0(S^{\mr{Sp}_2}, \uQ)_e   & q = 0 \\
        \vspace{2mm}
    \left[\mr{H}^1(S^{\mr{GL}_2}, \uQ)_e 
        \otimes \mr{H}^0(S^{\mr{Sp}_2}, \uQ)_e\right] \\
        \vspace{2mm}\hspace{5mm} \oplus \left[ \mr{H}^0(S^{\mr{GL}_2}, \uQ)_e 
        \otimes \mr{H}^1(S^{\mr{Sp}_2}, \uQ)_e\right]   & q = 1 \\
        \vspace{2mm}
    \mr{H}^1(S^{\mr{GL}_2}, \uQ)_e 
        \otimes \mr{H}^1(S^{\mr{Sp}_2}, \uQ)_e                                              & q = 2 \\
        \vspace{2mm}
    \left[\mr{H}^0(S^{\mr{GL}_2}, \wt{\mathcal{M}}_{(4,-2)})_{s_{232}} 
        \otimes \mr{H}^0(S^{\mr{Sp}_2},\uQ)_{s_{232}}\right]   \\
        \vspace{2mm}\hspace{5mm}
        \oplus \left[\mr{H}^0(S^{\mr{GL}_2}, \wt{\mathcal{M}}_{(2,-2)})_{s_{213}} 
        \otimes \mr{H}^0(S^{\mr{Sp}_2}, \wt{\mathcal{M}}_2)_{s_{213}}\right]               & q = 3 \\
        \vspace{2mm}
    \left[\mr{H}^1(S^{\mr{GL}_2}, \wt{\mathcal{M}}_{(4,-2)})_{s_{232}} 
        \otimes \mr{H}^0(S^{\mr{Sp}_2}, \uQ)_{s_{232}}\right]        \\
        \vspace{2mm}\hspace{5mm} 
        \oplus \left[\mr{H}^0(S^{\mr{GL}_2}, \wt{\mathcal{M}}_{(4,-2)})_{s_{232}} 
        \otimes \mr{H}^1(S^{\mr{Sp}_2}, \uQ)_{s_{232}}\right] \\
        \vspace{2mm}\hspace{5mm}
        \oplus \left[\mr{H}^1(S^{\mr{GL}_2}, \wt{\mathcal{M}}_{(2,-2)})_{s_{213}} 
          \otimes \mr{H}^0(S^{\mr{Sp}_2},\wt{\mathcal{M}}_2)_{s_{213}}\right] \\
          \vspace{2mm}\hspace{5mm}
        \oplus \left[\mr{H}^0(S^{\mr{GL}_2}, \wt{\mathcal{M}}_{(2,-2)})_{s_{213}} 
          \otimes \mr{H}^1(S^{\mr{Sp}_2},\wt{\mathcal{M}}_2)_{s_{213}}\right] \\
          \vspace{2mm}\hspace{5mm} 
        \oplus \left[\mr{H}^0(S^{\mr{GL}_2}, \wt{\mathcal{M}}_{(2,-3)})_{s_{2132}} 
          \otimes \mr{H}^0(S^{\mr{Sp}_2}, \wt{\mathcal{M}}_{2})_{s_{2132}}\right] \\
          \vspace{2mm}\hspace{5mm} 
        \oplus \left[\mr{H}^0(S^{\mr{GL}_2}, \wt{\mathcal{M}}_{(4,-3)})_{s_{2321}} 
          \otimes \mr{H}^0(S^{\mr{Sp}_2}, \uQ)_{s_{2321}}\right]                        & q = 4 \\
        \end{cases} \\
        &
          \begin{cases}
            \vspace{2mm}
    \left[\mr{H}^1(S^{\mr{GL}_2}, \wt{\mathcal{M}}_{(4,-2)})_{s_{232}} 
        \otimes \mr{H}^1(S^{\mr{Sp}_2}, \uQ)_{s_{232}}\right] \\
        \vspace{2mm}\hspace{5mm}
        \oplus \left[\mr{H}^1(S^{\mr{GL}_2}, \wt{\mathcal{M}}_{(2,-2)})_{s_{213}} 
          \otimes \mr{H}^1(S^{\mr{Sp}_2},\wt{\mathcal{M}}_2)_{s_{213}}\right] \\
          \vspace{2mm}\hspace{5mm} 
        \oplus \left[\mr{H}^1(S^{\mr{GL}_2}, \wt{\mathcal{M}}_{(2,-3)})_{s_{2132}} 
          \otimes \mr{H}^0(S^{\mr{Sp}_2}, \wt{\mathcal{M}}_{2})_{s_{2132}}\right] \\
          \vspace{2mm}\hspace{5mm} 
        \oplus \left[\mr{H}^0(S^{\mr{GL}_2}, \wt{\mathcal{M}}_{(2,-3)})_{s_{2132}} 
          \otimes \mr{H}^1(S^{\mr{Sp}_2}, \wt{\mathcal{M}}_{2})_{s_{2132}}\right] \\
          \vspace{2mm}\hspace{5mm} 
        \oplus \left[\mr{H}^1(S^{\mr{GL}_2}, \wt{\mathcal{M}}_{(4,-3)})_{s_{2321}} 
          \otimes \mr{H}^0(S^{\mr{Sp}_2}, \uQ)_{s_{2321}}\right] \\
          \vspace{2mm}\hspace{5mm} 
        \oplus \left[\mr{H}^0(S^{\mr{GL}_2}, \wt{\mathcal{M}}_{(4,-3)})_{s_{2321}} 
          \otimes \mr{H}^1(S^{\mr{Sp}_2}, \uQ)_{s_{2321}}\right] \\
          \vspace{2mm}\hspace{5mm} 
      \oplus \left[\mr{H}^0(S^{\mr{GL}_2}, \wt{\mathcal{M}}_{(0,-4)})_{s_{21323}} 
        \otimes \mr{H}^0(S^{\mr{Sp}_2}, \wt{\mathcal{M}}_{2})_{s_{21323}}\right]        & q = 5 \\
        \vspace{2mm}
    \left[\mr{H}^1(S^{\mr{GL}_2}, \wt{\mathcal{M}}_{(2,-3)})_{s_{2132}} 
        \otimes \mr{H}^1(S^{\mr{Sp}_2}, \wt{\mathcal{M}}_{2})_{s_{2132}}\right]  \\
        \vspace{2mm}\hspace{5mm} 
        \oplus \left[\mr{H}^1(S^{\mr{GL}_2}, \wt{\mathcal{M}}_{(4,-3)})_{s_{2321}} 
          \otimes \mr{H}^1(S^{\mr{Sp}_2}, \uQ)_{s_{2321}}\right] \\
          \vspace{2mm}\hspace{5mm} 
      \oplus \left[\mr{H}^1(S^{\mr{GL}_2}, \wt{\mathcal{M}}_{(0,-4)})_{s_{21323}} 
        \otimes \mr{H}^0(S^{\mr{Sp}_2}, \wt{\mathcal{M}}_{2})_{s_{21323}}\right] \\
        \vspace{2mm}\hspace{5mm} 
      \oplus \left[\mr{H}^0(S^{\mr{GL}_2}, \wt{\mathcal{M}}_{(0,-4)})_{s_{21323}} 
        \otimes \mr{H}^1(S^{\mr{Sp}_2}, \wt{\mathcal{M}}_{2})_{s_{21323}}\right]        & q = 6 \\
        \vspace{2mm}
    \mr{H}^1(S^{\mr{GL}_2}, \wt{\mathcal{M}}_{(0,-4)})_{s_{21323}} 
        \otimes \mr{H}^1(S^{\mr{Sp}_2}, \wt{\mathcal{M}}_{2})_{s_{21323}}               & q = 7 \\
        \vspace{2mm}
    0 & \text{otherwise}
   \end{cases}
\end{align*}
where we use the K\"unneth Theorem
\[
  \mr{H}^q(S^{\mr{GL}_2} \times S^{\mr{Sp}_2}, \wt{\mathcal{M}_{(a,b)}}) \cong \bigoplus_{n+m = q} \mr{H}^n(S^{\mr{GL}_2}, \wt{\mathcal{M}_{a}}) \otimes \mr{H}^m(S^{\mr{Sp}_4}, \wt{\mathcal{M}_{b}})
\]
and $\mr{H}^q(S^{\mr{GL}_2 \times \mr{Sp}_2}, \wt{\mathcal{M}}) \cong \mr{H}^q(S^{\mr{GL}_2} \times S^{\mr{Sp}_2}, \wt{\mathcal{M}})$; 
Indeed, 
although the locally symmetric space associated with a product of groups does not necessarily decompose into a product of locally symmetric spaces in a strict sense, 
the corresponding arithmetic subgroups are commensurable to the product of arithmetic subgroups of each factor. 
Since we are considering cohomology with $\Q$-coefficients, 
which is invariant under commensurability, 
the decomposition holds.
We have the following fact that
\begin{align*}
  \mr{H}^q(S^{\mr{GL}_2}, \wt{\mathcal{M}}_{(0,l)}) &= 
   \begin{cases}
    \Q & q = 0  \text{ and } l \ \text{ is even} \\
    0 & \text{otherwise}
   \end{cases} \\
   \mr{H}^q(S^{\mr{Sp}_2}, \uQ) &=
\begin{cases}
  \Q & q = 0 \\
  0 & \text{otherwise}
\end{cases}
\end{align*}
In the case $k \neq 0$,
\begin{align*}
\mr{H}^q(S^{\mr{GL}_2}, \wt{\mathcal{M}}_{(k,l)}) \otimes \C &\cong
  \begin{cases}
     \mathcal{S}_{k+2} \oplus \mathcal{E}_{k+2} & q = 1  \text{ and }\ \frac{k}{2} \not \equiv l \pmod 2 \\
     \mathcal{S}_{k+2} & q = 1  \text{ and } \frac{k}{2} \equiv l \pmod 2 \\
    0 & \text{otherwise}
  \end{cases} \\
\mr{H}^q(S^{\mr{SL}_2}, \wt{\mathcal{M}}_{k}) \otimes \C &= \mr{H}^q(S^{\mr{Sp}_2}, \wt{\mathcal{M}}_{k}) \otimes \C \cong
  \begin{cases}
     \mathcal{S}_{k+2} \oplus \overline{\mathcal{S}_{k+2}} \oplus \mathcal{E}_{k+2} & q = 1 \\
    0 & \text{otherwise}
  \end{cases}
\end{align*}
The last isomorphism is called the Eichler-Shimura isomorphism.

Using this fact, we get
\begin{align*}
  \mr{H}^q(\partial_{P_{\{\al_2\}}}, \uQ)
   =
   \begin{cases}
    \Q_e & q = 0 \\
     \vspace{2mm}
    \mathcal{S}_{6, \Q, s_{232}} & q = 4 \\
     \vspace{2mm}
    \left[(\mathcal{S}_{4, \Q} \oplus \mathcal{E}_{4, \Q}) \otimes (\mathcal{S}_{4, \Q} \oplus \overline{\mathcal{S}_{4, \Q}} \oplus \mathcal{E}_{4, \Q}) \right]_{s_{213}} \\
      \vspace{2mm} \hspace{5mm}
      \oplus (\mathcal{S}_{6, \Q} \oplus \mathcal{E}_{6, \Q})_{s_{2321}} & q = 5 \\
    \left[\mathcal{S}_{4, \Q} \otimes (\mathcal{S}_{4, \Q} \oplus \overline{\mathcal{S}_{4, \Q}} \oplus \mathcal{E}_{4, \Q})\right]_{s_{2132}} \\
      \vspace{2mm} \hspace{5mm}
      \oplus (\mathcal{S}_{4, \Q} \oplus \overline{\mathcal{S}_{4, \Q}} \oplus \mathcal{E}_{4, \Q})_{s_{21323}} & q = 6 \\
      \vspace{2mm}
    0 & \text{otherwise}
   \end{cases} \\
\end{align*}
where $\mathcal{S}_{k, \Q}$ denotes the space over $\Q$ such that $\mathcal{S}_{k, \Q} \otimes \C = \mathcal{S}_k$.
It follows from the fact that $\mathcal{S}_{4} = \mathcal{S}_6 = 0, \mathcal{E}_{4,\Q} = \mathcal{E}_{6,\Q} = \Q$.
\begin{align*}
  \mr{H}^q(\partial_{P_{\{\al_2\}}}, \uQ)
  =
   \begin{cases}
    \Q_e & q = 0 \\
    \Q_{s_{213}} \oplus \Q_{s_{2321}} & q = 5 \\
    \Q_{s_{21323}} & q = 6 \\
    0 & \text{otherwise}
   \end{cases}
  \end{align*}

\paragraph{(3) Case $I = \{\alpha_3\}$:}
In this case, the Levi component is $\mr{M} \cong \mr{GL}_3$.
The set of non-vanishing Kostant representatives is 
\[
  \overline{\mathcal{W}^{P_{\{\al_3\}}}} = \{e, s_3, *s_{32}, *s_{321}, *s_{323}, *s_{3213}, s_{32132}, s_{321323}\}.
\] 
We have $\mr{H}^q(S^{\mr{GL}_3}, \wt{\mathcal{M}}) = 0$ for $q > 3$.
The cohomology of the face is 

\begin{align*}
  \mr{H}^{q}(\partial_{\mr{P_{\{\al_3\}}}}, \uQ) 
  &=
   \begin{cases}
    \vspace{2mm}
    \mr{H}^0(S^{\mr{GL}_3}, \uQ)_e & q = 0 \\
      \vspace{2mm}
    \mr{H}^1(S^{\mr{GL}_3}, \uQ)_e
       \oplus \mr{H}^0(S^{\mr{GL}_3}, \wt{\mathcal{M}}_{(0,2,-2)})_{s_{3}} & q = 1 \\
       \vspace{2mm}
    \mr{H}^2(S^{\mr{GL}_3}, \uQ)_e 
       \oplus \mr{H}^1(S^{\mr{GL}_3}, \wt{\mathcal{M}}_{(0,2,-2)})_{s_{3}}
       \oplus \mr{H}^0(S^{\mr{GL}_3}, \wt{\mathcal{M}}_{(1,2,-3)})_{s_{32}} & q = 2 \\
       \vspace{2mm}
    \mr{H}^3(S^{\mr{GL}_3}, \uQ)_e 
       \oplus \mr{H}^2(S^{\mr{GL}_3}, \wt{\mathcal{M}}_{(0,2,-2)})_{s_{3}}
       \oplus \mr{H}^1(S^{\mr{GL}_3}, \wt{\mathcal{M}}_{(1,2,-3)})_{s_{32}} \\
       \vspace{2mm} \hspace{5mm} \oplus \mr{H}^0(S^{\mr{GL}_3}, \wt{\mathcal{M}}_{(0,3,-4)})_{s_{321}}
                                 \oplus \mr{H}^0(S^{\mr{GL}_3}, \wt{\mathcal{M}}_{(3,0,-3)})_{s_{323}} & q = 3 \\
                                 \vspace{2mm}
    \mr{H}^3(S^{\mr{GL}_3}, \wt{\mathcal{M}}_{(0,2,-2)})_{s_{3}}
       \oplus \mr{H}^2(S^{\mr{GL}_3}, \wt{\mathcal{M}}_{(1,2,-3)})_{s_{32}} \\
       \vspace{2mm} \hspace{5mm} \oplus \mr{H}^1(S^{\mr{GL}_3}, \wt{\mathcal{M}}_{(0,3,-4)})_{s_{321}}
                                 \oplus \mr{H}^1(S^{\mr{GL}_3}, \wt{\mathcal{M}}_{(3,0,-3)})_{s_{323}} \\
                                 \vspace{2mm} \hspace{5mm} \oplus \mr{H}^0(S^{\mr{GL}_3}, \wt{\mathcal{M}}_{(2,1,-4)})_{s_{3213}} & q = 4 \\
                                 \vspace{2mm}
    \mr{H}^3(S^{\mr{GL}_3}, \wt{\mathcal{M}}_{(1,2,-3)})_{s_{32}}
       \oplus \mr{H}^2(S^{\mr{GL}_3}, \wt{\mathcal{M}}_{(0,3,-4)})_{s_{321}} \\
       \vspace{2mm} \hspace{5mm} \oplus \mr{H}^2(S^{\mr{GL}_3}, \wt{\mathcal{M}}_{(3,0,-3)})_{s_{323}}
                                 \oplus \mr{H}^1(S^{\mr{GL}_3}, \wt{\mathcal{M}}_{(2,1,-4)})_{s_{3213}} \\
                                 \vspace{2mm} \hspace{5mm} \oplus \mr{H}^0(S^{\mr{GL}_3}, \wt{\mathcal{M}}_{(2,0,-4)})_{s_{32132}} & q = 5 \\
                                 \vspace{2mm}
    \mr{H}^3(S^{\mr{GL}_3}, \wt{\mathcal{M}}_{(0,3,-4)})_{s_{321}}
       \oplus \mr{H}^3(S^{\mr{GL}_3}, \wt{\mathcal{M}}_{(3,0,-3)})_{s_{323}} \\
       \vspace{2mm} \hspace{5mm} \oplus \mr{H}^2(S^{\mr{GL}_3}, \wt{\mathcal{M}}_{(2,1,-4)})_{s_{3213}}
                                 \oplus \mr{H}^1(S^{\mr{GL}_3}, \wt{\mathcal{M}}_{(2,0,-4)})_{s_{32132}} \\
                                 \hspace{5mm} \oplus \mr{H}^0(S^{\mr{GL}_3}, \wt{\mathcal{M}}_{(0,0,-4)})_{s_{321323}} & q = 6 \\
                                 \vspace{2mm}
    \mr{H}^3(S^{\mr{GL}_3}, \wt{\mathcal{M}}_{(2,1,-4)})_{s_{3213}}
       \oplus \mr{H}^2(S^{\mr{GL}_3}, \wt{\mathcal{M}}_{(2,0,-4)})_{s_{32132}} \\
       \vspace{2mm} \hspace{5mm} \oplus \mr{H}^1(S^{\mr{GL}_3}, \wt{\mathcal{M}}_{(0,0,-4)})_{s_{321323}} & q = 7 \\
                                 \vspace{2mm}
    \mr{H}^3(S^{\mr{GL}_3}, \wt{\mathcal{M}}_{(2,0,-4)})_{s_{32132}} 
       \oplus \mr{H}^2(S^{\mr{GL}_3}, \wt{\mathcal{M}}_{(0,0,-4)})_{s_{321323}} & q = 8 \\
       \vspace{2mm}
    \mr{H}^3(S^{\mr{GL}_3}, \wt{\mathcal{M}}_{(0,0,-4)})_{s_{321323}}     & q = 9 \\
      \vspace{2mm}
    0 & \text{otherwise}
   \end{cases}
\end{align*}

From the facts that \cite{Mil}
\begin{align*}
&\bullet 
  \mr{H}^q(S^{\mr{GL}_3(\Z)}, \wt{\mathcal{M}_{(a,b,c)}}) = 
  \begin{cases}
    0 & a + 2b + 3c \equiv 1 \bmod 2 \\
    \mr{H}^q(\mr{SL}_3(\Z), \mathcal{M}_{(a,b)}) & a + 2b + 3c \equiv 0 \bmod 2
  \end{cases}\\
&\bullet 
  \mr{H}_{!}^q(S^{\mr{SL}_3}, \uQ_e) = \mr{H}_{cusp}^q(S^{\mr{SL}_3}, \uQ_e) = 0 \hspace{5mm} \text{for all $q$} \\
&\bullet 
  \mr{H}_{Eis}^q(S^{\mr{SL}_3}, \uQ_e) =
   \begin{cases}
    \Q & q = 0 \\
    0 & \text{otherwise}
   \end{cases}\\
&\bullet  
 \mr{H}_{!}^q(S^{\mr{SL}_3}, \wt{\mathcal{M}_{\lambda \not = 0}}) = 0 \hspace{6mm} (\lambda \not = \lambda^*)  \\
& \text{where $\lambda^* = -w_0(\lambda)$ with the longest element $w_0$ in the Weyl group.} \\
& \text{In $\mr{SL}_3$ case, if $\lambda = (a+b)\varepsilon_1 + b\varepsilon_2$, then $\lambda^* = (a+b)\varepsilon_1 + a\varepsilon_2$.} \\  
&\bullet 
  \text{For even $a>0,$} \\
  & \mr{H}_{Eis}^q(S^{\mr{SL}_3}, \wt{\mathcal{M}_{((a,0))}}) = 
  \mr{H}_{Eis}^q(S^{\mr{SL}_3}, \wt{\mathcal{M}_{((0,a))}}) =
  \begin{cases}
    \mathcal{S}_{a+2, \Q} & q = 3 \\
    0 & \text{otherwise}
  \end{cases} \hspace{6mm} \\
&\bullet 
  \text{For odd $b>0,$} \\
  & \mr{H}_{Eis}^q(S^{\mr{SL}_3}, \wt{\mathcal{M}_{((b,0))}}) = 
  \mr{H}_{Eis}^q(S^{\mr{SL}_3}, \wt{\mathcal{M}_{((0,b))}}) =
  \begin{cases}
    \mathcal{S}_{b+3, \Q} \oplus \Q & q = 2 \\
    0 & \text{otherwise}
  \end{cases} \hspace{6mm} \\
&\bullet 
  \text{For even $a>0$ and odd $b>0,$} \\
  & \mr{H}_{Eis}^q(S^{\mr{SL}_3}, \wt{\mathcal{M}_{((a,b))}}) = 
  \mr{H}_{Eis}^q(S^{\mr{SL}_3}, \wt{\mathcal{M}_{((b,a))}}) =
  \begin{cases}
    \mathcal{S}_{a+b+3, \Q} & q = 2 \\
    \mathcal{S}_{a+2, \Q} & q = 3 \\
    0 & \text{otherwise}
  \end{cases} \hspace{6mm} \\  
&\bullet 
  \mathcal{S}_{4} = \mathcal{S}_6 = 0
\end{align*}
we can derive that
\begin{align*}
  \mr{H}^{q}(\partial_{P_{\{\al_3\}}}, \uQ) 
&=
  \begin{cases}
   \Q_e & q = 0 \\
   \mathcal{S}_{4, \Q, s_{3}} \oplus \mathcal{S}_{6, \Q, s_{32}} & q = 4 \\
   \mathcal{S}_{4, \Q, s_{32}} \oplus (\mathcal{S}_{6, \Q} \oplus \Q)_{s_{321}} \\
     \hspace{5mm} \oplus (\mathcal{S}_{6, \Q} \oplus \Q)_{s_{323}} & q = 5 \\ 
   \mathcal{S}_{6, \Q, s_{3213}} \oplus \Q_{s_{321323}} & q = 6 \\
   \mathcal{S}_{4, \Q, s_{3213}} & q = 7 \\ 
   \mathcal{S}_{4, \Q, s_{32132}} & q = 8 \\
   0 & \text{otherwise}
  \end{cases}  \\
&=
  \begin{cases}
    \Q_e & q = 0 \\
    \Q_{s_{321}} \oplus \Q_{s_{323}} & q = 5 \\
    \Q_{s_{321323}} & q = 6 \\
    0 & \text{otherwise}
  \end{cases}
\end{align*}

\subsubsection*{\underline{\bf{Summary of the $E_1^{0,q}$ terms}}}
Summing the contributions from all maximal parabolic subgroups, 
we obtain the $E_1^{0,q}$

\begin{align*}
  E_1^{0,q} = 
   \begin{cases}
    \Q_{\al_1,e} \oplus \Q_{\al_2,e} \oplus \Q_{\al_3,e} & q = 0 \\
    \Q_{\al_1,e} & q = 2 \\
    \Q_{\al_1, s_{12}} \oplus \Q_{\al_1, s_{12321}} \oplus\Q_{\al_2, s_{213}} \oplus \Q_{\al_2, s_{2321}} \oplus \Q_{\al_{3}, s_{321}} \oplus \Q_{\al_3, s_{323}} & q = 5 \\
    \Q_{\al_1, s_{123}} \oplus \Q_{\al_2, s_{21323}} \oplus \Q_{\al_3, s_{321323}} & q = 6 \\
    \Q_{\al_1, s_{12321}} & q = 7 \\ 
    0 & \text{otherwise}
   \end{cases}
\end{align*}
where the subscript $\al_i$ indicates that
 the object is obtained from the cohomology on $\partial_{P_{\{\al_{i}\}}}$, 
 and symbols such as $s_i$ denote elements of the Weyl group used therein.

\subsubsection{$p = 1$}

For $p = 1$, 
the $E_1$-term is the direct sum of the cohomology of the faces corresponding to the rank 2 parabolic subgroups:
\[
  E_{1}^{1,q} = \bigoplus_{i = 1}^{3} \mr{H}^q(\partial_{P_{\pi \backslash \{{\al_i}\}}}, \uQ).
\]

We compute each face $\mr{H}^q(\partial_{P_{\pi \backslash \{{\al_i}\}}}, \uQ)$.

  \paragraph{(1) Case $I = \{\alpha_1, \alpha_2\}$:}
  In this case, the Levi component is $\mr{M} \cong \mr{GL}_1 \times \mr{GL}_1 \times \mr{Sp}_2$.
  The set of non-vanishing Kostant representatives is
  \[
    \overline{\mathcal{W}^{P_{\{\al_1, \al_2\}}}} = \{e, s_{121}, s_{232}, s_{1213}, s_{2132}, s_{12321}, s_{21323}, s_{12132132}\}. 
  \]
    We have $\mr{H}^q(S^{\mr{Sp}_2}, \wt{\mathcal{M}}) = 0$ for $q > 1$.
    The cohomology of the face is
  \begin{align*}
    \mr{H}^q(\partial_{P_{\{\al_1, \al_2\}}}) &=
      \begin{cases}
        \mr{H}^0(S^{\mr{Sp}_2}, \uQ)_e & q = 0 \\
        \mr{H}^1(S^{\mr{Sp}_2}, \uQ)_e & q = 1 \\
        \mr{H}^0(S^{\mr{Sp}_2}, \wt{\mathcal{M}_{2}})_{s_{121}} 
          \oplus \mr{H}^0(S^{\mr{Sp}_2}, \uQ)_{s_{232}} & q = 3 \\
        \mr{H}^1(S^{\mr{Sp}_2}, \wt{\mathcal{M}_{2}})_{s_{121}}
          \oplus \mr{H}^1(S^{\mr{Sp}_2}, \uQ)_{s_{232}} \\
          \hspace{5mm} \oplus \mr{H}^0(S^{\mr{Sp}_2}, \wt{\mathcal{M}_{2}})_{s_{1213}}
          \oplus \mr{H}^0(S^{\mr{Sp}_2}, \wt{\mathcal{M}_{2}})_{s_{2132}} & q = 4 \\
        \mr{H}^1(S^{\mr{Sp}_2}, \wt{\mathcal{M}_{2}})_{s_{1213}}
          \oplus \mr{H}^1(S^{\mr{Sp}_2}, \wt{\mathcal{M}_{2}})_{s_{2132}} \\
          \hspace{5mm} \oplus \mr{H}^0(S^{\mr{Sp}_2}, \uQ)_{s_{12321}}
          \oplus \mr{H}^0(S^{\mr{Sp}_2}, \wt{\mathcal{M}_{2}})_{s_{21323}} & q = 5 \\
        \mr{H}^1(S^{\mr{Sp}_2}, \uQ)_{s_{12321}}
          \oplus \mr{H}^1(S^{\mr{Sp}_2}, \wt{\mathcal{M}_{2}})_{s_{21323}} & q = 6 \\
        \mr{H}^0(S^{\mr{Sp}_2}, \uQ)_{s_{12132132}} & q = 8 \\
        \mr{H}^1(S^{\mr{Sp}_2}, \uQ)_{s_{12132132}} & q = 9 \\
        0 & \text{otherwise}
      \end{cases} \\
      &=
      \begin{cases}
        \Q_e & q = 0 \\
        \Q_{s_{232}} & q = 3 \\
        (\mathcal{S}_{4, \Q} \oplus \overline{\mathcal{S}_{4, \Q}} \oplus \mathcal{E}_{4, \Q})_{s_{121}} & q = 4 \\
        (\mathcal{S}_{4, \Q} \oplus \overline{\mathcal{S}_{4, \Q}} \oplus \mathcal{E}_{4, \Q})_{s_{1213}} 
          \oplus (\mathcal{S}_{4, \Q} \oplus \overline{\mathcal{S}_{4, \Q}} \oplus \mathcal{E}_{4, \Q})_{s_{2132}}
          \oplus \Q_{s_{12321}} & q = 5 \\
        (\mathcal{S}_{4, \Q} \oplus \overline{\mathcal{S}_{4, \Q}} \oplus \mathcal{E}_{4, \Q})_{s_{21323}} & q = 6 \\
        \Q_{s_{12132132}} & q = 8 \\
        0 & \text{otherwise}
      \end{cases} \\
      &=
      \begin{cases}
        \Q_e & q = 0 \\
        \Q_{s_{232}} & q = 3 \\
        \Q_{s_{121}} & q = 4 \\
        \Q_{s_{1213}} \oplus \Q_{s_{2132}} \oplus \Q_{s_{12321}} & q = 5 \\
        \Q_{s_{21323}} & q = 6 \\
        \Q_{s_{12132132}} & q = 8 \\
        0 & \text{otherwise}
      \end{cases}
  \end{align*}

  \paragraph{(2) Case $I = \{\alpha_1, \alpha_3\}$:}
  In this case, the Levi component is $\mr{M} \cong \mr{GL}_1 \times \mr{GL}_2$.
  The set of non-vanishing Kostant representatives is 
  \[
    \overline{\mathcal{W}^{P_{\{\al_1, \al_3\}}}} = \{e, s_3, *s_{32}, *s_{132}, s_{1321}, *s_{1323}, s_{12321}, s_{13213}, s_{32132}, s_{123213}, s_{321323}, *s_{1232132}\}.
  \]
    We have $\mr{H}^q(S^{\mr{SL}_2}, \wt{\mathcal{M}}) = 0$ for $q > 1$.
 
  \begin{align*}
    \mr{H}^q(\partial_{P_{\{\al_1, \al_3\}}}) &= 
      \begin{cases}
        \vspace{2mm}
        \mr{H}^0(S^{\mr{GL}_2}, \uQ)_e & q = 0 \\
        \vspace{2mm}
        \mr{H}^1(S^{\mr{GL}_2}, \uQ)_e 
          \oplus \mr{H}^0(S^{\mr{GL}_2}, \wt{\mathcal{M}_{(2,-1)}})_{s_{3}} & q = 1 \\
          \vspace{2mm}
        \mr{H}^1(S^{\mr{GL}_2}, \wt{\mathcal{M}_{(2,-1)}})_{s_{3}}
          \oplus \mr{H}^0(S^{\mr{GL}_2}, \wt{\mathcal{M}_{(2,-2)}})_{s_{32}} & q = 2 \\
          \vspace{2mm}
        \mr{H}^1(S^{\mr{GL}_2}, \wt{\mathcal{M}_{(2,-2)}})_{s_{32}} 
          \oplus \mr{H}^0(S^{\mr{GL}_2}, \wt{\mathcal{M}_{(4,-1)}})_{s_{132}} & q = 3 \\
          \vspace{2mm}
        \mr{H}^1(S^{\mr{GL}_2}, \wt{\mathcal{M}_{(4,-1)}})_{s_{132}}
          \oplus \mr{H}^0(S^{\mr{GL}_2}, \wt{\mathcal{M}_{(4,2)}})_{s_{1321}} \\
          \vspace{2mm} \hspace{5mm}
          \oplus \mr{H}^0(S^{\mr{GL}_2}, \wt{\mathcal{M}_{(4,-1)}})_{s_{1323}} & q = 4 \\
          \vspace{2mm}
        \mr{H}^1(S^{\mr{GL}_2}, \wt{\mathcal{M}_{(4,2)}})_{s_{1321}} 
          \oplus \mr{H}^1(S^{\mr{GL}_2}, \wt{\mathcal{M}_{(4,-1)}})_{s_{1323}} \\
          \vspace{2mm} \hspace{5mm}
          \oplus \mr{H}^0(S^{\mr{GL}_2}, \uQ)_{s_{12321}} 
          \oplus \mr{H}^0(S^{\mr{GL}_2}, \wt{\mathcal{M}_{(4,-2)}})_{s_{13213}} \\
          \vspace{2mm} \hspace{5mm}
          \oplus \mr{H}^0(S^{\mr{GL}_2}, \wt{\mathcal{M}_{(0,-4)}})_{s_{32132}} & q = 5 \\
          \vspace{2mm}
        \mr{H}^1(S^{\mr{GL}_2}, \uQ)_{s_{12321}}
          \oplus \mr{H}^1(S^{\mr{GL}_2}, \wt{\mathcal{M}_{(4,-2)}})_{s_{13213}} \\
          \vspace{2mm} \hspace{5mm} 
          \oplus \mr{H}^1(S^{\mr{GL}_2}, \wt{\mathcal{M}_{(0,-4)}})_{s_{32132}}
          \oplus \mr{H}^0(S^{\mr{GL}_2}, \wt{\mathcal{M}_{(2,-1)}})_{s_{123213}} \\
          \vspace{2mm} \hspace{5mm}
          \oplus \mr{H}^0(S^{\mr{GL}_2}, \wt{\mathcal{M}_{(0,-4)}})_{s_{321323}} & q = 6 \\
          \vspace{2mm}
        \mr{H}^1(S^{\mr{GL}_2}, \wt{\mathcal{M}_{(2,-1)}})_{s_{123213}}
          \oplus \mr{H}^1(S^{\mr{GL}_2}, \wt{\mathcal{M}_{(0,-4)}})_{s_{321323}} \\
          \vspace{2mm} \hspace{5mm}
          \oplus \mr{H}^0(S^{\mr{GL}_2}, \wt{\mathcal{M}_{(2,-2)}})_{s_{1232132}} & q = 7 \\
          \vspace{2mm}
          \mr{H}^1(S^{\mr{GL}_2}, \wt{\mathcal{M}_{(2,-2)}})_{s_{1232132}} & q = 8 \\
        0 & \text{otherwise}
      \end{cases} \\
      &=
      \begin{cases}
        \Q_e & q = 0 \\
        \mathcal{S}_{4,\Q, s_{3}} & q = 2 \\
        (\mathcal{S}_{4, \Q} \oplus \mathcal{E}_{4, \Q})_{s_{32}} & q = 3 \\
        (\mathcal{S}_{6, \Q} \oplus \mathcal{E}_{6, \Q})_{s_{132}} & q = 4 \\
        \mathcal{S}_{6,\Q,s_{1321}} 
          \oplus (\mathcal{S}_{6, \Q} \oplus \mathcal{E}_{6, \Q})_{s_{1323}}
          \oplus \Q_{s_{12321}}
          \oplus \Q_{s_{32132}} & q = 5 \\
        \mathcal{S}_{6,\Q, s_{13213}}
          \oplus \Q_{s_{321323}} & q = 6 \\
        \mathcal{S}_{4,\Q, s_{123213}} & q = 7 \\
        (\mathcal{S}_{4, \Q} \oplus \mathcal{E}_{4, \Q})_{s_{1232132}} & q = 8 \\
        0 & \text{otherwise}
      \end{cases} \\
      &=
      \begin{cases}
        \Q_e & q = 0 \\
        \Q_{s_{32}} & q = 3 \\
        \Q_{s_{132}} & q = 4 \\
        \Q_{s_{1323}} \oplus \Q_{s_{12321}} \oplus \Q_{s_{32132}} & q = 5 \\
        \Q_{s_{321323}} & q = 6 \\
        \Q_{s_{1232132}} & q = 8 \\
        0 & \text{otherwise}
      \end{cases}
  \end{align*}

  \paragraph{(3) Case $I = \{\alpha_2, \alpha_3\}$:}
  In this case, the Levi component is $\mr{M} \cong \mr{GL}_2 \times \mr{GL}_1$.
  The set of non-vanishing Kostant representatives is 
  \[
    \overline{\mathcal{W}^{P_{\{\al_2, \al_3\}}}} = \{e, s_3, *s_{213}, s_{232}, s_{2132}, *s_{2321}, s_{2323}, *s_{3213}, s_{21323}, *s_{23213}, s_{32132}, s_{321323}\}.
  \]
    Following a similar argument to Case (2), 
    the cohomology of the face is
  \begin{align*}
    \mr{H}^q(\partial_{P_{\{\al_2, \al_3\}}}) &=
      \begin{cases}
        \vspace{2mm}
        \mr{H}^0(S^{\mr{GL}_2}, \uQ)_e & q = 0 \\
        \vspace{2mm}
        \mr{H}^1(S^{\mr{GL}_2}, \uQ)_e 
          \oplus \mr{H}^0(S^{\mr{GL}_2}, \uQ)_{s_{3}} & q = 1 \\
          \vspace{2mm}
        \mr{H}^1(S^{\mr{GL}_2}, \uQ)_{s_{3}} & q = 2 \\
          \vspace{2mm}
        \mr{H}^0(S^{\mr{GL}_2}, \wt{\mathcal{M}}_{(2,-2)})_{s_{213}}
          \oplus \mr{H}^0(S^{\mr{GL}_2}, \wt{\mathcal{M}}_{(4,-2)})_{s_{232}} & q = 3 \\
          \vspace{2mm}
        \mr{H}^1(S^{\mr{GL}_2}, \wt{\mathcal{M}}_{(2,-2)})_{s_{213}}
          \oplus \mr{H}^1(S^{\mr{GL}_2}, \wt{\mathcal{M}}_{(4,-2)})_{s_{232}} \\
          \vspace{2mm} \hspace{5mm}
          \oplus \mr{H}^0(S^{\mr{GL}_2}, \wt{\mathcal{M}}_{(2,-3)})_{s_{2132}}
          \oplus \mr{H}^0(S^{\mr{GL}_2}, \wt{\mathcal{M}}_{(4,-3)})_{s_{2321}} \\
          \vspace{2mm} \hspace{5mm} 
          \oplus \mr{H}^0(S^{\mr{GL}_2}, \wt{\mathcal{M}}_{(4,-2)})_{s_{2323}}
          \oplus \mr{H}^0(S^{\mr{GL}_2}, \wt{\mathcal{M}}_{(2,-2)})_{s_{3213}} & q = 4 \\
          \vspace{2mm}
        \mr{H}^1(S^{\mr{GL}_2}, \wt{\mathcal{M}}_{(2,-3)})_{s_{2132}} 
          \oplus \mr{H}^1(S^{\mr{GL}_2}, \wt{\mathcal{M}}_{(4,-3)})_{s_{2321}} \\
          \vspace{2mm} \hspace{5mm}
          \oplus \mr{H}^1(S^{\mr{GL}_2}, \wt{\mathcal{M}}_{(4,-2)})_{s_{2323}}
          \oplus \mr{H}^1(S^{\mr{GL}_2}, \wt{\mathcal{M}}_{(2,-2)})_{s_{3213}} \\
          \vspace{2mm} \hspace{5mm} 
          \oplus \mr{H}^0(S^{\mr{GL}_2}, \wt{\mathcal{M}}_{(0,-4)})_{s_{21323}}
          \oplus \mr{H}^0(S^{\mr{GL}_2}, \wt{\mathcal{M}}_{(4,-3)})_{s_{23213}} \\
          \vspace{2mm} \hspace{5mm}
          \oplus \mr{H}^0(S^{\mr{GL}_2}, \wt{\mathcal{M}}_{(2,-3)})_{s_{32132}} & q = 5 \\
          \vspace{2mm}
        \mr{H}^1(S^{\mr{GL}_2}, \wt{\mathcal{M}}_{(0,-4)})_{s_{21323}}
          \oplus \mr{H}^1(S^{\mr{GL}_2}, \wt{\mathcal{M}}_{(4,-3)})_{s_{23213}} \\
          \vspace{2mm} \hspace{5mm}
          \oplus \mr{H}^1(S^{\mr{GL}_2}, \wt{\mathcal{M}}_{(2,-3)})_{s_{32132}}
          \oplus \mr{H}^0(S^{\mr{GL}_2}, \wt{\mathcal{M}}_{(0,-4)})_{s_{321323}} & q = 6 \\
          \vspace{2mm}
        \mr{H}^1(S^{\mr{GL}_2}, \wt{\mathcal{M}}_{(0,-4)})_{s_{321323}} & q = 7 \\
        0 & \text{otherwise}
      \end{cases} \\
      &=
      \begin{cases}
        \Q_e & q = 0 \\
        \Q_{s_{3}} & q = 1 \\
        (\mathcal{S}_{4,\Q} \oplus \mathcal{E}_{4, \Q})_{s_{213}} \oplus \mathcal{S}_{6,\Q, s_{232}} & q = 4 \\
        \mathcal{S}_{4,\Q, s_{2132}}
          \oplus (\mathcal{S}_{6, \Q} \oplus \mathcal{E}_{6, \Q})_{s_{2321}} \\
          \hspace{5mm}
          \oplus \mathcal{S}_{6,\Q, s_{2323}}
          \oplus (\mathcal{S}_{4, \Q} \oplus \mathcal{E}_{4, \Q})_{s_{3213}}
          \oplus \Q_{s_{21323}} & q = 5 \\
        (\mathcal{S}_{6, \Q} \oplus \mathcal{E}_{6, \Q})_{s_{23213}}
          \oplus \mathcal{S}_{4,\Q, s_{32132}}
          \oplus \Q_{s_{321323}} & q = 6 \\
        0 & \text{otherwise}
      \end{cases} \\
      &=
      \begin{cases}
        \Q_e & q = 0 \\
        \Q_{s_{3}} & q = 1 \\
        \Q_{s_{213}} & q = 4 \\
        \Q_{s_{2321}} \oplus \Q_{s_{3213}} \oplus \Q_{s_{21323}} & q = 5 \\
        \Q_{s_{23213}} \oplus \Q_{s_{321323}} & q = 6 \\
        0 & \text{otherwise}
      \end{cases}
  \end{align*}

\subsubsection{\underline{\bf{Summary of the $E_{1}^{1,q}$ terms}}}
Collecting the results from all rank 2 parabolic subgroups, 
we obtain
        
\begin{align*}
  E_{1}^{1,q} = 
    \begin{cases}
      \Q_{\al_{12},e} 
        \oplus \Q_{\al_{13},e} 
        \oplus \Q_{\al_{23},e}  & q = 0 \\
      \Q_{\al_{23},s_{3}} & q = 1 \\
      \Q_{\al_{12},s_{232}} \oplus \Q_{\al_{13}, s_{32}} & q = 3 \\
      \Q_{\al_{12},s_{121}} \oplus \Q_{\al_{13}, s_{132}} \oplus \Q_{\al_{23}, s_{213}} & q = 4 \\
      \Q_{\al_{12},s_{1213}} 
        \oplus \Q_{\al_{12},s_{2132}} 
        \oplus \Q_{\al_{12},s_{12321}} \\
        \hspace{5mm}
        \oplus \Q_{\al_{13}, s_{1323}}
        \oplus \Q_{\al_{13},s_{12321}}
        \oplus \Q_{\al_{13},s_{32132}} \\
        \hspace{5mm}
        \oplus \Q_{\al_{23}, s_{2321}} 
        \oplus \Q_{\al_{23}, s_{3213}}
        \oplus \Q_{\al_{23},s_{21323}} & q = 5 \\
      \Q_{\al_{12},s_{21323}}
        \oplus \Q_{\al_{13},s_{321323}}
        \oplus \Q_{\al_{23}, s_{23213}}
        \oplus \Q_{\al_{23},s_{321323}} & q = 6 \\
      \Q_{\al_{12},s_{12132132}} \oplus \Q_{\al_{13},s_{1232132}} & q = 8 \\
      0 & \text{otherwise}
    \end{cases}
\end{align*}
where the subscript $\al_{i, j}$ indicates that the object is obtained from the cohomology on $\partial_{P_{\{\al_{i}, \al_{j}\}}}$.

\subsubsection{$p = 2$}

For $p = 2$, the cohomology of the face is
\[
  E_{1}^{2,q} = \mr{H}^q(\partial_B, \uQ) = \mr{H}^q(\partial_\pi, \uQ).
\]

In this case, the Levi component is the maximal $\Q$-split torus $\mr{T}$ of $\mr{Sp}_6$. 
Therefore the associated locally symmetric space $S^{\mr{T}}$ is a finite set, 
and so $\mr{H}^q(S^{\mr{T}}, \wt{\mathcal{M}}) = 0$ for all $q > 0$. 
In particular, since $S^{\mr{T}}$ consists of a single element, we have $\mr{H}^0(S^{\mr{T}}, \uQ) = \Q$.

The set of non-vanishing Kostant representatives is
\begin{align*}
  \overline{\mathcal{W}^{P_{\pi}}} &= \{e, s_{3}, s_{121}, s_{232}, s_{1213}, s_{1321}, s_{2132}, s_{2323}, s_{12321}, s_{13213}, s_{21323}, s_{32132},\\
  & s_{123213}, s_{321323}, s_{12132132}, s_{121321323}\},
\end{align*}
the cohomology of the face is
\begin{equation*}
  E_1^{2,q} = 
    \begin{cases}
      \Q_e & q = 0 \\
      \Q_{s_3} & q = 1 \\
      \Q_{s_{121}} \oplus \Q_{s_{232}} & q = 3 \\
      \Q_{s_{1213}} \oplus \Q_{s_{1321}} \oplus \Q_{s_{2132}} \oplus \Q_{s_{2323}} & q = 4 \\
      \Q_{s_{12321}} \oplus \Q_{s_{13213}} \oplus \Q_{s_{21323}} \oplus \Q_{s_{32132}} & q = 5 \\
      \Q_{s_{123213}} \oplus \Q_{s_{321323}} & q = 6 \\
      \Q_{s_{12132132}} & q = 8 \\
      \Q_{s_{121321323}} & q = 9 \\
      0 & \text{otherwise}
    \end{cases}
\end{equation*}

\begin{figure}[h]
  \centering
  \begin{tikzpicture}[xscale=2.5, yscale=0.9] 
    \draw[->, thin] (-0.05,0) -- (2.6,0) node[right] {$p$};
    \draw[->, thin] (0,-0.05) -- (0,9.5) node[above] {$q$};
  
    \foreach \x in {0,1,2}{
      \draw[thick] (\x,0) -- (\x,-0.2) node[below] {\x};
    }
  
    \foreach \y in {0,1,...,9}{
      \draw[thick] (0,\y) -- (-0.2,\y) node[left] {\y};
    }
  
    \foreach \x in {0,1,2}{
      \foreach \y in {0,1,...,9}{
        \fill[gray!40] (\x,\y) circle (1pt);
      }
    }

    \fill[black] (0,0) circle (4pt);
    \fill[black] (0,2) circle (4pt);
    \fill[black] (0,5) circle (4pt);
    \fill[black] (0,6) circle (4pt);
    \fill[black] (0,7) circle (4pt);

    \fill[black] (1,0) circle (4pt);
    \fill[black] (1,1) circle (4pt);
    \fill[black] (1,3) circle (4pt);
    \fill[black] (1,4) circle (4pt);
    \fill[black] (1,5) circle (4pt);
    \fill[black] (1,6) circle (4pt);
    \fill[black] (1,8) circle (4pt);

    \fill[black] (2,0) circle (4pt);
    \fill[black] (2,1) circle (4pt);
    \fill[black] (2,3) circle (4pt);
    \fill[black] (2,4) circle (4pt);
    \fill[black] (2,5) circle (4pt);
    \fill[black] (2,6) circle (4pt);
    \fill[black] (2,8) circle (4pt);
    \fill[black] (2,9) circle (4pt);
   
    \draw[->, very thick, black] (0.2,0) -- (0.8,0) node[midway, above] {$d_1$};
    \draw[->, very thick, black] (0.2,5) -- (0.8,5) node[midway, above] {$d_1$};
    \draw[->, very thick, black] (0.2,6) -- (0.8,6) node[midway, above] {$d_1$};

    \draw[->, very thick, black] (1.2,0) -- (1.8,0) node[midway, above] {$d_1$};
    \draw[->, very thick, black] (1.2,1) -- (1.8,1) node[midway, above] {$d_1$};
    \draw[->, very thick, black] (1.2,3) -- (1.8,3) node[midway, above] {$d_1$};
    \draw[->, very thick, black] (1.2,4) -- (1.8,4) node[midway, above] {$d_1$};
    \draw[->, very thick, black] (1.2,5) -- (1.8,5) node[midway, above] {$d_1$};
    \draw[->, very thick, black] (1.2,6) -- (1.8,6) node[midway, above] {$d_1$};
    \draw[->, very thick, black] (1.2,8) -- (1.8,8) node[midway, above] {$d_1$};

  \end{tikzpicture}
  \caption{$E_{1}$-page}
  \label{E1}
  \end{figure}
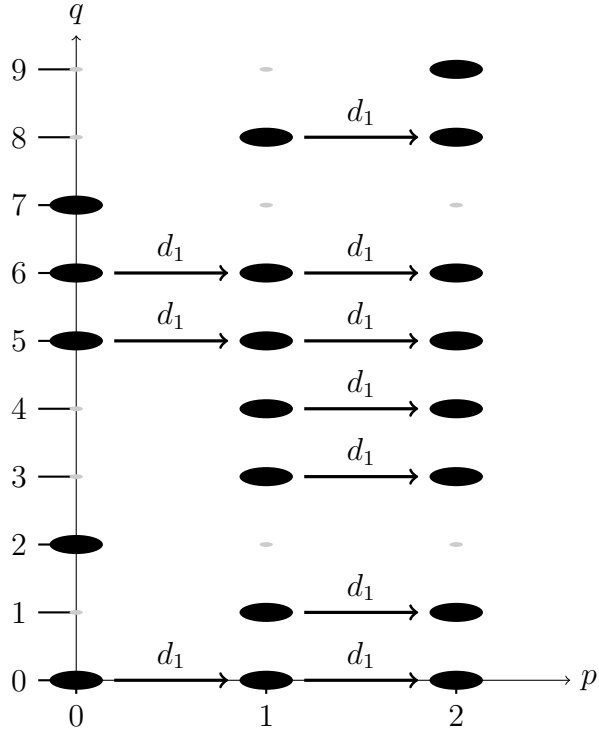

The structure of the $E_{1}$-page is summarized in Figure \ref{E1}.
Each dot represents a position where the cohomology group $E_{1}^{p,q}$ is non-vanishing,
and the arrows indicate the action of the first differentials $d_1^{p,q}$.

\subsection{$E_{2}$-page}

To obtain $E_2$-terms, 
it is necessary to consider the differentials $d_{1}^{p,q} \colon E_{1}^{p,q} \rightarrow E_{1}^{p+1,q}$;

\subsubsection{At the level $q = 0$}

We consider
\[
  0 \rightarrow E_{1}^{0,0} \overset{d_{1}^{0,0}}{\longrightarrow} E_{1}^{1,0} \overset{d_{1}^{1,0}}{\longrightarrow} E_{1}^{2,0} \rightarrow 0
\]
We have 
\begin{align*}
  E_{1}^{0,0} &= \Q_{\al_1,e} \oplus \Q_{\al_2,e} \oplus \Q_{\al_3,e} \\
  E_{1}^{1, 0} &= \Q_{\al_{12}, e} \oplus \Q_{\al_{13}, e} \oplus \Q_{\al_{23}, e} \\ 
  E_{1}^{2,0} &= \Q_e
\end{align*}
The first differential $d_{1}^{0,0} \colon \Q_{\al_1,e} \oplus \Q_{\al_2,e} \oplus \Q_{\al_3,e} \to \Q_{\al_{12}, e} \oplus \Q_{\al_{13}, e} \oplus \Q_{\al_{23}, e}$
is given by 
$(a_1, a_2, a_3) \mapsto (a_1 - a_2, a_1 - a_3, a_2 - a_3)$,
and the second differential
$d_{1}^{1,0} \colon \Q_{\al_{12}, e} \oplus \Q_{\al_{13}, e} \oplus \Q_{\al_{23}} \to \Q_e$
is given by
$(b_1, b_2, b_3) \mapsto (-b_1 + b_2 - b_3)$.
The structure of these differentials is shown in the diagram below.

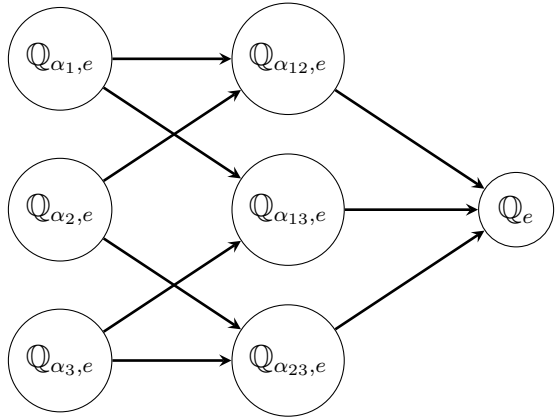
\begin{figure}[ht]

  \begin{tikzpicture}[node distance=1.2cm and 2.5cm, auto]

    \node[circle, draw] (A1) at (0,2) {$\Q_{\al_1,e}$};
    \node[circle, draw] (A2) at (0,0) {$\Q_{\al_2,e}$};
    \node[circle, draw] (A3) at (0,-2) {$\Q_{\al_3,e}$};

    \node[circle, draw] (B1) at (3,2) {$\Q_{\al_{12}, e}$};
    \node[circle, draw] (B2) at (3,0) {$\Q_{\al_{13}, e}$};
    \node[circle, draw] (B3) at (3,-2) {$\Q_{\al_{23}, e}$};

    \node[circle, draw] (C) at (6,0) {$\Q_e$};

    \draw[->, >=stealth, line width=1pt] (A1) -- (B1);
    \draw[->, >=stealth, line width=1pt] (A1) -- (B2);
    \draw[->, >=stealth, line width=1pt] (A2) -- (B1);
    \draw[->, >=stealth, line width=1pt] (A2) -- (B3);
    \draw[->, >=stealth, line width=1pt] (A3) -- (B2);
    \draw[->, >=stealth, line width=1pt] (A3) -- (B3);

    \draw[->, >=stealth, line width=1pt] (B1) -- (C);
    \draw[->, >=stealth, line width=1pt] (B2) -- (C);
    \draw[->, >=stealth, line width=1pt] (B3) -- (C);

  \end{tikzpicture}
\end{figure}

By analyzing the kernels and images of these maps, we obtain
\begin{align*}
  \mr{Ker}(d_{1}^{0,0}) &= \{(a_1, a_2, a_3) \in \Q^3 \mid a_1 = a_2 = a_3 \} = \Q, \\
  \mr{Im}(d_{1}^{0,0}) &= \Q^3 / \mr{Ker}(d_{1}^{0,0}) = \Q^2, \\
  \mr{Im}(d_{1}^{1,0}) &= \Q, \\
  \mr{Ker}(d_{1}^{1,0}) &= \Q^{3 - \mr{dim}_{\Q}(\mr{Im}(d_{1}^{1,0}))} = \Q^2.
\end{align*}
It follows that
\begin{align*}
  E_{2}^{0,0} &= \mr{Ker}(d_{1}^{0,0}) = \Q, \\
  E_{2}^{1,0} &= \mr{Ker}(d_{1}^{1,0}) / \mr{Im}(d_{1}^{0,0}) = \Q^2 / \Q^2 = 0, \\
  E_{2}^{2,0} &= \mr{Coker}(d_{1}^{1,0}) = \Q / \mr{Im}(d_{1}^{1,0}) = \Q / \Q = 0.
\end{align*}

\subsubsection{At the level $q = 1$}

We consider
\[
  0 \to E_{1}^{1,1} \overset{d_{1}^{1,1}}{\longrightarrow} E_{1}^{2,1} \to 0
\]
We have
\begin{align*}
  E_{1}^{1,1} &= \Q_{\al_{23},s_3} \\
  E_{1}^{2,1} &= \Q_{s_{3}}.
\end{align*}
The differential $d_{1}^{1,1}$ is an isomorphism. 
Therefore, we get
\begin{align*}
  E_2^{0,1} = E_2^{1,1} = E_2^{2,1} = 0
\end{align*}

\subsubsection{At the level $q = 2$}

We consider
\[
  0 \to E_{1}^{0,2} \overset{d_{1}^{0,2}}{\longrightarrow} 0
\]
We have
\begin{align*}
  E_{1}^{0,2} = \Q_{\al_{1},e}
\end{align*}
Therefore, we get
\begin{align*}
  E_2^{0,2} = \Q \\
  E_2^{1,2} = E_2^{2,2} = 0
\end{align*}

\subsubsection{At the level $q = 3$}

We consider
\[
  0 \to E_{1}^{1,3} \overset{d_{1}^{1,3}}{\longrightarrow} E_{1}^{2,3} \to 0
\]
We have
\begin{align*}
  E_{1}^{1,3} &= \Q_{\al_{12},s_{232}} \oplus \Q_{\al_{13}, s_{32}} \\
  E_{1}^{2,3} &= \Q_{s_{121}} \oplus \Q_{s_{232}}.
\end{align*}
There is a map 
\[
  \Q_{\al_{12},s_{232}} \oplus \Q_{\al_{13}, s_{32}} \to \Q_{s_{232}}; (a, b)  \mapsto -a + b.
\] 
Both the kernel and the image of this map are $\Q$.
Therefore, we get
\begin{align*}
  E_2^{0,3} &= 0, \\
  E_2^{1,3} &= \Q, \\
  E_2^{2,3} &= \Q.
\end{align*}

\subsubsection{At the level $q = 4$}

We consider
\[
  0 \to E_{1}^{1,4} \overset{d_{1}^{1,4}}{\longrightarrow} E_{1}^{2,4} \to 0
\]
We have
\begin{align*}\
  E_{1}^{1,4} &= \Q_{\al_{12},s_{121}} \oplus \Q_{\al_{13}, s_{132}} \oplus \Q_{\al_{23}, s_{213}} \\
  E_{1}^{2,4} &= \Q_{s_{1213}} \oplus \Q_{s_{1321}} \oplus \Q_{s_{2132}} \oplus \Q_{s_{2323}}.
\end{align*}
The differential $d_1^{1,4}$ consists of the maps
\begin{align*}
 \Q_{\al_{12}, s_{121}} \to \Q_{s_{1321}}, \\
 \Q_{\al_{13}, s_{132}} \to \Q_{s_{2132}}, \\
 \Q_{\al_{23}, s_{213}} \to \Q_{s_{1213}}.
\end{align*} 
Thus, $\mr{Ker}(d_1^{1,4}) = 0$ and $\mr{Im}(d_1^{1,4}) = \Q^3$.
Therefore we get
\begin{align*}
  E_2^{0,4} &= 0, \\
  E_2^{1,4} &= 0, \\
  E_2^{2,4} &= \Q.
\end{align*}

\subsubsection{At the level $q = 5$}

We consider
\[
  0 \to E_{1}^{0,5} \overset{d_{1}^{0,5}}{\longrightarrow} E_{1}^{1,5} \overset{d_{1}^{1,5}}{\longrightarrow} E_{1}^{2,5} \to 0
\]
We have
\begin{align*}
  E_{1}^{0,5} &= \Q_{\al_{1},s_{12}}
   \oplus \Q_{\al_{1},s_{12321}}
   \oplus \Q_{\al_{2},s_{213}}
   \oplus \Q_{\al_{2},s_{2321}}
   \oplus \Q_{\al_{3},s_{321}}
   \oplus \Q_{\al_{3},s_{323}} \\
  E_{1}^{1,5} &= \Q_{\al_{12},s_{1213}} 
   \oplus \Q_{\al_{12},s_{2132}} 
   \oplus \Q_{\al_{12},s_{12321}} \\
   & \hspace{5mm} 
   \oplus \Q_{\al_{13},s_{1323}}
   \oplus \Q_{\al_{13},s_{12321}}
   \oplus \Q_{\al_{13},s_{32132}} \\
   & \hspace{5mm}
   \oplus \Q_{\al_{23},s_{2321}}
   \oplus \Q_{\al_{23},s_{3213}}
   \oplus \Q_{\al_{23},s_{21323}} \\
  E_{1}^{2,5} &= \Q_{s_{12321}} 
   \oplus \Q_{s_{13213}} 
   \oplus \Q_{s_{21323}} 
   \oplus \Q_{s_{32132}}.
\end{align*}

First, the differential $d_1^{0,5}$ consists of
\begin{align*}
  &\Q_{\al_{1},s_{12}}
  \oplus \Q_{\al_{1},s_{12321}}
  \oplus \Q_{\al_{2},s_{213}}
  \oplus \Q_{\al_{2},s_{2321}} 
  \oplus \Q_{\al_{3},s_{321}}
  \oplus \Q_{\al_{3},s_{323}} \\
  &\to
   \Q_{\al_{12},s_{1213}} 
  \oplus \Q_{\al_{12},s_{2132}} 
  \oplus \Q_{\al_{12},s_{12321}} \\
  & \hspace{5mm} 
  \oplus \Q_{\al_{13},s_{1323}}
  \oplus \Q_{\al_{13},s_{12321}}
  \oplus \Q_{\al_{13},s_{32132}} \\
  & \hspace{5mm}
  \oplus \Q_{\al_{23},s_{2321}}
  \oplus \Q_{\al_{23},s_{3213}}
  \oplus \Q_{\al_{23},s_{21323}}: \\
  &(a,b,c,d,e,f) \\
  &\mapsto (-c, a, b-d, -f, b-e, a, d-e,c,-f)
\end{align*}
Thus, $\mr{Ker}(d_1^{0,5}) = \Q$ and $\mr{Im}(d_1^{0,5}) = \Q^5$.

Next, the differential $d_1^{1,5}$ consists of
\begin{align*}
  &\Q_{\al_{12},s_{1213}} 
  \oplus \Q_{\al_{12},s_{2132}} 
  \oplus \Q_{\al_{12},s_{12321}} \\
  & \hspace{5mm} 
  \oplus \Q_{\al_{13},s_{1323}}
  \oplus \Q_{\al_{13},s_{12321}}
  \oplus \Q_{\al_{13},s_{32132}} \\
  & \hspace{5mm}
  \oplus \Q_{\al_{23},s_{2321}}
  \oplus \Q_{\al_{23},s_{3213}}
  \oplus \Q_{\al_{23},s_{21323}} \\
  &\to
  \Q_{s_{12321}} 
  \oplus \Q_{s_{13213}} 
  \oplus \Q_{s_{21323}} 
  \oplus \Q_{s_{32132}}: \\
  &(a', b', c', d', e', f', g', h', i') \\
  &\mapsto (-c' + e' - g', -a' - h', d'-i', -b' + f').
\end{align*}
Thus, $\mr{Ker}(d_1^{1,5}) = \Q^5$ and $\mr{Im}(d_1^{1,5}) = \Q^4$.
 Therefore, we obtain
 \begin{align*}
  E_{2}^{0,5} &= \mr{Ker}(d_{1}^{0,5}) = \Q \\
  E_{2}^{1,5} &= \mr{Ker}(d_{1}^{1,5})/\mr{Im}(d_1^{0,5}) = 0 \\
  E_{2}^{2,5} &= \mr{Coker}(d_1^{1,5}) = 0
\end{align*}

\subsubsection{At the level $q = 6$}

We consider
\[
  0 \to E_{1}^{0,6} \overset{d_{1}^{0,6}}{\longrightarrow} E_{1}^{1,6} \overset{d_{1}^{1,6}}{\longrightarrow} E_{1}^{2,6} \to 0
\]
We have
\begin{align*}
  E_{1}^{0,6} &= 
   \Q_{\al_1, s_{123}}
   \oplus \Q_{\al_2, s_{21323}} 
   \oplus \Q_{\al_3, s_{321323}} \\
  E_{1}^{1,6} &= \Q_{\al_{12},s_{21323}}
   \oplus \Q_{\al_{13},s_{321323}}
   \oplus \Q_{\al_{23},s_{23213}}
   \oplus \Q_{\al_{23},s_{321323}}  \\
  E_{1}^{2,6} &= \Q_{s_{123213}} 
   \oplus \Q_{s_{321323}}.
\end{align*}

First, the differential $d_1^{0,6}$ consists of 
\begin{align*}
  &\Q_{\al_1, s_{123}} \oplus \Q_{\al_2, s_{21323}} \oplus \Q_{\al_3, s_{321323}} 
   \to \Q_{\al_{12},s_{21323}} \oplus \Q_{\al_{13},s_{321323}} \oplus \Q_{\al_{23},s_{321323}} \\
  & \hspace{15mm} (a,b,c) \hspace{15mm} \mapsto \hspace{15mm} (a+b, a+c, c-b).
\end{align*}
Thus, we get $\mr{Ker}(d_1^{0,6}) = \Q, \mr{Im}(d_1^{0,6}) = \Q^2$.

The second differential $d_1^{1,6}$ consists of
\begin{align*}
 &f_1: \Q_{\al_{12},s_{21323}} \oplus \Q_{\al_{13},s_{321323}} \oplus \Q_{\al_{23},s_{321323}} \to \Q_{s_{321323}}, (a,b,c) \mapsto -a+b-c, \\
 &f_2: \Q_{\al_{23}, s_{23213}} \to \Q_{s_{123213}}.
\end{align*}
Thus, $\mr{Ker}(d_1^{1,6}) = \Q^2, \mr{Im}(d_1^{1,6}) = \Q^2$.

Therefore, we obtain
\begin{align*}
  E_{2}^{0,6} &= \mr{Ker}(d_{1}^{0,6}) = \Q, \\
  E_{2}^{1,6} &= \mr{Ker}(d_{1}^{1,6})/\mr{Im}(d_1^{0,6}) = 0, \\
  E_{2}^{2,6} &= \mr{Coker}(d_1^{1,6}) = 0.
\end{align*}

\subsubsection{At the level $q = 7$}

We consider
\[
  0 \to E_{1}^{0,7} \overset{d_{1}^{0,7}}{\longrightarrow} 0
\]

Therefore, we get
\begin{align*}
  E_{2}^{0,7} = E_1^{0,7} = \Q \\
  E_{2}^{1,7} = E_{1}^{2,7} = 0
\end{align*}

\subsubsection{At the level $q = 8$}

We consider
\[
  0 \to E_{1}^{1,8} \overset{d_{1}^{1,8}}{\longrightarrow} E_{1}^{2,8} \to 0
\]
We have
\begin{align*}
  E_1^{1,8} &= \Q_{\al_{12},s_{12132132}} \oplus \Q_{\al_{13},s_{1232132}} \\
  E_1^{2,8} &= \Q_{s_{12132132}}
\end{align*}

The differential $d_1^{1,8}$ consists of 
\[ 
 \Q_{\al_{12},s_{12132132}} \oplus \Q_{\al_{13},s_{1232132}} \to \Q_{s_{12132132}}: (a, b) \mapsto -a + b.
\]
Therefore, we get
\begin{align*}
  E_{2}^{0,8} &= 0 \\
  E_{2}^{1,8} &= \Q \\
  E_{2}^{2,8} &= 0.
\end{align*}

\subsubsection{At the level $q = 9$}

We have
\[
  E_2^{p,9} = E_1^{p,9} \ (p = 0,1,2).
\]

\begin{figure}[h]
  \centering
  \begin{tikzpicture}[xscale=2.5, yscale=0.9] 
    \draw[->, thin] (-0.05,0) -- (2.6,0) node[right] {$p$};
    \draw[->, thin] (0,-0.05) -- (0,9.5) node[above] {$q$};
  
    \foreach \x in {0,1,2}{
      \draw[thick] (\x,0) -- (\x,-0.2) node[below] {\x};
    }
  
    \foreach \y in {0,1,...,9}{
      \draw[thick] (0,\y) -- (-0.2,\y) node[left] {\y};
    }
  
    \foreach \x in {0,1,2}{
      \foreach \y in {0,1,...,9}{
        \fill[gray!40] (\x,\y) circle (1pt);
      }
    }

    \fill[black] (0,0) circle (4pt);
    \fill[black] (0,2) circle (4pt);
    \fill[black] (0,5) circle (4pt);
    \fill[black] (0,6) circle (4pt);
    \fill[black] (0,7) circle (4pt);
    \fill[black] (1,3) circle (4pt);
    \fill[black] (2,3) circle (4pt);
    \fill[black] (2,4) circle (4pt);
    \fill[black] (1,8) circle (4pt);
    \fill[black] (2,9) circle (4pt);

    \draw[->, very thick, black] (0.2,5.0) -- (1.8,4.0) node[midway, above] {$d_2$};

  \end{tikzpicture}
  \caption{$E_{2}$-page}
  \label{E2}
  \end{figure}
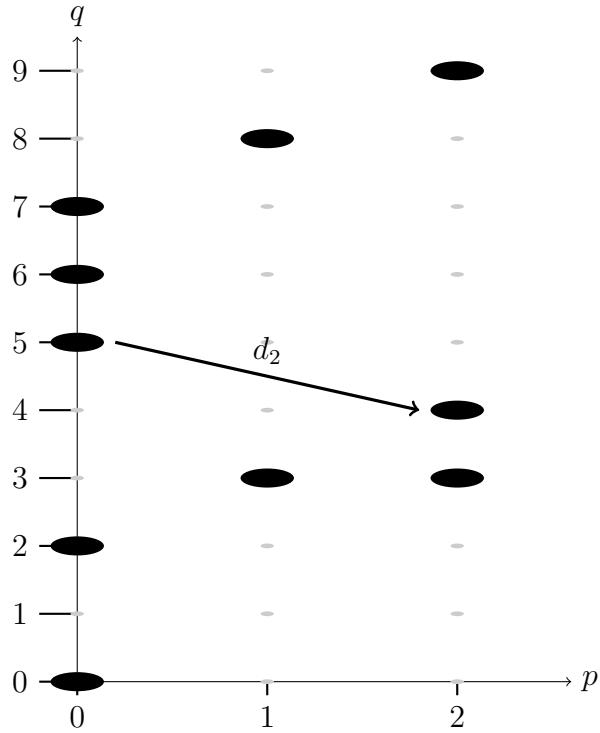

  \subsection{$E_{3}$-page}

  To obtain the $E_3$-page, 
  it is necessary to consider the differential $d_{2}^{p,q} \colon E_{2}^{p,q} \rightarrow E_{2}^{p+2,q-1}$.
  As illustrated in Figure \ref{E2}, 
  the only potentially non-trivial differential occurs when $(p, q) = (0, 5)$.
  For all other $(p, q)$, the differentials vanish, and thus $E_{3}^{p,q} = E_{2}^{p,q}$.

  We consider 
  \[
  0 \to E_{2}^{0,5} \overset{d_{2}^{0,5}}{\longrightarrow} E_{2}^{2,4} \to 0
  \]
  We have 
  \begin{align*}
    E_{2}^{0,5} &= \Q \\
    E_{2}^{2,4} &= \Q_{s_{2323}}.
  \end{align*}
  The differential $d_{2}^{0,5}$ is the zero map 
  since $E_{2}^{0,5}$ is derived from $\Q_{\al_{1}, s_{12321}}$ or $\Q_{\al_{2}, s_{2321}}$ or $\Q_{\al_{3}, s_{321}}$ 
  and none of the elements maps to $s_{2323}$.
  Therefore
  \begin{align*}
    E_{3}^{0,5} &= \mr{Ker}(d_{2}^{0,5}) = \Q, \\
    E_{3}^{2,4} &= \mr{Coker}(d_{2}^{0,5}) = \Q.
  \end{align*}

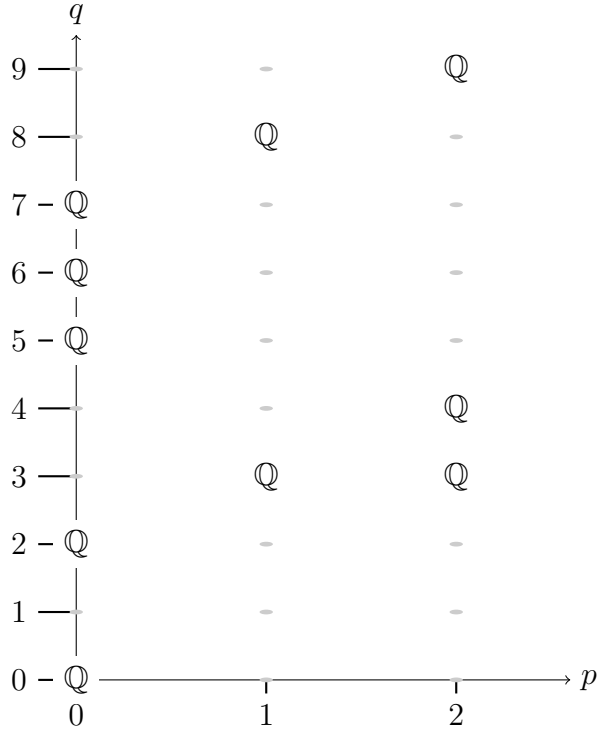
\begin{figure}[h]
  \centering
  \begin{tikzpicture}[xscale=2.5, yscale=0.9] 
    \draw[->, thin] (-0.05,0) -- (2.6,0) node[right] {$p$};
    \draw[->, thin] (0,-0.05) -- (0,9.5) node[above] {$q$};
  
    \foreach \x in {0,1,2}{
      \draw[thick] (\x,0) -- (\x,-0.2) node[below] {\x};
    }
  
    \foreach \y in {0,1,...,9}{
      \draw[thick] (0,\y) -- (-0.2,\y) node[left] {\y};
    }
  
    \foreach \x in {0,1,2}{
      \foreach \y in {0,1,...,9}{
        \fill[gray!40] (\x,\y) circle (1pt);
      }
    }

    \node[fill=white] at (0,0) {$\Q$};
    \node[fill=white] at (0,2) {$\Q$};
    \node[fill=white] at (0,5) {$\Q$};
    \node[fill=white] at (0,6) {$\Q$};
    \node[fill=white] at (0,7) {$\Q$};
    \node[fill=white] at (1,3) {$\Q$};
    \node[fill=white] at (1,8) {$\Q$};
    \node[fill=white] at (2,3) {$\Q$};
    \node[fill=white] at (2,4) {$\Q$};
    \node[fill=white] at (2,9) {$\Q$};

  \end{tikzpicture}
  \caption{$E_{3}$-page}
  \label{E3}
  \end{figure}

  The summary of the $E_{3}$-page is shown in Figure \ref{E3}.
The only difference compared to the $E_{2}$-page is 
the cancellation of the terms at $(p, q) = (0, 7)$ and $(2, 6)$.

\subsection{Boundary cohomology of $\mr{Sp}_6(\Z)$}

 From the relation
 \[
  \mr{H}^k(\partial S, \uQ) = \bigoplus_{p+q = k} E_{3}^{p,q},
 \]
we obtain the following theorem.

\begin{mainthm}
  The boundary cohomology of the orbifold $S$ of the arithmetic group $\mr{Sp}_{6}(\Z)$ with trivial coefficients is described as follows.
  \begin{equation*}
    \mr{H}^{q}(\partial S_{\Gamma}, \uQ) = 
    \begin{cases}
      \Q & q = 0,2,4,7,9,11 \\
      \Q^2 & q = 5,6 \\
      0 & \text{otherwise}
    \end{cases}
  \end{equation*}
\end{mainthm}

\begin{Rem}\rm
  While the computation is explicit for trivial coefficients, 
  the case of non-trivial coefficients is significantly more involved. 
  This difficulty stems primarily from the limited information currently available on the interior (inner) cohomology of the Levi factors, 
  such as $\mr{SL}_{3}$. 
  Although the Eisenstein cohomology for these groups is well-understood \cite{sl3}, 
  a complete determination of the $E_1$-page for general coefficients would require full knowledge of the interior cohomology, 
  which remains a subject of ongoing research.

 \end{Rem}

\newpage
\appendix
\section{Detailed Structure of Levi Quotients} \label{str of levi}

In this appendix, we provide the diagrammatic representation of the Levi quotients for each standard $\Q$-parabolic subgroup $P_{I}$. The nodes removed from the $\mr{C}_3$ Dynkin diagram are denoted by $\times$.

\tikzset{
  dotnode/.style={shape=circle, draw, fill=black, inner sep=0pt, minimum size=4pt},
  crossnode/.style={shape=cross out, draw, inner sep=0pt, minimum size=4pt}
}

\begin{center}
\begin{tikzpicture}
  \node[dotnode, label=below:$\alpha_1$] (1) at (0,0) {};
  \node[dotnode, label=below:$\alpha_2$] (2) at (1.2,0) {};
  \node[dotnode, label=below:$\alpha_3$] (3) at (2.4,0) {};
  \draw (1) -- (2);
  \draw ([yshift=0.7pt]3.west) -- node[midway, inner sep=0.5pt, draw=none] {$<$} ([yshift=0.7pt]2.east);
  \draw ([yshift=-0.7pt]3.west) -- ([yshift=-0.7pt]2.east);
\end{tikzpicture}
\ ($\mr{C}_3$)
\end{center}

\subsubsection*{Rank 1 ($|I|=1$)}

$\bullet P_{\{\alpha_1\}}$: $\mr{M}_{P_{\{\alpha_1\}}} = \mr{GL}_1 \times \mr{Sp}_4$.
\[
\begin{tikzpicture}
  \node[crossnode, label=below:$\alpha_1$] (1) at (0,0) {};
  \node[dotnode, label=below:$\alpha_2$] (2) at (1.2,0) {};
  \node[dotnode, label=below:$\alpha_3$] (3) at (2.4,0) {};
  \draw ([yshift=0.7pt]3.west) -- node[midway, inner sep=0.5pt, draw=none] {$<$} ([yshift=0.7pt]2.east);
  \draw ([yshift=-0.7pt]3.west) -- ([yshift=-0.7pt]2.east);
\end{tikzpicture}
\]

$\bullet P_{\{\alpha_2\}}$: $\mr{M}_{P_{\{\alpha_2\}}} = \mr{GL}_2 \times \mr{Sp}_2$.
\[
\begin{tikzpicture}
  \node[dotnode, label=below:$\alpha_1$] (1) at (0,0) {};
  \node[crossnode, label=below:$\alpha_2$] (2) at (1.2,0) {};
  \node[dotnode, label=below:$\alpha_3$] (3) at (2.4,0) {};
\end{tikzpicture}
\]

$\bullet P_{\{\alpha_3\}}$: $\mr{M}_{P_{\{\alpha_3\}}} = \mr{GL}_3$.
\[
\begin{tikzpicture}
  \node[dotnode, label=below:$\alpha_1$] (1) at (0,0) {};
  \node[dotnode, label=below:$\alpha_2$] (2) at (1.2,0) {};
  \node[crossnode, label=below:$\alpha_3$] (3) at (2.4,0) {};
  \draw (1) -- (2);
\end{tikzpicture}
\]

\subsubsection*{Rank 2 ($|I|=2$)}

$\bullet P_{\{\alpha_1, \alpha_2\}}$: $\mr{M}_{P_{\{\alpha_1, \alpha_2\}}} = \mr{GL}_1 \times \mr{GL}_1 \times \mr{Sp}_2$.
\[
\begin{tikzpicture}
  \node[crossnode, label=below:$\alpha_1$] (1) at (0,0) {};
  \node[crossnode, label=below:$\alpha_2$] (2) at (1.2,0) {};
  \node[dotnode, label=below:$\alpha_3$] (3) at (2.4,0) {};
\end{tikzpicture}
\]

$\bullet P_{\{\alpha_1, \alpha_3\}}$: $\mr{M}_{P_{\{\alpha_1, \alpha_3\}}} = \mr{GL}_1 \times \mr{GL}_2$. 
\[
\begin{tikzpicture}
  \node[crossnode, label=below:$\alpha_1$] (1) at (0,0) {};
  \node[dotnode, label=below:$\alpha_2$] (2) at (1.2,0) {};
  \node[crossnode, label=below:$\alpha_3$] (3) at (2.4,0) {};
\end{tikzpicture}
\]

$\bullet P_{\{\alpha_2, \alpha_3\}}$: $\mr{M}_{P_{\{\alpha_2, \alpha_3\}}} = \mr{GL}_2 \times \mr{GL}_1$.
\[
\begin{tikzpicture}
  \node[dotnode, label=below:$\alpha_1$] (1) at (0,0) {};
  \node[crossnode, label=below:$\alpha_2$] (2) at (1.2,0) {};
  \node[crossnode, label=below:$\alpha_3$] (3) at (2.4,0) {};
\end{tikzpicture}
\]

\subsubsection*{Rank 3 ($|I|=3$)}

$\bullet P_{\pi}$: $\mr{M}_{P_{\pi}} = \mr{GL}_1 \times \mr{GL}_1 \times \mr{GL}_1$.
\[
\begin{tikzpicture}
  \node[crossnode, label=below:$\alpha_1$] (1) at (0,0) {};
  \node[crossnode, label=below:$\alpha_2$] (2) at (1.2,0) {};
  \node[crossnode, label=below:$\alpha_3$] (3) at (2.4,0) {};
\end{tikzpicture}
\]

\section{Weyl group of type $\mr{C}_3$} \label{app:weyl_table}

In this appendix, 
we list the elements of the Weyl group $\mathcal{W}$ of type $\mr{C}_3$. 
The following table provides the length $l(w)$ and the images of simple roots under $w^{-1}$. 
For convenience, 
we denote $k = \alpha_1 + \alpha_2$, 
$f = \alpha_2 + \alpha_3$, 
$g = \alpha_1 + \alpha_2 + \alpha_3$, 
$h = 2\alpha_2 + \alpha_3$, 
$i = \alpha_1 + 2\alpha_2 + \alpha_3$, 
and $j = 2\alpha_1 + 2\alpha_2 + \alpha_3$.

  \begin{table}[htbp]
    \centering
    \caption{Weyl group elements $\mr{C}_3$ and $w^{-1}(\al_i)$.}
    \label{table:3}
    \begin{tabular}{cccccc}
    \toprule 
    $w$      & $w^{-1}$  & $\ell(w)$ & $w^{-1}(\alpha_1)$ & $w^{-1}(\alpha_2)$ & $w^{-1}(\alpha_3)$ \\
    \midrule
    \midrule
    $e$       & $e$       & $0$ & $\al_1$  & $\al_2$  & $\al_3$ \\ 
    \midrule
    $s_{1}$   & $s_1$     & $1$ & $-\al_1$ & $k$      & $\al_3$ \\ 
    $s_{2}$   & $s_2$     & $1$ & $k$      & $-\al_2$ & $h$ \\ 
    $s_{3}$   & $s_3$     & $1$ & $\al_1$  & $f$      & $-\al_3$ \\
    \midrule
    $s_{12}$  & $s_{21}$  & $2$ & $-k$     & $\al_1$  & $h$ \\ 
    $s_{13}$  & $s_{13}$  & $2$ & $-\al_1$ & $g$      & $-\al_3$ \\
    $s_{21}$  & $s_{12}$  & $2$ & $\al_2$  & $-k$     & $j$ \\ 
    $s_{23}$  & $s_{32}$  & $2$ & $g$      & $-f$     & $h$ \\
    $s_{32}$  & $s_{23}$  & $2$ & $k$      & $f$      & $-h$ \\
    \midrule
    $s_{121}$ & $s_{121}$ & $3$ & $-\al_2$ & $-\al_1$ & $j$ \\
    $s_{123}$ & $s_{321}$ & $3$ & $-g$     & $\al_1$  & $h$ \\
    $s_{132}$ & $s_{213}$ & $3$ & $-k$     & $i$      & $-h$ \\
    $s_{213}$ & $s_{132}$ & $3$ & $f$      & $-g$     & $j$ \\
    $s_{232}$ & $s_{232}$ & $3$ & $i$      & $-f$     & $\al_3$ \\
    $s_{321}$ & $s_{123}$ & $3$ & $\al_2$  & $g$      & $-j$ \\
    $s_{323}$ & $s_{323}$ & $3$ & $g$      & $\al_2$  & $-h$ \\
    \midrule
    $s_{1213}$ & $s_{1321}$ & $4$ & $-f$     & $-\al_1$ & $j$ \\
    $s_{1232}$ & $s_{2321}$ & $4$ & $-i$     & $k$      & $\al_3$ \\
    $s_{1321}$ & $s_{1213}$ & $4$ & $-\al_2$ & $i$      & $-j$ \\
    $s_{1323}$ & $s_{3213}$ & $4$ & $-g$     & $i$      & $-h$\\
    $s_{2132}$ & $s_{2132}$ & $4$ & $f$      & $-i$     & $j$\\
    $s_{2321}$ & $s_{1232}$ & $4$ & $i$      & $-g$     & $\al_3$\\
    $s_{2323}$ & $s_{2323}$ & $4$ & $i$      & $-\al_2$ & $-\al_3$\\
    $s_{3213}$ & $s_{1323}$ & $4$ & $f$      & $k$      & $-j$\\
    \midrule
    $s_{12132}$ & $s_{21321}$ & $5$ & $-f$     & $-k$     & $j$\\
    $s_{12321}$ & $s_{12321}$ & $5$ & $-i$     & $\al_2$  & $\al_3$\\
    $s_{12323}$ & $s_{23213}$ & $5$ & $-i$     & $g$      & $-\al_3$\\
    $s_{13213}$ & $s_{13213}$ & $5$ & $-f$     & $i$      & $-j$\\
    $s_{21321}$ & $s_{12132}$ & $5$ & $g$      & $-i$     & $h$\\
    $s_{21323}$ & $s_{32132}$ & $5$ & $\al_2$  & $-i$     & $j$\\
    $s_{23213}$ & $s_{12323}$ & $5$ & $i$      & $-k$     & $-\al_3$\\
    $s_{32132}$ & $s_{21323}$ & $5$ & $f$      & $\al_1$   & $-j$\\
    \midrule
    $s_{121321}$ & $s_{121321}$ & $6$ & $-g$     & $-\al_2$ & $h$\\
    $s_{121323}$ & $s_{232132}$ & $6$ & $-\al_2$ & $-g$     & $j$\\
    $s_{123213}$ & $s_{123213}$ & $6$ & $-i$     & $f$      & $-\al_3$\\
    $s_{132132}$ & $s_{213213}$ & $6$ & $-f$     & $g$      & $-j$\\
    $s_{213213}$ & $s_{132132}$ & $6$ & $k$      & $-i$     & $h$\\
    $s_{232132}$ & $s_{121323}$ & $6$ & $g$      & $-\al_1$ & $-h$\\
    $s_{321323}$ & $s_{321323}$ & $6$ & $\al_2$  & $\al_1$  & $-j$\\
    \midrule
    \end{tabular}
   \end{table}
   \newpage
   \begin{table}
   \begin{tabular}{cccccc}
    \toprule 
    $w$      & $w^{-1}$  & $\ell(w)$ & $w^{-1}(\alpha_1)$ & $w^{-1}(\alpha_2)$ & $w^{-1}(\alpha_3)$ \\
    \midrule
    \midrule
   $s_{1213213}$ & $s_{1232132}$ & $7$ & $-k$     & $-f$     & $h$\\
    $s_{1232132}$ & $s_{1213213}$ & $7$ & $-g$     & $f$      & $-h$\\
    $s_{1321323}$ & $s_{2321323}$ & $7$ & $-\al_2$ & $k$      & $-j$\\
    $s_{2132132}$ & $s_{2132132}$ & $7$ & $\al_1$  & $-g$     & $\al_3$\\
    $s_{2321323}$ & $s_{1321323}$ & $7$ & $k$      & $-\al_1$ & $-h$\\
    \midrule
    $s_{12132132}$  & $s_{12132132}$  & $8$ & $-\al_1$ & $-f$    & $\al_3$ \\
    $s_{12321323}$  & $s_{12321323}$  & $8$ & $-k$     & $\al_2$ & $-h$ \\
    $s_{21321323}$  & $s_{21321323}$  & $8$ & $\al_1$  & $-k$    & $-\al_3$ \\
    \midrule
    $s_{121321323}$ & $s_{121321323}$ & $9$ & $-\al_1$ & $-\al_2$ & $-\al_3$ \\
    \bottomrule
    \end{tabular}
    \end{table}

 \section{Weight Coefficients for $w \cdot \lambda$} \label{app:weight_tables}

    This appendix provides the explicit coefficients of the twisted weights $w \cdot \lambda$ in terms of the fundamental dominant weights $\gamma_i^I$ of each Levi quotient. We express the highest weight as $\lambda = n_1 \gamma_1 + n_2 \gamma_2 + n_3 \gamma_3$.
    
    \subsection{General coefficients}
    
    The following tables list the coefficients for general $n_1, n_2, n_3$. These formulas provide the foundation for computing the $E_1$-terms for any irreducible representation $\mathcal{M}_\lambda$.
    
\subsubsection*{Rank $1$($\ |I|=1\ $)}
\ 

$\bullet P_{\{\al_1\}}$:
$\mr{M}_{P_{\{\al_1\}}} \cong \mr{GL}_1 \times \mr{Sp}_4$\\
\textbf{Basis}: $\{\gamma_1^{\{\al_1\}} = \varepsilon_1, \  \gamma_2^{\{\al_1\}} = \varepsilon_2, \ \gamma_3^{\{\al_1\}} = \varepsilon_2+\varepsilon_3\}$

\begin{longtable}{llll}
\toprule
$w$ & \textbf{Coeff for $\gamma_1^{\{\al_1\}}$} & \textbf{Coeff for $\gamma_2^{\{\al_1\}}$} & \textbf{Coeff for $\gamma_3^{\{\al_1\}}$} \\
\midrule
\endfirsthead
\endfoot
$e$         & $n_{1} + n_{2} + n_{3}$       & $n_{2}$             & $n_{3}$ \\
$s_1$       & $n_{2} + n_{3} - 1$           & $n_{1} + n_{2} + 1$ & $n_{3}$ \\
$s_{12}$    & $n_{3} - 2$                   & $n_{1}$             & $n_{2} + n_{3} + 1$ \\
$s_{123}$   & $- n_{3} - 4$                 & $n_{1}$             & $n_{2} + n_{3} + 1$ \\
$s_{1232}$  & $- n_{2} - n_{3} - 5$         & $n_{1} + n_{2} + 1$ & $n_{3}$ \\
$s_{12321}$ & $- n_{1} - n_{2} - n_{3} - 6$ & $n_{2}$             & $n_{3}$ \\
\bottomrule
\end{longtable}

$\bullet P_{\{\alpha_2\}}$:
$\mr{M}_{P_{\{\alpha_2\}}} = \mr{SL}_2 \times \mr{GL}_1 \times \mr{Sp}_2 = \mr{GL}_2 \times \mr{Sp}_2$.\\
\textbf{Basis}: $\{\gamma_1^{\{\al_2\}} = \frac{1}{2}(\varepsilon_1 - \varepsilon_2), \ \gamma_2^{\{\al_2\}} = \varepsilon_1 + \varepsilon_2, \ \gamma_3^{\{\al_2\}} = \varepsilon_3\}$

\begin{longtable}{llll}
\toprule
$w$ & \textbf{Coeff for $\gamma_1^{\{\al_2\}}$} & \textbf{Coeff for $\gamma_2^{\{\al_2\}}$} & \textbf{Coeff for $\gamma_3^{\{\al_2\}}$} \\
\midrule
\endfirsthead
$e$           & $n_1$                   & $\frac{n_1}{2} + n_2 + n_3$                          & $n_3$ \\
$s_2$         & $n_1 + n_2 + 1$         & $\frac{n_1}{2} + \frac{n_2}{2} + n_3 - \frac{1}{2}$  & $n_2 + n_3 + 1$ \\
$s_{21}$      & $n_2$                   & $\frac{n_2}{2} + n_3 - 1$                            & $n_1 + n_2 + n_3 + 2$ \\
$s_{23}$      & $n_1 + n_2 + 2n_3 + 3$  & $\frac{n_1}{2} + \frac{n_2}{2} - \frac{3}{2}$        & $n_2 + n_3 + 1$ \\
$s_{213}$     & $n_2 + 2n_3 + 2$        & $\frac{n_2}{2} - 2$                                  & $n_1 + n_2 + n_3 + 2$ \\
$s_{232}$     & $n_1 + 2n_2 + 2n_3 + 4$ & $\frac{n_1}{2} - 2$                                  & $n_3$ \\
$s_{2132}$    & $n_2 + 2n_3 + 2$        & $-\frac{n_2}{2} - 3$                                 & $n_1 + n_2 + n_3 + 2$ \\
$s_{2321}$    & $n_1 + 2n_2 + 2n_3 + 4$ & $-\frac{n_1}{2} - 3$                                 & $n_3$ \\
$s_{21321}$   & $n_1 + n_2 + 2n_3 + 3$  & $-\frac{n_1}{2} - \frac{n_2}{2} - \frac{7}{2}$       & $n_2 + n_3 + 1$ \\
$s_{21323}$   & $n_2$                   & $-\frac{n_2}{2} - n_3 - 4$                           & $n_1 + n_2 + n_3 + 2$ \\
$s_{213213}$  & $n_1 + n_2 + 1$         & $-\frac{n_1}{2} - \frac{n_2}{2} - n_3 - \frac{9}{2}$ & $n_2 + n_3 + 1$ \\
$s_{2132132}$ & $n_1$                   & $-\frac{n_1}{2} - n_2 - n_3 - 5$                     & $n_3$ \\
\bottomrule
\end{longtable}

$\bullet P_{\{\alpha_3\}}$:
$\mr{M}_{P_{\{\alpha_3\}}} = \mr{SL}_3 \times \mr{GL}_1 = \mr{GL}_3$.\\
\textbf{Basis}: $\{\gamma_1^{\{\al_3\}} = \varepsilon_1, \ \gamma_2^{\{\al_3\}} = \varepsilon_1 + \varepsilon_2, \ \gamma_3^{\{\al_3\}} = \varepsilon_1 + \varepsilon_2 + \varepsilon_3\}$

\begin{longtable}{llll}
\toprule
$w$ & \textbf{Coeff for $\gamma_1^{\{\al_3\}}$} & \textbf{Coeff for $\gamma_2^{\{\al_3\}}$} & \textbf{Coeff for $\gamma_3^{\{\al_3\}}$} \\
\midrule
\endfirsthead
\endfoot
$e$            & $n_1$                   & $n_2$                       & $n_3$ \\
$s_3$          & $n_1$                   & $n_2 + 2n_3 + 2$            & $-n_3 - 2$ \\
$s_{32}$       & $n_1 + n_2 + 1$         & $n_2 + 2n_3 + 2$            & $-n_2 - n_3 - 3$ \\
$s_{321}$      & $n_2$                   & $n_1 + n_2 + 2n_3 + 3$      & $-n_1 - n_2 - n_3 - 4$ \\
$s_{323}$      & $n_1 + n_2 + 2n_3 + 3$  & $n_2$                       & $-n_2 - n_3 - 3$ \\
$s_{3213}$     & $n_2 + 2n_3 + 2$        & $n_1 + n_2 + 1$             & $-n_1 - n_2 - n_3 - 4$ \\
$s_{32132}$    & $n_2 + 2n_3 + 2$        & $n_1$                       & $-n_1 - n_2 - n_3 - 4$ \\
$s_{321323}$   & $n_2$                   & $n_1$                       & $-n_1 - n_2 - n_3 - 4$ \\
\bottomrule
\end{longtable}

\subsubsection*{Rank $2$($\ |I|=2$\ )}
\ 

$\bullet P_{\{\alpha_1, \alpha_2\}}$:
$\mr{M}_{P_{\{\alpha_1, \alpha_2\}}} = \mr{GL}_1 \times \mr{GL}_1 \times \mr{Sp}_2$.\\
\textbf{Basis}: $\{\gamma_1^{\{\al_1, \al_2\}} = \varepsilon_1, \ \gamma_2^{\{\al_1, \al_2\}} = \varepsilon_2, \ \gamma_3^{\{\al_1, \al_2\}} = \varepsilon_3\}$

\begin{longtable}{llll}
\toprule
\textbf{Kostant Rep ($w$)} & \textbf{Coeff for $\gamma_1^{\{\al_1, \al_2\}}$} & \textbf{Coeff for $\gamma_2^{\{\al_1, \al_2\}}$} & \textbf{Coeff for $\gamma_3^{\{\al_1, \al_2\}}$} \\
\midrule
\endfirsthead
\endfoot
$e$            & $n_1 + n_2 + n_3$        & $n_2 + n_3$                 & $n_3$ \\
$s_1$          & $n_2 + n_3 - 1$          & $n_1 + n_2 + n_3 + 1$       & $n_3$ \\
$s_2$          & $n_1 + n_2 + n_3$        & $n_3 - 1$                   & $n_2 + n_3 + 1$ \\
$s_{12}$       & $n_3 - 2$                & $n_1 + n_2 + n_3 + 1$       & $n_2 + n_3 + 1$ \\
$s_{21}$       & $n_2 + n_3 - 1$          & $n_3 - 1$                   & $n_1 + n_2 + n_3 + 2$ \\
$s_{23}$       & $n_1 + n_2 + n_3$        & $-n_3 - 3$                  & $n_2 + n_3 + 1$ \\
$s_{121}$      & $n_3 - 2$                & $n_2 + n_3$                 & $n_1 + n_2 + n_3 + 2$ \\
$s_{123}$      & $-n_3 - 4$               & $n_1 + n_2 + n_3 + 1$       & $n_2 + n_3 + 1$ \\
$s_{213}$      & $n_2 + n_3 - 1$          & $-n_3 - 3$                  & $n_1 + n_2 + n_3 + 2$ \\
$s_{232}$      & $n_1 + n_2 + n_3$        & $-n_2 - n_3 - 4$            & $n_3$ \\
$s_{1213}$     & $-n_3 - 4$               & $n_2 + n_3$                 & $n_1 + n_2 + n_3 + 2$ \\
$s_{1232}$     & $-n_2 - n_3 - 5$         & $n_1 + n_2 + n_3 + 1$       & $n_3$ \\
$s_{2132}$     & $n_3 - 2$                & $-n_2 - n_3 - 4$            & $n_1 + n_2 + n_3 + 2$ \\
$s_{2321}$     & $n_2 + n_3 - 1$          & $-n_1 - n_2 - n_3 - 5$      & $n_3$ \\
$s_{12132}$    & $-n_2 - n_3 - 5$         & $n_3 - 1$                   & $n_1 + n_2 + n_3 + 2$ \\
$s_{12321}$    & $-n_1 - n_2 - n_3 - 6$   & $n_2 + n_3$                 & $n_3$ \\
$s_{21321}$    & $n_3 - 2$                & $-n_1 - n_2 - n_3 - 5$      & $n_2 + n_3 + 1$ \\
$s_{21323}$    & $-n_3 - 4$               & $-n_2 - n_3 - 4$            & $n_1 + n_2 + n_3 + 2$ \\
$s_{121321}$   & $-n_1 - n_2 - n_3 - 6$   & $n_3 - 1$                   & $n_2 + n_3 + 1$ \\
$s_{121323}$   & $-n_2 - n_3 - 5$         & $-n_3 - 3$                  & $n_1 + n_2 + n_3 + 2$ \\
$s_{213213}$   & $-n_3 - 4$               & $-n_1 - n_2 - n_3 - 5$      & $n_2 + n_3 + 1$ \\
$s_{1213213}$  & $-n_1 - n_2 - n_3 - 6$   & $-n_3 - 3$                  & $n_2 + n_3 + 1$ \\
$s_{2132132}$  & $-n_2 - n_3 - 5$         & $-n_1 - n_2 - n_3 - 5$      & $n_3$ \\
$s_{12132132}$ & $-n_1 - n_2 - n_3 - 6$   & $-n_2 - n_3 - 4$            & $n_3$ \\
\bottomrule
\end{longtable}

$\bullet P_{\{\alpha_1, \alpha_3\}}$:
$\mr{M}_{P_{\{\alpha_1, \alpha_3\}}} = \mr{GL}_1 \times \mr{SL}_2 \times \mr{GL}_1 = \mr{GL}_1 \times \mr{GL}_2 $. \\
\textbf{Basis}: $\{\gamma_1^{\{\al_1, \al_3\}} = \varepsilon_1, \ \gamma_2^{\{\al_1, \al_3\}} = \frac{1}{2}(\varepsilon_2-\varepsilon_3), \ \gamma_3^{\{\al_1, \al_3\}} = \varepsilon_2 + \varepsilon_3\}$

\begin{longtable}{llll}
\toprule
\textbf{Kostant Rep ($w$)} & \textbf{Coeff for $\gamma_1^{\{\al_1, \al_3\}}$} & \textbf{Coeff for $\gamma_2^{\{\al_1, \al_3\}}$} & \textbf{Coeff for $\gamma_3^{\{\al_1, \al_3\}}$} \\
\midrule
\endfirsthead
\endfoot
$e$                & $n_1 + n_2 + n_3$                 & $n_2$                                                & $\frac{n_2}{2} + n_3$ \\
$s_1$              & $n_2 + n_3 - 1$                   & $n_1 + n_2 + 1$                                      & $\frac{n_1}{2} + \frac{n_2}{2} + n_3 + \frac{1}{2}$ \\
$s_3$              & $n_1 + n_2 + n_3$                 & $n_2 + 2n_3 + 2$                                     & $\frac{n_2}{2} - 1$ \\
$s_{12}$           & $n_3 - 2$                         & $n_1$                                                & $\frac{n_1}{2} + n_2 + n_3 + 1$ \\
$s_{13}$           & $n_2 + n_3 - 1$                   & $n_1 + n_2 + 2n_3 + 3$                               & $\frac{n_1}{2} + \frac{n_2}{2} - \frac{1}{2}$ \\
$s_{32}$           & $n_1 + n_2 + n_3$                 & $n_2 + 2n_3 + 2$                                     & $-\frac{n_2}{2} - 2$ \\
$s_{123}$          & $-n_3 - 4$                        & $n_1$                                                & $\frac{n_1}{2} + n_2 + n_3 + 1$ \\
$s_{132}$          & $n_3 - 2$                         & $n_1 + 2n_2 + 2n_3 + 4$                              & $\frac{n_1}{2} - 1$ \\
$s_{321}$          & $n_2 + n_3 - 1$                   & $n_1 + n_2 + 2n_3 + 3$                               & $-\frac{n_1}{2} - \frac{n_2}{2} - \frac{5}{2}$ \\
$s_{323}$          & $n_1 + n_2 + n_3$                 & $n_2$                                                & $-\frac{n_2}{2} - n_3 - 3$ \\
$s_{1232}$         & $-n_2 - n_3 - 5$                  & $n_1 + n_2 + 1$                                      & $\frac{n_1}{2} + \frac{n_2}{2} + n_3 + \frac{1}{2}$ \\
$s_{1321}$         & $n_3 - 2$                         & $n_1 + 2n_2 + 2n_3 + 4$                              & $-\frac{n_1}{2} - 2$ \\
$s_{1323}$         & $-n_3 - 4$                        & $n_1 + 2n_2 + 2n_3 + 4$                              & $\frac{n_1}{2} - 1$ \\
$s_{3213}$         & $n_2 + n_3 - 1$                   & $n_1 + n_2 + 1$                                      & $-\frac{n_1}{2} - \frac{n_2}{2} - n_3 - \frac{7}{2}$ \\
$s_{12321}$        & $-n_1 - n_2 - n_3 - 6$            & $n_2$                                                & $\frac{n_2}{2} + n_3$ \\
$s_{12323}$        & $-n_2 - n_3 - 5$                  & $n_1 + n_2 + 2n_3 + 3$                               & $\frac{n_1}{2} + \frac{n_2}{2} - \frac{1}{2}$ \\
$s_{13213}$        & $-n_3 - 4$                        & $n_1 + 2n_2 + 2n_3 + 4$                              & $-\frac{n_1}{2} - 2$ \\
$s_{32132}$        & $n_3 - 2$                         & $n_1$                                                & $-\frac{n_1}{2} - n_2 - n_3 - 4$ \\
$s_{123213}$       & $-n_1 - n_2 - n_3 - 6$            & $n_2 + 2n_3 + 2$                                     & $\frac{n_2}{2} - 1$ \\
$s_{132132}$       & $-n_2 - n_3 - 5$                  & $n_1 + n_2 + 2n_3 + 3$                               & $-\frac{n_1}{2} - \frac{n_2}{2} - \frac{5}{2}$ \\
$s_{321323}$       & $-n_3 - 4$                        & $n_1$                                                & $-\frac{n_1}{2} - n_2 - n_3 - 4$ \\
$s_{1232132}$      & $-n_1 - n_2 - n_3 - 6$            & $n_2 + 2n_3 + 2$                                     & $-\frac{n_2}{2} - 2$ \\
$s_{1321323}$      & $-n_2 - n_3 - 5$                  & $n_1 + n_2 + 1$                                      & $-\frac{n_1}{2} - \frac{n_2}{2} - n_3 - \frac{7}{2}$ \\
$s_{12321323}$     & $-n_1 - n_2 - n_3 - 6$            & $n_2$                                                & $-\frac{n_2}{2} - n_3 - 3$ \\
\bottomrule
\end{longtable}

$\bullet P_{\{\alpha_2, \alpha_3\}}$:
$\mr{M}_{P_{\{\alpha_2, \alpha_3\}}} = \mr{SL}_2 \times \mr{GL}_1 \times \mr{GL}_1 = \mr{GL}_2 \times \mr{GL}_1$.\\
\textbf{Basis}: $\{\gamma_1^{\{\al_2, \al_3\}} = \frac{1}{2}(\varepsilon_1-\varepsilon_2), \ \gamma_2^{\{\al_2, \al_3\}} = \varepsilon_1 + \varepsilon_2, \ \gamma_3^{\{\al_2, \al_3\}} = \varepsilon_3\}$

\begin{longtable}{llll}
\toprule
\textbf{Kostant Rep ($w$)} & \textbf{Coeff for $\gamma_1^{\{\al_2, \al_3\}}$} & \textbf{Coeff for $\gamma_2^{\{\al_2, \al_3\}}$} & \textbf{Coeff for $\gamma_3^{\{\al_2, \al_3\}}$} \\
\midrule
\endfirsthead
\endfoot
$e$            & $n_1$                   & $\frac{n_1}{2} + n_2 + n_3$                          & $n_3$ \\
$s_2$          & $n_1 + n_2 + 1$         & $\frac{n_1}{2} + \frac{n_2}{2} + n_3 - \frac{1}{2}$  & $n_2 + n_3 + 1$ \\
$s_3$          & $n_1$                   & $\frac{n_1}{2} + n_2 + n_3$                          & $-n_3 - 2$ \\
$s_{21}$       & $n_2$                   & $\frac{n_2}{2} + n_3 - 1$                            & $n_1 + n_2 + n_3 + 2$ \\
$s_{23}$       & $n_1 + n_2 + 2n_3 + 3$  & $\frac{n_1}{2} + \frac{n_2}{2} - \frac{3}{2}$        & $n_2 + n_3 + 1$ \\
$s_{32}$       & $n_1 + n_2 + 1$         & $\frac{n_1}{2} + \frac{n_2}{2} + n_3 - \frac{1}{2}$  & $-n_2 - n_3 - 3$ \\
$s_{213}$      & $n_2 + 2n_3 + 2$        & $\frac{n_2}{2} - 2$                                  & $n_1 + n_2 + n_3 + 2$ \\
$s_{232}$      & $n_1 + 2n_2 + 2n_3 + 4$ & $\frac{n_1}{2} - 2$                                  & $n_3$ \\
$s_{321}$      & $n_2$                   & $\frac{n_2}{2} + n_3 - 1$                            & $-n_1 - n_2 - n_3 - 4$ \\
$s_{323}$      & $n_1 + n_2 + 2n_3 + 3$  & $\frac{n_1}{2} + \frac{n_2}{2} - \frac{3}{2}$        & $-n_2 - n_3 - 3$ \\
$s_{2132}$     & $n_2 + 2n_3 + 2$        & $-\frac{n_2}{2} - 3$                                 & $n_1 + n_2 + n_3 + 2$ \\
$s_{2321}$     & $n_1 + 2n_2 + 2n_3 + 4$ & $-\frac{n_1}{2} - 3$                                 & $n_3$ \\
$s_{2323}$     & $n_1 + 2n_2 + 2n_3 + 4$ & $\frac{n_1}{2} - 2$                                  & $-n_3 - 2$ \\
$s_{3213}$     & $n_2 + 2n_3 + 2$        & $\frac{n_2}{2} - 2$                                  & $-n_1 - n_2 - n_3 - 4$ \\
$s_{21321}$    & $n_1 + n_2 + 2n_3 + 3$  & $-\frac{n_1}{2} - \frac{n_2}{2} - \frac{7}{2}$       & $n_2 + n_3 + 1$ \\
$s_{21323}$    & $n_2$                   & $-\frac{n_2}{2} - n_3 - 4$                           & $n_1 + n_2 + n_3 + 2$ \\
$s_{23213}$    & $n_1 + 2n_2 + 2n_3 + 4$ & $-\frac{n_1}{2} - 3$                                 & $-n_3 - 2$ \\
$s_{32132}$    & $n_2 + 2n_3 + 2$        & $-\frac{n_2}{2} - 3$                                 & $-n_1 - n_2 - n_3 - 4$ \\
$s_{213213}$   & $n_1 + n_2 + 1$         & $-\frac{n_1}{2} - \frac{n_2}{2} - n_3 - \frac{9}{2}$ & $n_2 + n_3 + 1$ \\
$s_{232132}$   & $n_1 + n_2 + 2n_3 + 3$  & $-\frac{n_1}{2} - \frac{n_2}{2} - \frac{7}{2}$       & $-n_2 - n_3 - 3$ \\
$s_{321323}$   & $n_2$                   & $-\frac{n_2}{2} - n_3 - 4$                           & $-n_1 - n_2 - n_3 - 4$ \\
$s_{2132132}$  & $n_1$                   & $-\frac{n_1}{2} - n_2 - n_3 - 5$                     & $n_3$ \\
$s_{2321323}$  & $n_1 + n_2 + 1$         & $-\frac{n_1}{2} - \frac{n_2}{2} - n_3 - \frac{9}{2}$ & $-n_2 - n_3 - 3$ \\
$s_{21321323}$ & $n_1$                   & $-\frac{n_1}{2} - n_2 - n_3 - 5$                     & $-n_3 - 2$ \\
\bottomrule
\end{longtable}

\subsubsection*{Rank $3$($\ |I|=3\ )$}
\ 

$\bullet P_{\pi}$:
$\mr{M}_{P_{\pi}} = \mr{GL}_1 \times \mr{GL}_1 \times \mr{GL}_1 = \mr{T}$.\\
\textbf{Basis}: $\{\gamma_1^{\pi}=\varepsilon_1, \gamma_2^{\pi}=\varepsilon_2, \gamma_3^{\pi} = \varepsilon_3\}$

\begin{longtable}{llll}
\toprule
\textbf{Weyl Element ($w$)} & \textbf{Coeff for $\gamma_1^I$} & \textbf{Coeff for $\gamma_2^I$} & \textbf{Coeff for $\gamma_3^I$} \\
\midrule
\endfirsthead
\endfoot
$e$             & $n_1 + n_2 + n_3$        & $n_2 + n_3$           & $n_3$ \\
$s_1$           & $n_2 + n_3 - 1$          & $n_1 + n_2 + n_3 + 1$ & $n_3$ \\
$s_2$           & $n_1 + n_2 + n_3$        & $n_3 - 1$             & $n_2 + n_3 + 1$ \\
$s_3$           & $n_1 + n_2 + n_3$        & $n_2 + n_3$           & $-n_3 - 2$ \\
$s_{12}$        & $n_3 - 2$                & $n_1 + n_2 + n_3 + 1$ & $n_2 + n_3 + 1$ \\
$s_{13}$        & $n_2 + n_3 - 1$          & $n_1 + n_2 + n_3 + 1$ & $-n_3 - 2$ \\
$s_{21}$        & $n_2 + n_3 - 1$          & $n_3 - 1$             & $n_1 + n_2 + n_3 + 2$ \\
$s_{23}$        & $n_1 + n_2 + n_3$        & $-n_3 - 3$            & $n_2 + n_3 + 1$ \\
$s_{32}$        & $n_1 + n_2 + n_3$        & $n_3 - 1$             & $-n_2 - n_3 - 3$ \\
$s_{121}$       & $n_3 - 2$                & $n_2 + n_3$           & $n_1 + n_2 + n_3 + 2$ \\
$s_{123}$       & $-n_3 - 4$               & $n_1 + n_2 + n_3 + 1$ & $n_2 + n_3 + 1$ \\
$s_{132}$       & $n_3 - 2$                & $n_1 + n_2 + n_3 + 1$ & $-n_2 - n_3 - 3$ \\
$s_{213}$       & $n_2 + n_3 - 1$          & $-n_3 - 3$            & $n_1 + n_2 + n_3 + 2$ \\
$s_{232}$       & $n_1 + n_2 + n_3$        & $-n_2 - n_3 - 4$      & $n_3$ \\
$s_{321}$       & $n_2 + n_3 - 1$          & $n_3 - 1$             & $-n_1 - n_2 - n_3 - 4$ \\
$s_{323}$       & $n_1 + n_2 + n_3$        & $-n_3 - 3$            & $-n_2 - n_3 - 3$ \\
$s_{1213}$      & $-n_3 - 4$               & $n_2 + n_3$           & $n_1 + n_2 + n_3 + 2$ \\
$s_{1232}$      & $ - n_2 - n_3 - 5$       & $n_1 + n_2 + n_3 + 1$ & $n_3$ \\
$s_{1321}$      & $n_3 - 2$                & $n_2 + n_3$           & $- n_1 - n_2 - n_3 - 4$ \\
$s_{1323}$      & $-n_3 - 4$               & $n_1 + n_2 + n_3 + 1$ & $-n_2 - n_3 - 3$ \\
$s_{2132}$      & $n_3 - 2$                & $-n_2 - n_3 - 4$      & $n_1 + n_2 + n_3 + 2$ \\
$s_{2321}$      & $n_2 + n_3 - 1$          & $-n_1 - n_2 - n_3 -5$ & $n_3$ \\
$s_{2323}$      & $n_1 + n_2 + n_3$        & $-n_2 - n_3 - 4$      & $-n_3 - 2$ \\
$s_{3213}$      & $n_2 + n_3 - 1$          & $-n_3 - 3$            & $-n_1 - n_2 - n_3 - 4$ \\
$s_{12132}$     & $-n_2 - n_3 - 5$         & $n_3 - 1$             & $n_1 + n_2 + n_3 + 2$ \\
$s_{12321}$     & $-n_1 - n_2 - n_3 - 6$   & $n_2 + n_3$           & $n_3$ \\
$s_{12323}$     & $-n_2 - n_3 - 5$         & $n_1 + n_2 + n_3 + 1$ & $-n_3 - 2$ \\
$s_{13213}$     & $-n_3 - 4$               & $n_2 + n_3$           & $-n_1 - n_2 - n_3 - 4$ \\
$s_{21321}$     & $n_3 - 2$                & $-n_1 - n_2 - n_3 -5$ & $n_2 + n_3 + 1$ \\
$s_{21323}$     & $-n_3 -4$                & $-n_2 - n_3 -4$       & $n_1 + n_2 + n_3 + 2$ \\
$s_{23213}$     & $n_2 + n_3 -1$           & $-n_1 - n_2 - n_3 -5$ & $-n_3 - 2$ \\
$s_{32132}$     & $n_3 - 2$                & $-n_2 - n_3 - 4$      & $-n_1 - n_2 - n_3 - 4$ \\
$s_{121321}$    & $-n_1 - n_2 - n_3 - 6$   & $n_3 - 1$             & $n_2 + n_3 + 1$ \\
$s_{121323}$    & $-n_2 -n_3 - 5$          & $-n_3 - 3$            & $n_1 + n_2 + n_3 + 2$ \\
$s_{123213}$    & $-n_1 - n_2 - n_3 - 6$   & $n_2 + n_3$           & $-n_3 - 2$ \\
$s_{132132}$    & $-n_2 - n_3 - 5$         & $n_3 - 1$             & $-n_1 - n_2 - n_3 - 4$ \\
$s_{213213}$    & $-n_3 - 4$               & $-n_1 - n_2 - n_3 -5$ & $n_2 + n_3 + 1$ \\
$s_{232132}$    & $n_3 -2$                 & $-n_1 - n_2 - n_3 -5$ & $-n_2 - n_3 - 3$ \\
$s_{321323}$    & $-n_3 - 4$               & $-n_2 - n_3 - 4$      & $-n_1 - n_2 - n_3 - 4$ \\
$s_{1213213}$   & $-n_1 - n_2 - n_3 - 6$   & $-n_3 - 3$            & $n_2 + n_3 + 1$ \\
$s_{1232132}$   & $-n_1 - n_2 - n_3 - 6$   & $n_3 - 1$             & $-n_2 - n_3 - 3$ \\
$s_{1321323}$   & $-n_2 - n_3 - 5$         & $-n_3 - 3$            & $-n_1 - n_2 - n_3 - 4$ \\
$s_{2132132}$   & $-n_2 - n_3 -5$          & $-n_1 - n_2 - n_3 -5$ & $n_3$ \\
$s_{2321323}$   & $-n_3 - 4$               & $-n_1 - n_2 - n_3 -5$ & $-n_2 - n_3 - 3$ \\
$s_{12132132}$  & $-n_1 - n_2 - n_3 - 6$   & $-n_2 - n_3 - 4$      & $n_3$ \\
$s_{12321323}$  & $-n_1 - n_2 - n_3 - 6$   & $- n_3 - 3$           & $-n_2 - n_3 - 3$ \\
$s_{21321323}$  & $-n_2 - n_3 - 5$         & $-n_1 - n_2 - n_3 -5$ & $-n_3 - 2$ \\
$s_{121321323}$ & $-n_1 - n_2 - n_3 - 6$   & $-n_2 - n_3 - 4$      & $-n_3 - 2$ \\
\bottomrule
\end{longtable}

\subsection{Specialization to the trivial representation}
    
  In the specific case of the trivial representation where $n_1 = n_2 = n_3 = 0$, 
  the coefficients simplify to the following values. 
  These constants are used to evaluate the parity conditions in Section 3 and to determine the dimensions of the cohomology groups in Section 4.
    
  \subsubsection*{Rank $1$($\ |I|=1\ $)}
  \ 
  
  
  $\bullet P_{\{\al_1\}}$:
  $\mr{M}_{P_{\{\al_1\}}} \cong \mr{GL}_1 \times \mr{Sp}_4$\\
  \textbf{Basis}: $\{\gamma_1^{\{\al_1\}} = \varepsilon_1, \  \gamma_2^{\{\al_1\}} = \varepsilon_2, \ \gamma_3^{\{\al_1\}} = \varepsilon_2+\varepsilon_3\}$
  
  \begin{longtable}{llll}
  \toprule
  $w$ & \textbf{Coeff for $\gamma_1^{\{\al_1\}}$} & \textbf{Coeff for $\gamma_2^{\{\al_1\}}$} & \textbf{Coeff for $\gamma_3^{\{\al_1\}}$} \\
  \midrule
  \endfirsthead
  \endfoot
  $e$         & $0$  & $0$ & $0$ \\
  $s_1$       & $-1$ & $1$ & $0$ \\
  $s_{12}$    & $-2$ & $0$ & $1$ \\
  $s_{123}$   & $-4$ & $0$ & $1$ \\
  $s_{1232}$  & $-5$ & $1$ & $0$ \\
  $s_{12321}$ & $-6$ & $0$ & $0$ \\
  \bottomrule
  \end{longtable}
  
  $\bullet P_{\{\alpha_2\}}$:
  $\mr{M}_{P_{\{\alpha_2\}}} = \mr{SL}_2 \times \mr{GL}_1 \times \mr{Sp}_2 = \mr{GL}_2 \times \mr{Sp}_2$.\\
  \textbf{Basis}: $\{\gamma_1^{\{\al_2\}} = \frac{1}{2}(\varepsilon_1 - \varepsilon_2), \ \gamma_2^{\{\al_2\}} = \varepsilon_1 + \varepsilon_2, \ \gamma_3^{\{\al_2\}} = \varepsilon_3\}$
  
  \begin{longtable}{llll}
  \toprule
  $w$ & \textbf{Coeff for $\gamma_1^{\{\al_2\}}$} & \textbf{Coeff for $\gamma_2^{\{\al_2\}}$} & \textbf{Coeff for $\gamma_3^{\{\al_2\}}$} \\
  \midrule
  \endfirsthead
  $e$           & $0$ & $0$            & $0$ \\
  $s_2$         & $1$ & $-\frac{1}{2}$ & $1$ \\
  $s_{21}$      & $0$ & $-1$           & $2$ \\
  $s_{23}$      & $3$ & $-\frac{3}{2}$ & $1$ \\
  $s_{213}$     & $2$ & $-2$           & $2$ \\
  $s_{232}$     & $4$ & $-2$           & $0$ \\
  $s_{2132}$    & $2$ & $-3$           & $2$ \\
  $s_{2321}$    & $4$ & $-3$           & $0$ \\
  $s_{21321}$   & $3$ & $-\frac{7}{2}$ & $1$ \\
  $s_{21323}$   & $0$ & $-4$           & $2$ \\
  $s_{213213}$  & $1$ & $-\frac{9}{2}$ & $1$ \\
  $s_{2132132}$ & $0$ & $-5$           & $0$ \\
  \bottomrule
  \end{longtable}
  
  $\bullet P_{\{\alpha_3\}}$:
  $\mr{M}_{P_{\{\alpha_3\}}} = \mr{SL}_3 \times \mr{GL}_1 = \mr{GL}_3$.\\
  \textbf{Basis}: $\{\gamma_1^{\{\al_3\}} = \varepsilon_1, \ \gamma_2^{\{\al_3\}} = \varepsilon_1 + \varepsilon_2, \ \gamma_3^{\{\al_3\}} = \varepsilon_1 + \varepsilon_2 + \varepsilon_3\}$
  
  \begin{longtable}{llll}
  \toprule
  $w$ & \textbf{Coeff for $\gamma_1^{\{\al_3\}}$} & \textbf{Coeff for $\gamma_2^{\{\al_3\}}$} & \textbf{Coeff for $\gamma_3^{\{\al_3\}}$} \\
  \midrule
  \endfirsthead
  \endfoot
  $e$            & $0$ & $0$ & $0$ \\
  $s_3$          & $0$ & $2$ & $-2$ \\
  $s_{32}$       & $1$ & $2$ & $-3$ \\
  $s_{321}$      & $0$ & $3$ & $-4$ \\
  $s_{323}$      & $3$ & $0$ & $-3$ \\
  $s_{3213}$     & $2$ & $1$ & $-4$ \\
  $s_{32132}$    & $2$ & $0$ & $-4$ \\
  $s_{321323}$   & $0$ & $0$ & $-4$ \\
  \bottomrule
  \end{longtable}

  \subsubsection*{Rank $2$($\ |I|=2$\ )}
  \ 
  
  $\bullet P_{\{\alpha_1, \alpha_2\}}$:
  $\mr{M}_{P_{\{\alpha_1, \alpha_2\}}} = \mr{GL}_1 \times \mr{GL}_1 \times \mr{Sp}_2$.\\
  \textbf{Basis}: $\{\gamma_1^{\{\al_1, \al_2\}} = \varepsilon_1, \ \gamma_2^{\{\al_1, \al_2\}} = \varepsilon_2, \ \gamma_3^{\{\al_1, \al_2\}} = \varepsilon_3\}$
  
  \begin{longtable}{llll}
  \toprule
  \textbf{Kostant Rep ($w$)} & \textbf{Coeff for $\gamma_1^{\{\al_1, \al_2\}}$} & \textbf{Coeff for $\gamma_2^{\{\al_1, \al_2\}}$} & \textbf{Coeff for $\gamma_3^{\{\al_1, \al_2\}}$} \\
  \midrule
  \endfirsthead
  \endfoot
  $e$            & $0$  & $0$  & $0$ \\
  $s_1$          & $-1$ & $1$  & $0$ \\
  $s_2$          & $0$  & $-1$ & $1$ \\
  $s_{12}$       & $-2$ & $1$  & $1$ \\
  $s_{21}$       & $-1$ & $-1$ & $2$ \\
  $s_{23}$       & $0$  & $-3$ & $1$ \\
  $s_{121}$      & $-2$ & $0$  & $2$ \\
  $s_{123}$      & $-4$ & $1$  & $1$ \\
  $s_{213}$      & $-1$ & $-3$ & $2$ \\
  $s_{232}$      & $0$  & $-4$ & $0$ \\
  $s_{1213}$     & $-4$ & $0$  & $2$ \\
  $s_{1232}$     & $-5$ & $1$  & $0$ \\
  $s_{2132}$     & $-2$ & $-4$ & $2$ \\
  $s_{2321}$     & $-1$ & $-5$ & $0$ \\
  $s_{12132}$    & $-5$ & $-1$ & $2$ \\
  $s_{12321}$    & $-6$ & $0$  & $0$ \\
  $s_{21321}$    & $-2$ & $-5$ & $1$ \\
  $s_{21323}$    & $-4$ & $-4$ & $2$ \\
  $s_{121321}$   & $-6$ & $-1$ & $1$ \\
  $s_{121323}$   & $-5$ & $-3$ & $2$ \\
  $s_{213213}$   & $-4$ & $-5$ & $1$ \\
  $s_{1213213}$  & $-6$ & $-3$ & $1$ \\
  $s_{2132132}$  & $-5$ & $-5$ & $0$ \\
  $s_{12132132}$ & $-6$ & $-4$ & $0$ \\
  \bottomrule
  \end{longtable}

  $\bullet P_{\{\alpha_1, \alpha_3\}}$:
  $\mr{M}_{P_{\{\alpha_1, \alpha_3\}}} = \mr{GL}_1 \times \mr{SL}_2 \times \mr{GL}_1 = \mr{GL}_1 \times \mr{GL}_2$. \\
  \textbf{Basis}: $\{\gamma_1^{\{\al_1, \al_3\}} = \varepsilon_1, \ \gamma_2^{\{\al_1, \al_3\}} = \frac{1}{2}(\varepsilon_2-\varepsilon_3), \ \gamma_3^{\{\al_1, \al_3\}} = \varepsilon_2 + \varepsilon_3\}$
  
  \begin{longtable}{llll}
  \toprule
  \textbf{Kostant Rep ($w$)} & \textbf{Coeff for $\gamma_1^{\{\al_1, \al_3\}}$} & \textbf{Coeff for $\gamma_2^{\{\al_1, \al_3\}}$} & \textbf{Coeff for $\gamma_3^{\{\al_1, \al_3\}}$} \\
  \midrule
  \endfirsthead
  \endfoot
  $e$            & $0$  & $0$  & $0$ \\
  $s_1$          & $-1$ & $1$  & $\frac{1}{2}$ \\
  $s_3$          & $0$  & $2$  & $-1$ \\
  $s_{12}$       & $-2$ & $0$  & $1$ \\
  $s_{13}$       & $-1$ & $3$  & $-\frac{1}{2}$ \\
  $s_{32}$       & $0$  & $2$  & $-2$ \\
  $s_{123}$      & $-4$ & $0$  & $1$ \\
  $s_{132}$      & $-2$ & $4$  & $-1$ \\
  $s_{321}$      & $-1$ & $3$  & $-\frac{5}{2}$ \\
  $s_{323}$      & $0$  & $0$  & $-3$ \\
  $s_{1232}$     & $-5$ & $1$  & $\frac{1}{2}$ \\
  $s_{1321}$     & $-2$ & $4$  & $-2$ \\
  $s_{1323}$     & $-4$ & $4$  & $-1$ \\
  $s_{3213}$     & $-1$ & $1$  & $-\frac{7}{2}$ \\
  $s_{12321}$    & $-6$ & $0$  & $0$ \\
  $s_{12323}$    & $-5$ & $3$  & $-\frac{1}{2}$ \\
  $s_{13213}$    & $-4$ & $4$  & $-2$ \\
  $s_{32132}$    & $-2$ & $0$  & $-4$ \\
  $s_{123213}$   & $-6$ & $2$  & $-1$ \\
  $s_{132132}$   & $-5$ & $3$  & $-\frac{5}{2}$ \\
  $s_{321323}$   & $-4$ & $0$  & $-4$ \\
  $s_{1232132}$  & $-6$ & $2$  & $-2$ \\
  $s_{1321323}$  & $-5$ & $1$  & $-\frac{7}{2}$ \\
  $s_{12321323}$ & $-6$ & $0$  & $-3$ \\
  \bottomrule
  \end{longtable}

  $\bullet P_{\{\alpha_2, \alpha_3\}}$:
  $\mr{M}_{P_{\{\alpha_2, \alpha_3\}}} = \mr{SL}_2 \times \mr{GL}_1 \times \mr{GL}_1 = \mr{GL}_2 \times \mr{GL}_1$.\\
  \textbf{Basis}: $\{\gamma_1^{\{\al_2, \al_3\}} = \frac{1}{2}(\varepsilon_1-\varepsilon_2), \ \gamma_2^{\{\al_2, \al_3\}} = \varepsilon_1 + \varepsilon_2, \ \gamma_3^{\{\al_2, \al_3\}} = \varepsilon_3\}$
  
  \begin{longtable}{llll}
  \toprule
  \textbf{Kostant Rep ($w$)} & \textbf{Coeff for $\gamma_1^{\{\al_2, \al_3\}}$} & \textbf{Coeff for $\gamma_2^{\{\al_2, \al_3\}}$} & \textbf{Coeff for $\gamma_3^{\{\al_2, \al_3\}}$} \\
  \midrule
  \endfirsthead
  \endfoot
  $e$            & $0$ & $0$             & $0$ \\
  $s_2$          & $1$ & $-\frac{1}{2}$  & $1$ \\
  $s_3$          & $0$ & $0$             & $-2$ \\
  $s_{21}$       & $0$ & $-1$            & $2$ \\
  $s_{23}$       & $3$ & $-\frac{3}{2}$  & $1$ \\
  $s_{32}$       & $1$ & $-\frac{1}{2}$  & $-3$ \\
  $s_{213}$      & $2$ & $-2$            & $2$ \\
  $s_{232}$      & $4$ & $-2$            & $0$ \\
  $s_{321}$      & $0$ & $-1$            & $-4$ \\
  $s_{323}$      & $3$ & $-\frac{3}{2}$  & $-3$ \\
  $s_{2132}$     & $2$ & $-3$            & $2$ \\
  $s_{2321}$     & $4$ & $-3$            & $0$ \\
  $s_{2323}$     & $4$ & $-2$            & $-2$ \\
  $s_{3213}$     & $2$ & $-2$            & $-4$ \\
  $s_{21321}$    & $3$ & $-\frac{7}{2}$  & $1$ \\
  $s_{21323}$    & $0$ & $-4$            & $2$ \\
  $s_{23213}$    & $4$ & $-3$            & $-2$ \\
  $s_{32132}$    & $2$ & $-3$            & $-4$ \\
  $s_{213213}$   & $1$ & $-\frac{9}{2}$  & $1$ \\
  $s_{232132}$   & $3$ & $-\frac{7}{2}$  & $-3$ \\
  $s_{321323}$   & $0$ & $-4$            & $-4$ \\
  $s_{2132132}$  & $0$ & $-5$            & $0$ \\
  $s_{2321323}$  & $1$ & $-\frac{9}{2}$  & $-3$ \\
  $s_{21321323}$ & $0$ & $-5$            & $-2$ \\
  \bottomrule
  \end{longtable}

  \subsubsection*{Rank $3$($\ |I|=3\ )$}
  \ 
  
  $\bullet P_{\pi}$:
  $\mr{M}_{P_{\pi}} = \mr{GL}_1 \times \mr{GL}_1 \times \mr{GL}_1 = \mr{T}$.\\
  \textbf{Basis}: $\{\gamma_1^{\pi}=\varepsilon_1, \gamma_2^{\pi}=\varepsilon_2, \gamma_3^{\pi} = \varepsilon_3\}$
  
  \begin{longtable}{llll}
  \toprule
  \textbf{Weyl Element ($w$)} & \textbf{Coeff for $\gamma_1$} & \textbf{Coeff for $\gamma_2$} & \textbf{Coeff for $\gamma_3$} \\
  \midrule
  \endfirsthead
  \endfoot
  $e$             & $0$  & $0$  & $0$ \\
  $s_1$           & $-1$ & $1$  & $0$ \\
  $s_2$           & $0$  & $-1$ & $1$ \\
  $s_3$           & $0$  & $0$  & $-2$ \\
  $s_{12}$        & $-2$ & $1$  & $1$ \\
  $s_{13}$        & $-1$ & $1$  & $-2$ \\
  $s_{21}$        & $-1$ & $-1$ & $2$ \\
  $s_{23}$        & $0$  & $-3$ & $1$ \\
  $s_{32}$        & $0$  & $-1$ & $-3$ \\
  $s_{121}$       & $-2$ & $0$  & $2$ \\
  $s_{123}$       & $-4$ & $1$  & $1$ \\
  $s_{132}$       & $-2$ & $1$  & $-3$ \\
  $s_{213}$       & $-1$ & $-3$ & $2$ \\
  $s_{232}$       & $0$  & $-4$ & $0$ \\
  $s_{321}$       & $-1$ & $-1$ & $-4$ \\
  $s_{323}$       & $0$  & $-3$ & $-3$ \\
  $s_{1213}$      & $-4$ & $0$  & $2$ \\
  $s_{1232}$      & $-5$ & $1$  & $0$ \\
  $s_{1321}$      & $-2$ & $0$  & $-4$ \\
  $s_{1323}$      & $-4$ & $1$  & $-3$ \\
  $s_{2132}$      & $-2$ & $-4$ & $2$ \\
  $s_{2321}$      & $-1$ & $-5$ & $0$ \\
  $s_{2323}$      & $0$  & $-4$ & $-2$ \\
  $s_{3213}$      & $-1$ & $-3$ & $-4$ \\
  $s_{12132}$     & $-5$ & $-1$ & $2$ \\
  $s_{12321}$     & $-6$ & $0$  & $0$ \\
  $s_{12323}$     & $-5$ & $1$  & $-2$ \\
  $s_{13213}$     & $-4$ & $0$  & $-4$ \\
  $s_{21321}$     & $-2$ & $-5$ & $1$ \\
  $s_{21323}$     & $-4$ & $-4$ & $2$ \\
  $s_{23213}$     & $-1$ & $-5$ & $-2$ \\
  $s_{32132}$     & $-2$ & $-4$ & $-4$ \\
  $s_{121321}$    & $-6$ & $-1$ & $1$ \\
  $s_{121323}$    & $-5$ & $-3$ & $2$ \\
  $s_{123213}$    & $-6$ & $0$  & $-2$ \\
  $s_{132132}$    & $-5$ & $-1$ & $-4$ \\
  $s_{213213}$    & $-4$ & $-5$ & $1$ \\
  $s_{232132}$    & $-2$ & $-5$ & $-3$ \\
  $s_{321323}$    & $-4$ & $-4$ & $-4$ \\
  $s_{1213213}$   & $-6$ & $-3$ & $1$ \\
  $s_{1232132}$   & $-6$ & $-1$ & $-3$ \\
  $s_{1321323}$   & $-5$ & $-3$ & $-4$ \\
  $s_{2132132}$   & $-5$ & $-5$ & $0$ \\
  $s_{2321323}$   & $-4$ & $-5$ & $-3$ \\
  $s_{12132132}$  & $-6$ & $-4$ & $0$ \\
  $s_{12321323}$  & $-6$ & $-3$ & $-3$ \\
  $s_{21321323}$  & $-5$ & $-5$ & $-2$ \\
  $s_{121321323}$ & $-6$ & $-4$ & $-2$ \\
  \bottomrule
  \end{longtable}

\bigskip
\noindent
\section*{Acknowledgement}
I would like to thank Professor Lin Weng for many helpful discussions and for his continuous support. 
He encouraged me to read \cite{sl3} and compute the $\mr{Sp}_6$ case,
which became the starting point of this work.
I am deeply grateful for his patient guidance throughout the preparation of this manuscript.

I would also like to thank Professor Ivan Horozov for his helpful correspondence.
After reading an earlier version of this paper, 
he suggested checking whether the boundary cohomology satisfies Poincar\'e duality.
His insightful advice led me to identify and correct a calculation error, ensuring the accuracy of the final result.

\end{document}